\documentclass[11pt]{article}

\usepackage{amsmath,amsthm}
\usepackage{mdframed}
\usepackage{amssymb}
\usepackage{amsfonts}
\usepackage{cite}
\usepackage[dvips]{graphicx}
\usepackage{float}

\def \Sinverse  { U_{I,\epsilon}(M^{+}) }

\newtheorem{theorem}{Theorem}[section]

\newtheorem{corollary}[theorem]{Corollary}
\newtheorem{definition}[theorem]{Definition}

\newtheorem{lemma}[theorem]{Lemma}
\newtheorem{prop}[theorem]{Proposition}

\newtheorem{remark}[theorem]{Remark}
\newtheorem{assumption}[theorem]{Assumption}

\def \bC {\mathbb C}

\def \fkz{\mathfrak{z}}
\def \fl {\mathfrak{l}}

\def \inv{^{-1}}

\def \fc {\mathfrak{c}}

\def \upper{^{(i)}}

\def \fa {\mathfrak{a}}
\def \p{\partial}

\def \bR {\mathbb R}
\def \bh {\mathbb H}

 \def \bk {\bf k}

 \def \fe {\mathfrak{e}}

               \def \bR{{\mathbb R}}
               \def \bZ{{\mathbb Z}}

              \def \E{{\mathcal E}}

\def \p{\partial}
\def \E{{\mathcal E}}
\def \mbf {\mathbf }
\def \v {\vskip 0.1in}
\def \n {\noindent}

\begin{document}

\begin{center}
  {\LARGE Relative Invariants, Contact Geometry and Open String Invariants}
  \end{center}

  \noindent
  \begin{center}
   {\large An-Min Li}\footnote{partially supported by a NSFC grant}\\[5pt]
      Department of Mathematics, Sichuan University\\
        Chengdu, PRC\\[5pt]
  {\large  Li Sheng}\footnote{partially supported by a NSFC grant}\\[5pt]
      Department of Mathematics, Sichuan University\\
        Chengdu, PRC\\[5pt]
\end{center}

\begin{abstract}
In this paper we propose a theory of contact invariants and open string invariants, which
are generalizations of the relative invariants. We introduce  two moduli spaces $\overline{\mathcal{M}}_{A}(M^{+},C,g,m+\nu,{\bf y},{\bf p},(\mathbf{k},\mathfrak{e}))$ and $\overline{\mathcal{M}}_{A}(M,L;g,m+\nu,{\bf y},{\bf p},\overrightarrow{\mu})$, prove the compactness of the moduli spaces and the existence of the  invariants.
\end{abstract}

  \tableofcontents
\v\v

\section{Introduction}

Open string invariant theory have been studied by many mathematicians and physicists (see \cite{AKV,AV,FL,KL,LLLZ}). This theory closely relates to relative Gromov-Witten theory and contact geometry. In this paper we propose a theory of contact invariants and open string invariants, which are generalizations of the relative invariants. We outline the idea as follows.
\v
{\bf 1).} Let $(M,\omega)$ be a compact symplectic manifold, $L\subset M$ be a compact Lagrangian sub-manifold.
Let $(x_1,\cdots, x_n)$ be a local coordinate system in $O\subset L$, there is a canonical coordinates
$(x_1,\cdots,x_n,y_1,\cdots,y_n)$ in $T^{*}L|_{O}.$ Suppose that given a Riemannian metric $g$ on $L$, and  consider the unit sphere bundle $\widetilde{M}$. Let $\Lambda$ be the Liouville form on $T^{*}L$, denote $\lambda= -\vartheta \mid_{\widetilde{M}}$. Then $(\widetilde{M}, \lambda)$ is a contact manifold with contact form $\lambda$. The Reeb vector field is given by $X=\sum \frac{y_i}{\|y\|}\frac{\partial}{\partial x_i}$. By Lagrangian Neighborhood Theorem we can write $M-L$ as
$$M^+=M_{0}\bigcup\{[0,\infty)\times
\widetilde{M}\},$$
where $M_0$ is a compact symplectic manifold with boundary.
We choose an almost complex structure $J$ on $M^+$ such that $J$ is tamed by $\omega$ and over the cylinder end $J$ is given by
$$J\mid_{\xi}= \tilde{J},\;\;\;JX=-\frac{\partial}{\partial a},\;\;J(\frac{\partial}{\partial a})=X,$$
where $\xi=ker \lambda$, $a$ is the canonical coordinate in $\mathbb R$,
and $\tilde{J}$ is a $d\lambda$-tame almost complex structure in $\xi$. Assume that the periodic orbit sets of $X$ are either non-degenerate or of Bott-type, and $\tilde{J}$ can be chosen such that
$L_X J=0$ along every periodic orbit.

\v
Let $\stackrel{\circ}{\Sigma}$ be a Riemann surface with a puncture point $p$. We use the cylinder coordinates $(s,t)$ near $p$, i.e.,we consider a neighborhood of $p$ as $(s_0,\infty)\times S^1$.  Let ${u}:\stackrel{\circ}{\Sigma}
\rightarrow M^+$ be a ${J}$-holomorphic map with finite energy. Suppose that
$$[\lim_{s\rightarrow \infty} u(s,S^1)]=\sum \mu^i[c_i],$$
where $[c_i],i=1,...,\fa,$ is a bases in $H_1(L,\mathbb Z)$ and $\mu^{i}\in \mathbb{Z}$. Then $u(s,t)$ converges to a periodic orbit $x\subset \widetilde{M}$ of the Reeb vector field $X$ as $s\rightarrow \infty$. We can view $\lim_{s\rightarrow \infty} u(s,S^1)$
as a loop in $L$, representing $\sum \mu^i[c_i]$. In this way, we can control the behaviour at infinity of ${J}$-holomorphic maps with finite energy.
\v
{\bf 2).} To compactify the moduli space we need study the ${J}$-holomorphic maps into $\{(-\infty,\infty)\times
\widetilde{M}\}$. There are two global vector fields on $\mathbb R\times \widetilde{M}$: $\frac{\p}{\p a}$ and $X$.
Similar to the situation of relative invariants, there is a $\mathbb R$ action, which induces a $\mathbb R$-action on the moduli space of $J$-holomorphic maps. We need mod this action. Since there is a vector field $X$ with $|X|=1$ on $\widetilde{M}$, the Reeb vector field, there is a one parameter group $\varphi_\theta$ action on $\widetilde{M}$ generated by $-X$. In particular, there is a $S^1$-action on every periodic orbit, corresponding to the freedom of the choosing origin of $S^1$. Along every periodic orbit we have $L_X \lambda =0$, on the other hand, we assume that on every periodic orbit $L_X J=0$, then we can mod this action.

\v
On the other hand, we choose the Li-Ruan's compactification in \cite{LR}, that is, we firstly let the Riemann surfaces degenerate in Delingne-Mumford space and then let $M^+$ degenerate compatibly
as in the situation of relative invariants. At any node, the Riemann surface degenerates independently with two parameters, which compatible with those freedoms of choosing the origins ( see section \S\ref{compact_theorem} for degeneration and section \S\ref{gluing_pregluing} for gluing). Then both blowups at interior and at infinity lead boundaries of codimension 2 or more in the moduli space.

\v
{\bf 3).} Another core technical issue in this paper is to define invariants using virtual techniques.
As we know, there had been several different approaches, such as Fukaya-Ono \cite{FO}, Li-Tian \cite{LT}, Liu-Tian(\cite{LiuT}),
Ruan(\cite{R2}), Siebert(\cite{S}) and etc.
In \cite{LR}, Li and Ruan provide a completely new approach to this issue: they show that the invariants can be defined via the integration on the {\em top} stratum virtually. In order to
achieve this goal, they provide {\em refined} estimates of differentiations for gluing parameters: $\partial/\partial r$.  In \cite{LR} the estimates for $\frac{\p}{\p r}$  is of order $r^{-2}$ when $r\to \infty$, that is enough to define invariants. In this paper these estimates achieve to be of exponential decay order $\exp(-\mathfrak{c}r)$. Then we can use the estimates to define the invariants and prove the smoothness of the moduli space. The same method can be applied to GW-invariant for compact symplectic manifold.

\v
In this paper we introduce  two moduli spaces $\overline{\mathcal{M}}_{A}(M^{+},C,g,m+\nu,{\bf y},{\bf p},(\mathbf{k},\mathfrak{e}))$ and $\overline{\mathcal{M}}_{A}(M,L;g,m+\nu,{\bf y},{\bf p},\overrightarrow{\mu})$, prove the compactness of the moduli spaces and the existence of the contact invariants $\Psi^{(C)}_{(A,g,m+\nu,\mathbf{k},\mathfrak{e})}(\alpha_1,...,
\alpha_{m} ; \beta_{m + 1},..., \beta_{m+\nu})$ and the
open string invariants $\Psi^{(L)}_{(A,g,m+\nu,\overrightarrow{\mu})}(\alpha_1,...,
\alpha_{m})$. In our next paper \cite{LS} we will prove the smoothness of the two moduli spaces. Our open string invariants $\Psi^{(L)}_{(A,g,m+\nu,\overrightarrow{\mu})}$ can be generalized to $L$ which is a disjoint union of compact Lagrangian sub-manifolds $L_1$,...,$L_d$.
\v
We consider a neighborhood of Lagrangian  sub-manifold $L$ as $\mathbb R\times \widetilde{M}$. By the same method above we can define a local open string invariant. We will discuss this problem and calculate some examples in our next paper.

\section{Symplectic manifolds with cylindrical ends}

\subsection{Contact manifolds}

Let $(\mathcal{Q},\lambda)$ be a $(2n-1)$-dimensional compact manifold
equipped with a contact form $\lambda$. We recall that a contact form
$\lambda$ is a 1-form on $\mathcal{Q}$ such that $\lambda\bigwedge
(d\lambda)^{n-1}$ is a volume form. Associated to $(\mathcal{Q},\lambda)$
we have the contact structure $\xi=\ker(\lambda)$, which is a
$(2n-2)$-dimensional subbundle of $T\mathcal{Q}$, and $(\xi, d\lambda|_{\xi})$
defines a symplectic vector bundle. Furthermore, there is a unique
nonvanishing vector field $X=X_{\lambda}$, called the Reeb vector field,
defined by the condition
$$i_X\lambda=1,\;\;\;\;i_Xd\lambda=0.$$
We have a canonical splitting of $TQ$,
$$ T\mathcal{Q}=\mathbb R X\oplus\xi,$$
where $\mathbb R X$ is the line bundle generated by $X$.

\v\v

\subsection{Neighbourhoods of Lagrangian submanifolds}\label{neighbor_L}

Let $(M,\omega)$ be a compact symplectic manifold, $L\subset M$ be a compact Lagrangian sub-manifold.
The following Theorem is well-known.
\begin{theorem}\label{Lagrangian neighbourhood theorem}
 Let $(M,\omega)$ be a symplectic manifold of dimension
$2n$, and $L$ be a compact Lagrangian submanifold. Then there exists a neighbourhood $U\subset T^*L$ of the zero section, a neighbourhood $V \subset M$ of $L$, and a diffeomorphism $\phi:U\rightarrow V$ such that
\begin{equation}
\phi^{*}\omega=-d \Lambda,\;\;\;\; \phi|_{L}=id,
\end{equation}
where $\Lambda$ is the canonical Liouville  form.
\end{theorem}

 Let $(x_1,\cdots, x_n)$ be a
local coordinate system on $O\subset L$, there is a canonical coordinates
$$(x_1,\cdots,x_n,y_1,\cdots,y_n)$$ on $T^{*}O=T^{*}L|_{O}.$ In terms of this coordinates
the  Liouville form can be written as
$$\Lambda =\sum y_{i}dx_{i}.$$
Let $\pi :T^{*}L \rightarrow L$ be the canonical projection. There is a global defined vector field $W$ in $T^*L$ such that  in the local coordinates of $\pi^{-1}(O),$  $W$ can be written as
\begin{equation}
W|_{\pi^{-1}(O)}=-\sum_{i=1}^{n}y_{i}\frac{\partial }{\partial y_{i}}.
\end{equation}
Suppose that given a Riemannian metric on $L$, in terms of the coordinates $x_1,...,x_n$, $g_{L}=\sum\limits_{i,j=1}^{n} g_{ij}dx_idx_j.$ It naturally induced a metric on $T^{*}L.$
Let $$V=\frac{W}{\|W\|}.$$
$V$ is a global defined vector field on $T^{*}L-L$.
\v
Denote by $S^{n-1}(1)$ (resp.$B_{1}(0)$) the Euclidean unit sphere (resp. the Euclidean unit ball). Consider the coordinates transformation between the sphere coordinates and the Cartesian coordinate
\begin{align}
\Psi:(0,1]\times S^{n-1}(1)&\rightarrow  B_{1}(0) \nonumber\\
(r,\theta_{1},\cdots,\theta_{n-1})&\rightarrow   (y_1,\cdots,y_{n}).
\end{align}
Consider the unit sphere bundle $\widetilde{M}$ and the unit ball bundle $\mathbb D_{1}(T^{*}L)$ in $T^{*}L,$ in terms of the coordinates $(x_{1},\cdots,x_{n},y_{1},\cdots,y_n)$
\begin{align}
 \widetilde{M}|_{\pi^{-1}(O)}&=\{ (x_{1},\cdots,x_{n},y_{1},\cdots,y_n)\in \pi^{-1}(O)\;|\; \sum_{i,j=1}^{n} g^{ij}(x)y_{i}y_{j}=1 \},\\
\mathbb D_{1}(T^{*}L)|_{\pi^{-1}(O)}&=\{ (x_{1},\cdots,x_{n},y_{1},\cdots,y_n)\in \pi^{-1}(O)\;|\; \sum_{i,j=1}^{n} g^{ij}(x)y_{i}y_{j}\leq 1 \}.
\end{align}
Denote $\lambda= -\Lambda \mid_{\widetilde{M}}$.
We have
$$\Lambda=-\|y\|\lambda.$$
$\lambda $ is a contact form, i.e., $(\widetilde{M}, \lambda)$ is a contact manifold. Put $\xi=\ker(\lambda)$. Then
$X\mid_{\widetilde{M}}=- \sum g^{ij}y_i\frac{\partial}{\partial x_j}$ is the Reeb vector field, and $V\mid_{\widetilde{M}}=\sum_{i=1}^{n}y_{i}\frac{\partial }{\partial y_{i}}.$

The map $\Psi$ induced a map $\tilde \Psi :(0,1]\times \widetilde{M} \rightarrow \mathbb D_{1}(T^{*}L).$
Through $\tilde{\Psi}$ we consider $\mathbb D_{1}(T^{*}L)|_{\pi^{-1}(O)}-L$ as
$(0,1]\times \widetilde{M} $. By Theorem \ref{Lagrangian neighbourhood theorem} we consider $M-L$ as
$$M^+= M_0^{+}\bigcup\{(0,1]\times \widetilde{M}\}$$
with the symplectic form
\begin{equation}\omega_{\phi} =- d\Lambda=
 \|y\| d\lambda + d \|y\| \wedge\lambda,
\end{equation}
where $M^+:= M-L$ and $M^{+}_{0}$ is a compact symplectic manifold with boundary.
\v
We choose the  {\em neck stretching  technique}.

\v

Let $\phi :[0,\infty)\rightarrow (0,\ell]$ be a smooth
function satisfying, for any $k>0,$
\begin{itemize}
\item[(1)] $\phi^{\prime}<0,
\;\;\phi(0)= \ell, \;\;\phi(a)\rightarrow 0 \;\;as\;\; a
\rightarrow \infty$,
\item[(2)] $\lim\limits_{a\to 0^{+}} \frac{\partial^{k} \phi}{\partial a^{k}}=0.$
\end{itemize}

Through $\phi$ we consider
$M^+$ to be $M^{+} = M_{0}^{+}\bigcup\{[0,\infty)\times
\widetilde{M}\} $ with symplectic form
$\omega_{\phi}|_{M_0^{+}}=\omega$, and over the cylinder $[0,\infty)\times \widetilde{M}$
\begin{equation}\omega_{\phi} =- d\Lambda=
\phi d\lambda + \phi^{\prime}d a\wedge \lambda.
\end{equation}
 Moreover, if we choose the origin of ${\mathbb
R}$
tending to $\infty $, we obtain ${\mathbb
R}\times \widetilde{M}$ in
the limit.

\v

Choose $\ell_0 < \ell $ and denote $$\Phi^+ =\left \{ \phi :[0,
\infty )\rightarrow (0, \ell_0] | \phi^{\prime}
< 0 \right \}.$$
Let $\ell_1 < \ell_2$ be two real numbers satisfying $
0< \ell_1 < \ell_2 \leq  \ell_0 .$ Let $\Phi_{\ell_1,\ell_2}$ be the
set of all smooth functions $\phi:{\mathbb
R}\rightarrow
(\ell_1,\ell_2)$ satisfying $$\phi^{\prime}<0,\;\;
\phi(a)\rightarrow \ell_1\;\:\;{\rm as} \;a\rightarrow \infty ,\;\;
\phi(a)\rightarrow \ell_2\;\:\;{\rm as} \;a\rightarrow -\infty.$$ To
simplify notations we use $\Phi $ to denote both $\Phi^+$ and
$\Phi_{\ell_1,\ell_2}$, in case this does not cause confusion.
\vskip 0.1in \noindent
We fixed $\phi\in \Phi$ and consider the symplectic manifold $\left(M^{+}_{0}\bigcup\{[0,\infty)\times \widetilde{M}\},\omega_{\phi}\right).$ For any different $\phi_{1}\in \Phi$, we have
$\omega_{\phi}|_{M^{+}_0}=\omega_{\phi_1}|_{M^{+}_0}$ and $\phi\circ \phi_{1}^{-1}$ is a symplectic diffeomorphic over cylinder part $\{(0,1]\times \widetilde{M}\}$.

\v

\subsection{Cylindrical almost complex structures}\label{subsection:7.2}

Let
\begin{equation}\label{cylinder}
M^+=M^{+}_{0}\bigcup\left\{[0,\infty)\times \widetilde{M}\right\}\end{equation}
be a symplectic manifold with cylindrical end,
where $\widetilde{M}$ be a compact contact manifold with contact form $\lambda$. Denote by $\omega_{\phi}$  the symplectic form of $M^{+}$ such that
$\omega_{\phi}|_{M_0^{+}}=\omega$, and over the cylinder $[0,\infty)\times \widetilde{M}$
\begin{equation}\omega_{\phi} =- d\Lambda=
\phi d\lambda + \phi^{\prime}d a\wedge \lambda.
\end{equation} We also consider $\mathbb R\times \widetilde{M}$. Denote by $N$ one of $M^+$ and $\mathbb R\times \widetilde{M}$.
\v
Put $\xi=\ker(\lambda)$, and denote by $X$ the Reeb vector field defined by
$$\lambda(X)=1,\;\;\; d\lambda (X)=0.$$
We choose a $d \lambda$-tame almost complex structure $\widetilde{J}$ for the symplectic vector bundle
$(\xi, d\lambda)\rightarrow \widetilde{M}$ such that
\begin{equation}
g_{\widetilde{J}(x)}(h,k)=\frac{1}{2}\left(d\lambda(x)(h,\widetilde{J}(x)k)+d\lambda(x)(k,\widetilde{J}(x)h)\right),
\end{equation}
for all $x\in \widetilde{M},\;h,k\in\xi_x$, defines a smooth fibrewise
metric for $\xi$.  We assume that we can choose $\widetilde{J}$ such that on every periodic orbit $L_X \widetilde{J}=0$.
Denote by $\Pi:T\widetilde{M}\rightarrow\xi$ the projection along $X$.
We define a Riemannian metric $\langle\;,\;\rangle$ on $\widetilde{M}$ by
\begin{equation}\langle h,k\rangle=\lambda(h)\lambda(k)+g_{\widetilde{J}}(\Pi h,\Pi k)
\end{equation}
for all $h,k\in T\widetilde{M}$.

\vskip 0.1in
Given a $\widetilde{J}$ as above there is an associated almost complex
structure $J$ on $\mathbb R\times \widetilde{M}$ defined by
\begin{equation}\label{complex structure}
J\mid_{\xi}= \tilde{J},\;\;\;JX=\frac{\partial}{\partial a},\;\;J(\frac{\partial}{\partial a})=-X,
\end{equation}
where $a$ is the canonical coordinate in $\mathbb R$. It is easy to check that $J$ defined by \eqref{complex structure}
is  $\omega_{\varphi}$-tame over the cylinder end. We can choose an almost complex structure $J$ on $M^+$ such that $J$ is tamed by $\omega$ and over the cylinder end  $J$ is given by \eqref{complex structure}.
\v
There is a canonical coordinate system for ${\mathbb
R}\times \widetilde{M}$ and for cylinder end of $M^+$, but still there
are some freedom of choosing the coordinates:
\v
When we write $M^+$ as \eqref{cylinder} we have chosen a coordinate $a$ over the cylinder part. We can choose different coordinate $\hat{a}$ over the
cylinder part such that
\begin{equation}\label{rescaling 1}
a=\hat{a}+ C
\end{equation}
for some constant $C>0$. Similarly, for $\mathbb
R \times \widetilde{M}$ we can choose $\hat{a}$ such that
\begin{equation}\label{rescaling 3}
a=\hat{a}\pm C
\end{equation}
for some constant $C>0$.
\v
 \vskip 0.1in \noindent For any $\phi \in \Phi $
\begin{equation}\label{omega_forms}
\langle v,w\rangle_{\omega_{\phi}} = \frac{1}{2}\left(
\omega_{\phi} (v,Jw) + \omega_{\phi} (w,Jv)\right) \;\;\;\;\;
\forall \;\; v, w \in TN
\end{equation}
defines a Riemannian metric
on $N$.  Note that $\langle \;,\;\rangle_{\omega_{\phi}}$ is not
complete. We choose another metric $\langle \;,\;\rangle$ on $N$
such that
\begin{equation}\label{omega_forms_on_M0}\langle \;,\;\rangle = \langle
\;,\;\rangle_{\omega_{\phi}} \;\;\;\;on \;\; M^{+}_{0}
\end{equation}
 and over
the tubes
\begin{equation} \label{omega_forms_on_tubes}
\langle(a,v),(b,w) \rangle= ab + \lambda (v)\lambda
(w) + g_{\widetilde{J}}(\Pi v, \Pi w),
\end{equation}
 where we
denote by $\Pi:T\widetilde{M}\rightarrow\xi$ the projection along
$X$. It is easy to see that $\langle \;,\;\rangle$ is a complete
metric on $N$. \vskip 0.1in \noindent

\subsection{$J$-holomorphic maps with finite energy}

Let $(\Sigma,i)$ be a
compact Riemann surface and $P\subset\Sigma$ be a finite
collection of puncture points. Denote $\stackrel{\circ}{\Sigma}
=\Sigma\backslash P.$ Let ${u}:\stackrel{\circ}{\Sigma}
\rightarrow N$ be a ${J}$-holomorphic map, i.e., ${u}$ satisfies
\begin{equation}\label{j_holomorphic_maps}
d{u}\circ i={J}\circ d{u}.
\end{equation}
Following \cite{HWZ1} we
impose an energy condition on $u$. For any $J$-holomorphic map
$u:\stackrel{\circ}{\Sigma}\rightarrow N$ and any $\phi \in \Phi $
the energy $E_{\phi}(u)$ is defined by
\begin{equation}\label{definition_of_energy}
E_{\phi}(u)=\int_{ {\Sigma}}u^{\ast}\omega_{\phi}.
\end{equation}
Let $z=e^{s+2 \pi \sqrt{-1}t}.$ One
computes over the cylindrical part
\begin{equation}\label{omega_forms_cylinder}
u^{\ast}\omega_{\phi}=
(\phi d\lambda \left((\pi\widetilde{u})_s,
(\pi\widetilde{u})_t)\right) - {\phi}^{\prime}(a^{2}_s + a^{2}_t
))ds\wedge dt,
\end{equation} which is a nonnegative integrand.
A $J$-holomorphic map $u:\stackrel{\circ}{\Sigma} \rightarrow N $
is called a finite energy $J$-holomorphic map if over the cylinder end
\begin{equation}\label{finite_energy_j_holomorphic_maps}
\sup_{\phi \in
\Phi }\left \{\int_{{\Sigma}}u^{*} \omega_{\phi}
\right \}+\int_{\Sigma}u^{*}d\lambda <\infty.
\end{equation}
  For a $J$-holomorphic
map $u:{\Sigma} \rightarrow {\mathbb
R}\times\widetilde{M}$  we write $u=(a, \widetilde{u})$ and define
\begin{equation}\label{definition_of_energy_on_complex_manifolds}
\widetilde{E}(u)=\int_{{\Sigma}}\widetilde{u}^{\ast}
 d\lambda .
\end{equation}
Denote
$$\widetilde{E}(s)=\int_s^{\infty}\int_{S^1}\widetilde{u}^{\ast}(
d\lambda).$$ Then $$\widetilde{E}(s)=\int_s^{\infty}\int_{S^1} |\Pi
\widetilde{u}_t |^2dsdt,$$
\begin{equation}\label{deriative_of_energy_on_Z}
\frac{d\widetilde{E}(s)}{ds}=-\int_{S^1} |\Pi\widetilde{u}_t |^2dt.
\end{equation}
Here and later we use $ |\cdot |$ denotes the norm with respect to the metric defined by \eqref{omega_forms_on_tubes}.

 \v\n
The following two lemmas are well-known ( see \cite{H}): \vskip 0.1in \noindent
\begin{lemma}\label{zero_energy_classification}
 \begin{itemize}
 \item[(1)] Let $u=(a,\widetilde{u}):\mathbb{C} \rightarrow {\mathbb R}
\times \widetilde{M}$ be a $J$-holomorphic map with finite energy.
If $\int_{ \mathbb{C}}\widetilde{u}^{\ast}(\pi^{\ast}d\lambda)=0$, then
$u$ is a constant.
\item[(2)] Let $u=(a,\widetilde{u}):{
\mathbb{R}}\times S^1 \rightarrow {\mathbb R}\times \tilde{M}$ be a
$J$-holomorphic map with finite energy. If $\int_{\mathbb{R}\times
S^1}\widetilde{u}^{\ast}(\pi^{\ast}d\lambda)=0$, then
$(a,\widetilde{u})=(kTs+c, kt+d)$, where $k\in {\mathbb Z}^{+}$, $c$ and
$d$ are constants.
\end{itemize}
\end{lemma}
\vskip 0.1in
 \noindent
\begin{lemma}\label{sequence_convergence_of_j_holomorphic_maps} Let
$u=(a,\widetilde{u}):\mathbb{{C}}-D_1 \rightarrow {\mathbb R} \times
\widetilde{M}$ be a nonconstant $J$-holomorphic map with finite
energy. Put $z=e^{s+2\pi \sqrt{-1}t}$. Then for any sequence
$s_i\rightarrow \infty $ , there is a subsequence, still denoted
by $s_i$, such that $$\lim_{i\rightarrow \infty}\widetilde
u(s_i,t)=x(kTt)$$ in $C^{\infty}(S^1)$ for some $kT$-periodic orbit
$x(kTt)$.
\end{lemma}

 \vskip 0.1in
 \noindent

\subsection{Periodic orbits of Bott-type}

Let $\mathcal{F}\subset \widetilde M$ be the locus of minimal periodic orbits with $\mathcal F=\bigcup_{i=1}^{\ell}\mathcal{F}_{i}$, where each $\mathcal{F}_{i}$ is a connected component of $\mathcal{F}$ with  minimal periodic $T_{i}$. We assume that
\begin{assumption}\label{bott_assumption}
\begin{itemize}
\item[(1)] every $\mathcal{F}_{i}$ is either non-degenerate or of Bott-type;
\item[(2)] let $\mathcal{F}_{i}$ be of Bott-type, then there exists  a free $S^{1}$-action on $\mathcal {F}_{i}$ such that $Z_{i}=\mathcal{F}_{i}/S^{1}$ is a closed, smooth manifold. Set $n_{i}=\dim(\mathcal{F}_{i});$
\item[(3)] for every periodic orbit $x,$ there is a smooth submanifold $\Re_{x}\subset \widetilde{M}$ of dimension $>2$ such that $d\lambda \mid _{\Re_{x}}=0$ and $x\subset\Re_{x}$, and the almost complex structure $\tilde{J}$ can be chosen such that $L_X \tilde{J}=0$ along $x$.
\end{itemize}
\end{assumption}
\v

Let $u=(a,\widetilde{u}):\mathbb{{C}}-D_1 \rightarrow {\mathbb R}  \times \widetilde{M}$ be a
$J$-holomorphic map with finite energy. Put $z=e^{s+2\pi it}$. Assume that there exists a sequence $s_{i}\to\infty$ such that $\tilde { u} (s_{i},t)\longrightarrow x(kTt)$ in $C^{\infty}(S^1,\widetilde {M})$ as $i\rightarrow \infty$ for some $k\in {\mathbb Z}$, where $T$ is the minimal periodic.    Following Hofer (see \cite{HWZ1}) we introduce a convenient local coordinates near the periodic orbit $x$. Since $S^1=\bR/\bZ$, we work in the covering space $ \bR \rightarrow S^1$.
\v
\begin{lemma}  \label{Darboux}
Let $(\widetilde{M},\lambda)$ be a $(2n-1)$-dimensional compact manifold, and let $x(kt)$ be a $k$-periodic orbit with the  minimal periodic  $T.$ Then there is an open neighborhood $U\subset S^1\times \mathbb R^{2n-2}$ of $S^1\times \{0\}$ with coordinates $(\vartheta, w_1,\cdots, w_{2n-2})$ and an open neighborhood $\Im \subset \widetilde{M}$ of $\{x(t)|t\in \mathbb R\}$ and a diffeomorphism $\psi : U\rightarrow \Im$ mapping $S^1\times \{0\}$ onto $\{x(t)|t\in \mathbb R\}$ such that
\begin{equation}
\psi^*\lambda = g\lambda_0,\end{equation}
where $\lambda_0=d\vartheta + \sum w_i dw_{n+i-1}$ and $g:U\rightarrow \mathbb R$ is a smooth function satisfying
\begin{equation}\label{standard_lambda-1}
g(\vartheta, 0)=T,\;\;dg(\vartheta,0)=0
\end{equation}
for all $\vartheta\in S^1$.
\end{lemma}

\begin{remark} We call the coordinate system $(\vartheta, {\bf w})$ in Lemma \ref{Darboux} a pseudo-Darboux coordinate system, and call the following transformation of two local pseudo-Darboux coordinate systems
\begin{equation}\label{canonical}
(\vartheta, {\bf w})\longrightarrow (\check{\vartheta},{\bf w}),\;\;\check{\vartheta}=\vartheta + \vartheta_0\end{equation}
a canonical coordinate transformation, where $\vartheta_0$ is a constant.
\end{remark}
 \v
The following theorem is well-known (see \cite{LR,BH,OW})
\begin{theorem}\label{exponential_estimates_theorem} Suppose that $\mathcal{F}$ satisfies Assumption \ref{bott_assumption}. Let $u:\mathbb C-D_1\rightarrow {\mathbb R}\times \tilde{M}$ be a J-holomorphic map with finite energy.
Put
$z=e^{s+2 \pi \sqrt{-1}t}$. Then $$\lim_{s\rightarrow \infty}\widetilde
u(s,t)=x(kTt) $$ in $C^{\infty}(S^1)$ for some $kT$-periodic orbit
$x,$ and there are constants $\ell_0$,
$\vartheta_0$  such that for any $0<\mathfrak{c}<\min\{\frac{1}{2},\frac{\mathcal C_{1}^2}{2}\}$ and for all ${\bf n}=(n_1,n_2)\in {\mathbb Z_{\geq 0}^2 }$
\begin{eqnarray}
\label{exponential_decay_a}
|\partial^{\bf n}[a(s,t)-kTs-\ell_0]|\leq \mathcal C_{\bf n} e^{-\mathfrak{c}|s|}\\
\label{exponential_decay_theta}
|\partial^{\bf n}[\vartheta(s,t)-kt-\vartheta_0]|\leq \mathcal C_{\bf n} e^{-\mathfrak{c}|s|} \\
\label{exponential_decay_y}|\partial^{\bf n} {\mathbf w}(s,t)|\leq \mathcal C_{\bf n}
e^{-\mathfrak{c} |s|},
\end{eqnarray}
 where $C_{\bf n}$ are constants. Here $(\vartheta,{\bf w})$ is a  pseudo-Darboux coordinate near the periodic orbit $x$.
\end{theorem}
\v

Let $x_{o}$ be a minimal periodic orbit.
  We choose a local  pseudo-Darboux  coordinate on an open set  $\Im\subset \widetilde M$ near $x_{o}$. Let  $\mathcal{O}$ be an open set  such that $\overline{\mathcal{O}}\subset \Im$ is compact and $x_{o}\subset \mathcal{O}.$ We fix a positive constant $  \mathsf{C} .$
Denote by $\mathfrak{S}$ the class of $J$-holomorphic maps with finite energy $u:[s'_0,\infty)\times S^1\rightarrow {\mathbb R}\times \tilde{M}$ satisfying
\begin{enumerate}
\item[(i)] $\sup_{\phi\in\Phi} E_{\phi}(u)\leq \frac{1}{2}\hbar$, where $\hbar$ is the constant in  Theorem \ref{lambda_lower_bound} and Lemma \ref{lower_bound_of_singular_points},
\item[(ii)] $\lim_{s\rightarrow \infty}\widetilde u(s,t)=x(kTt)$ in $C^{\infty}(S^1)$ with $x\subset \mathcal{O}$,
\item[(iii)] $\widetilde{u}([s'_0,\infty)\times S^1)$ lie in the pseudo-Darboux coordinate system on $\Im$,
\item[(iv)] there is a ball $D_{1}(s_0,t_{0})\subset[s'_0,\infty)\times S^1$  such that $\widetilde{u}(D_{1}(s_0,t_{0}))\subset \mathcal{O}$ and $|a(D_{1}(s_0,t_{0}))|\leq \mathsf{C}$.
\end{enumerate}
 In our next paper \cite{LS1} we will prove the following theorem.
\begin{theorem}\label{exponential_decay_uniform}
Let $u\in \mathfrak{S}$.
  Then   the constants $C_{\bf n}$ in Theorem \ref{exponential_estimates_theorem} depend only on $\mathcal C_{1}, s_0,{\mathsf C}$,  ${\bf n}$, $\hbar,$ $ \mathfrak{c}$ and $\mathcal{O}$. Moreover, we have
\begin{eqnarray}\label{control_ell}
|\ell_{0}|\leq  C_{1}  ,
\end{eqnarray}
where $C_{1}$ depends only on  $\hbar$, $ C_{\bf n} ,  \mathfrak{c}, s_{0},\mathsf{C}$ and $\mathcal{O}.$
\end{theorem}

 Following \cite{HWZ1} we introduce functions
\begin{equation}
 a^{\diamond}(s,t)=a(s,t)-ks,\;\; \vartheta^{\diamond}(s,t)=\vartheta(s,t)-kt.
\end{equation}
Denote
\begin{equation}
\pounds=(a^{\diamond},  \vartheta^{\diamond}).
\end{equation}
Set $b_{i}=0,\;b_{n-1+i}({\bf w})=w_{i},  \forall \; i=1,\cdots,n-1,$ and
$$e_{i}=\frac{\p}{\p w_{i}}-b_{i}\frac{\p}{\p \vartheta},\;\;\;\;i=1,\cdots,2n-2.  $$
Then $\xi=$span$\{e_1,\cdots,e_{2n-2}\}.$ Denote $\tilde Je_{i}=\sum\tilde J_{ij}e_{i}.$ In the basis $\p_{\vartheta},e_{1},\cdots,e_{2n-2},$ the Reeb vector field can be re-written as
\begin{equation}
X=\frac{1}{g}\p_{\vartheta} +\frac{1}{g^2}\left(\sum_{i\leq n-1}e_{n-1+i}(g)e_{i}-\sum_{i\geq n }e_{i-n+1}(g)e_{i}\right).
\end{equation}
Let $\tilde{\bf X }=\frac{1}{g^2}(e_{n}g,\cdots,e_{2n-2}g,-e_{1}g,\cdots, -e_{n-1}g ).$
 By \eqref{j_holomorphic_maps}, we have
$$\pounds_s + J\pounds_t= h,\;\;\;\; {\bf w}_{s}+\widetilde{J}{\bf w}_{t}+a_{t}\tilde{{\bf X}}-a_{s}\tilde J\tilde{{\bf X}} =0.$$
where $$h=(\Sigma b_i({\bf w})(w_{i})_{t}+(\vartheta_{t}+\Sigma b_i({\bf w})(w_{i})_{t})(\Sigma f_{i}w_{i}), -\Sigma b_i({\bf w})(w_{i})_{s}-(\vartheta_{t}+\Sigma b_i({\bf w})(w_{i})_{s})(\Sigma f_{i}w_{i}))$$
and $f_{i}=\int_{0}^1\p_{w_{i}}g(\vartheta,\tau {\bf w})d\tau$. Let $V=\pounds_{t}$ and $g=h_{t}.$ Denote
\begin{equation}
\hat{E}(\pounds)=\int_{\Sigma}\|\pounds_{s}\|^2+\|\pounds_{t}\|^2 dsdt.
\end{equation}

 In our next paper \cite{LS1} we will  prove the following  theorem.

\begin{theorem}\label{tube_exponential_decay} Suppose that $\widetilde{M}$  satisfies the Assumption \ref{bott_assumption}.
Let $u:[-R ,R ]\times S^1\rightarrow {\mathbb R}\times \widetilde{M}$ be a $J$-holomorphic maps with finite energy. Assume that
\begin{enumerate}
\item[(i)] $ \sup\limits_{\phi\in \Phi} E_{\phi}(u,-R\leq s\leq R)+\hat{E}(\pounds)\leq \frac{1}{2}\hbar$,
\item[(ii)] $\widetilde u([-R,R]\times S^1)$  lie in a pseudo-Darboux coordinate system $(\vartheta,{\bf w})$ on $\Im$,
\item[(iii)] $\sum\limits_{n_1,n_{2}\leq 2}\|\nabla^{\bf n} u(-R,\cdot)\|_{L^2(S^{1})}\leq C_{2},\;\;\;\;\sum\limits_{n_1,n_{2}\leq 2}\|\nabla^{\bf n}u(R,\cdot)\|_{L^2(S^{1})}\leq C_{2},$  where ${\bf n}=(n_1,n_2),$
\end{enumerate} Then there exist constants  $\mathcal{C}_{1}>0$ and   $0<\fc< \frac{1}{2} $ depending only on $\tilde J$ and $C_2$ such that
 \begin{align}
|\nabla {\bf w}|(s,t)\leq \mathcal{C}_{1}e^{-\fc(R-|s|)},\;\;|\nabla \pounds|\leq \mathcal{C}_{1}e^{-\fc(R- |s|)},\;\;\;\forall\; |s|\leq R-1,
\end{align}
 \end{theorem}

In this paper, we also need the following  implicit function theorem (see \cite{MS}).
\begin{lemma}\label{details_implicit_function_theorem}
Let  $X$ and $Y$ be Banach spaces, $U\subset X$ be an open set, and $\ell$ be a positive integer. If $F:U\longrightarrow Y$ is of class $C^{\ell}$. Let $x_{0}\in U$ be such that $D:=dF(x_{0}):X\longrightarrow Y$ is surjective and has a bounded linear right inverse $Q:Y\longrightarrow X.$ Choose positive constants $\hbar_{1} $ and $C$ such that $\|Q\|\leq C,$ $B_{\hbar_{1} }(x_{0},X)\subset U$ and
\begin{equation}\label{continue_of_differential}
\|dF(x)-D\|\leq \frac{1}{2C},\;\;\;\;\forall \;x\in B_{\hbar_{1}}(x_{0})
\end{equation}
where $B_{\hbar_{1} }(x_{0})=\{x\in X|\; \|x-x_{0}\| \leq \hbar_{1} \}.$ Suppose that $x_1\in X$ satisfies
\begin{equation}\label{small_value_of_F}
\|F(x_1)\|<\frac{\hbar_{1} }{4C},\;\;\;\; \|x_{1}-x_{0}\|\leq \frac{\hbar_{1}}{8}.
\end{equation}
Then there exists a unique $x\in X$ such that
\begin{equation}
F(x)=0,\;\;\;\;x-x_1\in Im\; Q,\;\;\;\;\|x-x_{0}\|\leq \hbar_{1} ,\;\;\;\; \|x-x_{1}\|\leq 2C\|F(x_{1})\|.
\end{equation}
Moreover, write $x:=x_0 + \xi + Q \circ f(\xi)$, $\xi\in ker D$, then $f$ is of class $C^{\ell}$.
\end{lemma}

\v

\section{Weighted sobolev norms}\label{weight_norm}

Consider ${\mathbb{R}}\times \widetilde{M}$ and $M^+=M^{+}_{0}\bigcup\left\{[0,\infty)\times \widetilde{M}\right\}.$ Let $N$ be one of ${\mathbb{R}}\times \widetilde{M}$, $M^+$. Suppose that  $\Sigma= \bigcup\limits_v \Sigma_{v}$ is Riemann surface with nodal points $\{q_{1},\cdots,q_{\mathfrak{I}}\}$, puncture points $\{p_{1},\cdots,p_{\nu}\}$ and  $u:{\Sigma}\rightarrow \bigcup N_i$ is a continuous map such that the restriction of $u$ to each smooth component is smooth, where $\bigcup N_i$ denotes the union of some copy of $N$. We choose cylinder coordinates $(s,t)$ on $\Sigma$ near each nodal point and each puncture point. We choose a local pseudo-Darboux coordinate system near each periodic orbit on $N$. Let $\stackrel{\circ}{\Sigma}=\Sigma-\{q_{1},\cdots,q_{\mathfrak{I}},p_{1},\cdots,p_{\nu} \}$.
\v

Over each tube the linearized operator
$D_{u}$ takes the following form
\begin{equation}
D_{u}=\frac{\partial}{\partial
s}+J_0\frac{\partial} {\partial t}+S = \bar{\partial}_{J} +
S.
\end{equation}
By exponential decay we have
$$
\left|\frac{\partial^{k+l} }{\partial s^k \partial t^{l} }S\right|\leq C_{k,l}e^{- \mathfrak{c} s}
$$
for some constant $C_{k,l}>0$ for $s$ big enough.
Therefore, the operator $H_s=J_0\frac{d}{dt}+S$ converges
to $H_{\infty}=J_0\frac{d}{dt}$. Obviously, the operator $D_u$ is
not Fredholm operator because over each   puncture and node the
operator $H_{\infty}=J\frac{D}{dt}$ has zero eigenvalue. The $\ker
H_{\infty}$ consists of constant vectors. To recover a Fredholm theory we use
weighted function spaces. We choose a weight $\alpha$ for each
end. Fix a positive function $W$ on $\Sigma$ which has order equal
to $e^{\alpha |s|} $ on each end, where $\alpha$ is a small
constant such that $0<\alpha<\mathfrak{c} $ and over each end
$H_{\infty}- \alpha = J_0\frac{d}{dt}- \alpha $ is invertible. We
will write the weight function simply as $e^{\alpha |s|}.$ Denote by $C(\Sigma;u^{\ast}T(\bigcup N_i))$ all tangent vector fields $h$ on $\bigcup N_i$ along $u$ satisfying
\v
{\bf (a)} $h\in C^{0}({\Sigma},u^*T(\bigcup N_i))$,
\v
{\bf (b)} the restriction of $h$ to each smooth component is smooth.
\v
 For any
section $h\in C(\Sigma;u^{\ast}T(\bigcup N_i))$ and section $\eta \in
\Omega^{0,1}(u^{\ast}T(\bigcup N_i))$ we define the norms
\begin{eqnarray}\label{def_1_p_alpha}
&&\|h\|_{1,p,\alpha}=\sum_{v}\left(\int_{\Sigma_v}( |h|^p+ |\nabla
h|^p)d\mu\right)^{1/p}
+\sum_{v}\left(\int_{\Sigma_v}e^{2\alpha|s|}(|h|^2+|\nabla h|^2)
d\mu\right)^{1/2} \\ \label{def_p_alpha}
&&\|\eta\|_{p,\alpha}=\sum_{v}\left(\int_{\Sigma_v}|\eta|^p d\mu
\right)^{1/p}+ \sum_v\left(\int_{\Sigma_v}e^{2\alpha|s|}|\eta|^2
d\mu\right)^{1/2}
\end{eqnarray} for $p\geq 2$, where all norms and
covariant derivatives are taken with respect to the  metric
$\langle\;\;\rangle$ on $u^{\ast}T(\bigcup N_i)$ defined in \eqref{omega_forms_on_tubes}, and
the metric on $\stackrel{\circ}{\Sigma}$. Denote
\begin{eqnarray}
&&{\mathcal C}(\Sigma;u^{\ast}T(\bigcup N_i))=\{h
\in C(\Sigma;u^{\ast}T(\bigcup N_i)); \|h\|_{1,p,\alpha}< \infty
\},\\
&&{\mathcal C}(u^{\ast}T(\bigcup N_i)\otimes \wedge^{0,1})
=\{\eta\in \Omega^{0,1}(u^{\ast}T(\bigcup N_i)); \|\eta\|_{p,\alpha}< \infty
\}.
\end{eqnarray}
 Denote by $W^{1,p,\alpha}(\Sigma;u^{\ast}T(\bigcup N_i))$ and
$L^{p,\alpha}(u^{\ast}T(\bigcup N_i)\otimes \wedge^{0,1})$ the completions of
${\mathcal C}(\Sigma;u^{\ast}T(\bigcup N_i))$ and ${\mathcal C}(u^{\ast}T(\bigcup N_i)\otimes
\wedge^{0,1}) $ with respect to the norms \eqref{def_1_p_alpha} and \eqref{def_p_alpha}
respectively. Then the operator $D_u:   W^{1,p,\alpha}\rightarrow L^{p,\alpha}$
is a Fredholm operator.
\vskip 0.1in
\noindent
\v
For each puncture point $p_j,j=1,...,\nu $, let $h_{j0}\in (T_{p_{j}}(\mathcal{F}_{\fe_j})\oplus (span\{\frac{\p}{\p a}\}) $.
For each bounded nodal $q_{i},$ denote  $\bh_{q_{i}}=T_{q_{i}}N$,  let $h_{(i+\nu)0}\in \bh_{q_{i}}$;  for each unbounded nodal $q_{i},$ denote  $\bh_{q_{i}}=(T_{q_{i}}(\mathcal{F}_{i_{j}})\oplus (span\{\frac{\p}{\p a}\})$ and let $h_{(i+\nu)0}\in \bh_{q_{i}}$, where $u:\Sigma\rightarrow N$ converges to $kT$ periodic orbit $x(kTt)\subset \mathcal{F}_{i_j}$ as $z\rightarrow q_{i}$.
Put
$\bh= \left(\oplus_{j=1}^\nu (T_{p_{j}}(\mathcal{F}_{\fe_j})\oplus (span\{\frac{\p}{\p a}\})\right) \bigoplus \left(\oplus_{i=1}^{\mathfrak{I}} \bh_{q_{i}}\right) ),$ $h_{0}=(h_{1 0},...,h_{\nu 0},h_{(1+\nu)0},...,h_{  (\mathfrak{I}+\nu) 0}).$

$h_0$ may be considered as a vector field in the coordinate neighborhood.
We fix a cutoff function $\rho$:
\[
\rho(s)=\left\{
\begin{array}{ll}
1, & if\ |s|\geq d, \\
0, & if\ |s|\leq \frac{d}{2}
\end{array}
\right.
\]
where $d$ is a large positive number. Put
$$\hat{h}_0=\rho h_0.$$
Then for $d$ big enough $\hat{h}_0$ is a section in $C^{\infty}(\Sigma; u^{\ast}TN)$
supported in the tube $\{(s,t)||s|\geq \frac{d}{2}, t \in S ^1\}$.
Denote
$${\mathcal W}^{1,p,\alpha}=\{h+\hat{h}_0 | h \in
W^{1,p,\alpha},h_0 \in \bh\}.$$
\vskip 0.1in
\noindent
 We define the weighted Sobolev  norm  on ${\mathcal W}^{1,p,\alpha}$ by $$\|( h, \hat{h}_{0})\|_{\Sigma,1,p,\alpha}=
 \|h\|_{\Sigma,1,p,\alpha} + |h_{0}| .$$
Obviously, the operator $D_u:   \mathcal W^{1,p,\alpha}\rightarrow L^{p,\alpha}$
is also a Fredholm operator.
\v


\section{  Moduli spaces of $J$-holomorphic maps}
\v
\subsection{Boundary conditions}

Consider the  symplectic manifold with cylindrical end
$$M^+=M^{+}_{0}\bigcup\left\{[0,\infty)\times \widetilde{M}\right\}.$$
Let $((\Sigma,{\bf j}); {\bf y}, {\bf p})$ be a connected semistable curve with $m$ marked points
${\bf y}=(y_1,...,y_m)$ and $\nu$ puncture points ${\bf p}=(p_1,...,p_{\nu})$,
and $u:{\Sigma} \rightarrow M^{+}$ be a
$J$-holomorphic map.  Let  $ \Sigma=\bigcup\limits_{v=1}^d(\Sigma_{v}, j_v) $ where $(\Sigma_{v},j_v)$ is a smooth Riemann surface and $\pi_{v}:\Sigma_{v} \rightarrow \Sigma$ is a continuous map.
 To describe the boundary conditions we consider two different cases separately:

\v
{\bf Case A .} Moduli space of $J$-holomorphic maps in contact geometry.

\begin{definition} Let ${\mathbf p}=(p_{1},\cdots,p_{\nu})$ be the puncture points.
We assign two weights $(\bk,\fe)$ to ${\mathbf p}$:
\begin{itemize}
\item[(1)]  $\mathbf{k}:{\mathbf p} \rightarrow   \mathbb Z_{>0} $ assigning a $k_{i}$ to each puncture point $p_{i}$, denote ${\mathbf k}=(k_{1},\cdots,k_{\nu}).$
\item[(2)]   $\mathfrak{e}: {\mathbf p} \rightarrow   \mathbb Z_{>0} $ assigning a number $\mathfrak{e_i}$, $1\leq \mathfrak{e_i}\leq \ell$ to each puncture point $p_{i}$.
\end{itemize}
\end{definition}
We call a $J$-holomorphic map $u$ satisfies $(\mathbf{k},\mathfrak{e})$ boundary condition if \v
{ \em
  $u(z)$ converges to a $k_{i}\cdot T_{\fe_{i}}$-periodic
orbit $x(k_{i}\cdot T_{\fe_{i}} t)\subset \mathcal{F}_{ \mathfrak{e_i}}$ as $z$ tends to $p_{i}$.
}
\v

\v
{\bf Case B.} Moduli space of $J$-holomorphic maps in $(M,L)$. \\
As we show in section \S\ref{neighbor_L} that $M-L$ can be considered as $M^+=M^{+}_{0}\bigcup\left\{[0,\infty)\times \widetilde{M}\right\}.$
Let $[c_{i}],i=1,\cdots,\fa$ is a bases in $H_{1}(L;\mathbb Z).$
\begin{definition} Let ${\mathbf p}=(p_{1},\cdots,p_{\nu})$ be the order puncture points.
We assign a weight $\overrightarrow{\mu}$ to ${\mathbf p}$: \\
   $\overrightarrow{\mu}:{\mathbf p} \rightarrow   \mathbb Z_{>0}^{\oplus \fa} $ assigning a $\mu_i=\sum_{l=1}^\fa \mu_{il}[c_{l}]$ to each puncture point $p_{i}$, where   $\mu_{il}\in \mathbb{Z}$.
   Choose the cylinder coordinates $(s_i,t_i)$ near $p_i$.
   We call a $J$-holomorphic map $u$ satisfies $(\overrightarrow{\mu})$ boundary condition if $u$ satisfies
 \begin{equation}
[\pi(\lim\limits_{s_i\rightarrow \infty} u(s_i,S^1))]=\mu_i,\;\;\forall\;1\leq i \leq \nu,
\end{equation}
where $\pi:T^{*}L\rightarrow L$ is  the canonical projection.
\end{definition}

\v

\v

\subsection{Homology}

We fix $A\in H^{2}(M^{+},\Re;\mathbb Z)$ satisfying $\p A=\sum [x(k_{i}T_{\fe_{i}}t)],$ where $x(k_{i}T_{\fe_{i}}t)\subset \mathcal{F}_{\fe_{i}}$ is a $k_{i}T_{\fe_{i}}$ periodic orbit. Consider a $J$-holomorphic map $u$ satisfying
\begin{equation}
[u_{*}(\Sigma)]=A.
\end{equation}
We show that the homology class $A$ give a bound of Energy. To simplify notation we let $(u, (\Sigma,{\bf j}),{\bf y},p)$ be a $J$-holomorphic map converging to a $kT$-periodic
orbit $x(kTt)$ as $z$ tends to $p$, where $x(kTt)$ lie in $\mathcal{F}_i$ and $T$ is the minimal periodic. For the case ${\bf A},$ by the assumption (3) of $\mathcal{F},$ we construct a connected surface $W\subset \Re_{x}$ with boundary $x(kTt)$ ( including the $kT$-periodic).
Then $u_{*}(\Sigma)\cup W$ is a closed surface in $M^{+}$ and
$$[u_{*}(\Sigma)\cup W]\in  H^{2}(M^{+};\mathbb Z).$$
Denote $A=[u_{*}(\Sigma)\cup W].$  By $d\lambda|_{\Re_{x}}=0 $ and $W\subset \widetilde M$ we have
\begin{equation}\label{energy_bound}
 \omega({A})=\int_{u_{*}(\Sigma)}\omega+\int_{W}\omega  =E_{\phi}(u)+\int_{W}d\lambda = E_{\phi}(u).
\end{equation}
 Let $W'\subset \Re_{x}$ be another surface with boundary $x,$ denote $ A'=[u_{*}(\Sigma)\cup W']\in  H^{2}(M^{+};\mathbb Z).$
We have $\omega(A)= \omega({A}')=E_{\phi}(u),$ that is, $E_{\phi}(u)$ is independent of the choice of $W$ in $\Re_{x}$.

\v

For case ${\bf B}$ let $A\in H^{2}(M,L;\mathbb Z)$ be a fixed homology class satisfying $\p A=\sum \mu_i.$ We have the same results.
\v

\subsection{Holomorphic blocks in $M^+$}\label{holomorphic_block_M}

Let $((\Sigma,{\bf j}); {\bf y}, {\bf p} )$ be a connected semistable curve with $m$ marked points
${\bf y}=(y_1,...,y_m)$ and $\nu $ puncture points ${\bf p} =(p_1 ,...,p_{\nu } )$.
Let  $u:{\Sigma} \rightarrow  M^{+}$ be a $J$-holomorphic map. Suppose that $u(z)$ converges to a $k_i \cdot T_{\fe_{i} }$-periodic orbit $x(k_iT_{\fe_{i}}t)\subset \mathcal{F}_{\fe_{i} }$ as $z$ tends to $p_i $.

\begin{definition}  A  J-holomorphic map $(u;((\Sigma,{\bf j}),{\bf y},{\bf p}))$ is said to be stable if for each $v$ one of the following conditions holds:
\begin{itemize}
\item[(1).] $u\circ \pi_{\Sigma_{v}}:\Sigma_{v}\rightarrow M^{+}$ is not a constant map.
\item[(2).] Let $val_v$ be the number of special points on $\Sigma_{v}$ which are nodal points, marked points or puncture points. Then
$val_v + 2g_{v}\geq 3.$
\end{itemize}
\end{definition}

\v

\begin{definition}\label{holomorphic_block_map_equiv_M}
Two   stable $J$-holomorphic maps $\Gamma =(u, (\Sigma,{\bf j}),{\bf y},{\bf p})$ and $\check{\Gamma} =(\check{u}, (\check{\Sigma},\check{\bf j}),{\bf \check{y}},{\bf \check{p}})$
is called equivalent if there exists a diffeomorphism $\varphi:\Sigma\rightarrow \check{\Sigma}$ such that  it can be lifted to bi-holomorphic isomorphisms $\varphi_{v w}:(\Sigma_{v},j_v)\rightarrow (\check\Sigma_{w},\check j_w)$ for
 each component $\Sigma_{v}$ of $\Sigma$, and
\begin{itemize}
\item[{\bf(1)}] $\varphi(y_i)= \check{y}_i$,  $\varphi(p_j)= \check{p_j}$ for any $1\leq i\leq m$, $1\leq j\leq \nu$,
\item[{\bf(2)}]  $\check{u}\circ \varphi= u$.
\item[{\bf(3)}] near every periodic orbit $x$, $\tilde{u}$ and $\check{\tilde{u}}\circ \varphi$ may differ by a canonical coordinate transformation \eqref{canonical}.
\end{itemize}
\end{definition}

\v
\begin{definition}  Put
\begin{align*} Aut(u;((\Sigma,{\bf j}),{\bf y},{\bf p}))=
\{\phi:\Sigma\rightarrow \Sigma| &\phi \mbox{ is an automorphism satisfying (1), (2)}  \\ &\mbox{ and (3) in Definition \ref{holomorphic_block_map_equiv_M}}  \} .\end{align*}
We call it the automorphism group of $(u;((\Sigma,{\bf j}),{\bf y}, {\bf p}))$.
\end{definition}
\v

The following Lemma is obvious.

\begin{lemma}
A $J$-holomorphic map $(u;((\Sigma,{\bf j}),{\bf y},{\bf p}))$ is stable if and only if $Aut(u;((\Sigma,{\bf j}),{\bf y},{\bf p}))$ is a finite group.
\end{lemma}

Denote by $\mathcal{M}_{A}(M^{+},C,g,m+\nu,{\bf y},{\bf p},(\mathbf{k},\mathfrak{e}))$ the moduli space of equivalence classes of all $J$-holomorphic curves in $M^{+}$ representing the homology class $A$ and satisfying $(\mathbf{k},\mathfrak{e})$ boundary condition.
\v
Fix $A\in H^{2}(M,L;\mathbb Z),$ denote by $\mathcal{M}_{A}(M,L;g,m+\nu,{\bf y},{\bf p},\overrightarrow{\mu})$ the moduli space of equivalence classes of all $J$-holomorphic curves in $M^{+}$ representing the homology class $A$ and satisfying $(\overrightarrow{\mu})$ boundary condition.

\begin{lemma}
There is a constant $C>0$ depending   on $A$ and $\mathbf{k}$ such that for any $b=(u;(\Sigma,j),{\bf y},{\bf p})\in \mathcal{M}_{A}(M^{+},C,g,m+\nu,{\bf y},{\bf p},(\mathbf{k},\mathfrak{e}))$ we have, over the cylinder end,
\begin{equation}
E_{\phi}(u)+\int_{\Sigma}u^*d\lambda \leq C.
\end{equation}
\end{lemma}
{\bf Proof.} By the Stokes formula we have  $\int_{\Sigma}u^*d\lambda \leq \sum_{i=1}^{\nu} k_{i}\cdot T_{\fe_{i}} $ over the cylinder end. Then the lemma follows from \eqref{energy_bound} \;\;\;\; $\Box$

Similarly, we also have
\begin{lemma}
There is a constant $C>0$ depending   on $A$ and $ \overrightarrow{\mu}$ such that for any $b=(u;(\Sigma,j),{\bf y},{\bf p})\in \mathcal{M}_{A}(M,L;g,m+\nu,{\bf y},{\bf p},\overrightarrow{\mu}) $ we have, over cylinder end,
\begin{equation}
E_{\phi}(u)+\int_{\Sigma}u^*d\lambda \leq C.
\end{equation}
\end{lemma}

We call $\mathcal{M}_{A}(M^{+},C,g,m+\nu,{\bf y},{\bf p},(\mathbf{k},\mathfrak{e}))$   a holomorphic block in $M^+$.
\v
Let $b=(u;(\Sigma,j),{\bf y},{\bf p})\in \mathcal{M}_{A}(M^{+},C,g,m+\nu,{\bf y},{\bf p},(\mathbf{k},\mathfrak{e}))$, $D_{u}:\mathcal{W}^{1,p,\alpha}\rightarrow L^{p,\alpha}$ be a Fredholm operator with $ind=dim (ker D_{u})-dim(coker D_{u})$. Put
$$Ind^C=ind +6(g-6)+2(m+\nu).$$
The virtual dimension of $\mathcal{M}_{A}(M^{+},C,g,m+\nu,{\bf y},{\bf p},(\mathbf{k},\mathfrak{e}))$ is $Ind^C$.
\v
Similarly, let $b=(u;(\Sigma,j),{\bf y},{\bf p})\in \mathcal{M}_{A}(M,L;g,m+\nu,{\bf y},{\bf p},\overrightarrow{\mu})$. Suppose that $\lim\limits_{s_j\rightarrow \infty} u(s_j,S^1)\subset \mathcal{F}_{i_{j}}$ and
 \begin{equation}\label{Le_rest}
[\pi(\lim\limits_{s_j\rightarrow \infty} u(s_j,S^1))]=\mu_j,\;\;\forall\;1\leq j \leq \nu.
\end{equation}
Then $D_{u}:\mathcal{W}^{1,p,\alpha}\rightarrow L^{p,\alpha}$ be a Fredholm operator with $ind=dim (ker D_{u})-dim(coker D_{u})$.
Put
$$Ind^L=ind +6(g-6)+2(m+\nu).$$
The virtual dimension of $\mathcal{M}_{A}(M,L;g,m+\nu,{\bf y},{\bf p},\overrightarrow{\mu})$ is $Ind^L$.
\vskip 0.1in
\noindent
It is possible that there are finite many combinations such that \eqref{Le_rest} holds. In this case $\mathcal{M}_{A}(M,L;g,m+\nu,{\bf y},{\bf p},\overrightarrow{\mu})$ is the union of some $\mathcal{M}_{A}(M^{+},C,g,m+\nu,{\bf y},{\bf p},.)$.
\v
\begin{remark} As we mod the $S^1$ action on every periodic orbit of Reeb vector field, the situation is very similar to the symplectic cutting. For example, let $b=(u; \Sigma,j, p)\in {\mathcal{M}}_{A}(M^{+},C;g,m+1,(k,\fe_{1}))$. Roughly speaking :
\v We collapse the $S^1$-action on the orbit $x(kT_{\fe_{1}}t)$
at infinity to get a ``manifold" $\widehat{M}^+$, locally ,  such that $Z_{\fe_{1}}$ is a ``submanifold" of $\widehat{M}^+$, choose a local pseudo-Daubaux coordinate $(a,\vartheta,{\bf w})$ around $x(kT_{\fe_{1}}t)$. Our estimates \eqref{exponential_decay_a},\eqref{exponential_decay_theta},
 \eqref{exponential_decay_y}
show that the puncture point $p$ can be ``removed", we get a $J$-holomorphic map $\bar{u}$ from $\Sigma $ into $\widehat{M}^+$. The condition
that $u$ converges
to a $k$-multiple periodic orbit at $p$ is naturally interpreted as
$\bar{u}$ being tangent to $Z_{\fe_{1}}$ at $p$ with order $k$.
\v
We will study the local geometry in other paper.
\end{remark}

\subsection{Holomorphic blocks in ${\mathbb{R}}\times \widetilde{M}$}\label{holomorphic_block_R}

Note that the space $\mathcal{M}_{A}(M^{+},C;g,m+\nu,({\bf k},\mathfrak{e}))$ is not large enough to compactify the moduli space of all $J$-holomorphic maps into $M^+$, we need consider ${\mathcal{M}}_{A}({\mathbb{R}}\times \widetilde{M},C; ({\bf k}^-,\fe^{-}),({\bf k}^+,\fe^{+})) $, which will be studied in this section.

Let $((\Sigma,{\bf j}); {\bf y}, {\bf p}^{+},{\bf p}^{-})$ be a connected semistable curve with $m$ marked points
${\bf y}=(y_1,...,y_m)$ and $\nu^{\pm}$ puncture points ${\bf p}^{+}=(p_1^{+},...,p_{\nu^{+}}^{+})$, ${\bf p}^{-}=(p_1^{-},...,p_{\nu^{-}}^{-})$,
and $u:{\Sigma} \rightarrow {\mathbb{R}}\times \widetilde{M}$ be a $J$-holomorphic map. Suppose that $u(z)$ converges to a $k_i^{\pm}\cdot T_{\fe_{i}^{\pm}}$-periodic orbit $x_{k_i^{\pm}}\subset \mathcal{F}_{\fe_{i}^{\pm}}$ as $z$ tends to $p_i^{\pm}$.

\v

Similar to the situation of relative invariants, there is a $\mathbb R$ action, which induces a $\mathbb R$-action on the moduli space of $J$-holomorphic maps. We need mod this action. Since there is a vector field $X$ with $|X|=1$ on $\widetilde{M}$, the Reeb vector field, there is a one parameter group $\varphi_\theta$ action on $\widetilde{M}$ generated by $-X$. In particular, there is a $S^1$-action on every periodic orbit, corresponding to the freedom of the choosing origin of $S^1$. Along every periodic orbit we have $L_X \lambda =0$. Recall that by (3) of {\bf Assumption} \ref{bott_assumption}, $L_X J=0$ along every periodic orbit. We can mod this action.
\begin{definition}\label{equivalent_R}
Two $J$-holomorphic maps $\Gamma =(u, (\Sigma,{\bf j}),{\bf y}, {\bf p}^{+},{\bf p}^{-})$ and $\check{\Gamma} =(\check{u}, (\check{\Sigma},\check{\bf j}),{\bf \check{y}}, \check{\bf p}^{+},\check{\bf p}^{-})$
are called equivalent if there exists a diffeomorphism $\varphi:\Sigma\rightarrow \check{\Sigma}$ such that
  it can be lifted to bi-holomorphic isomorphisms $\varphi_{v \check{v}}:(\Sigma_{v},j_v)\rightarrow (\check\Sigma_{\check v},j_{\check v})$ for
 each component $\Sigma_{v}$ of $\Sigma$, and
\begin{itemize}
\item[{\bf(1)}] $\varphi(y_i)= \check{y_i}$,  $\varphi(p_j^{+})= \check{p_j}^{+}$, $\varphi(p_l^{-})= \check{p_l}^{-}$ for any $1\leq i\leq m$, $1\leq j\leq \nu^{+}$, $1\leq l\leq \nu^{-};$  $u$ and $\check u\circ \varphi$ converges to the same periodic orbit $x_{k_{i}^{\pm}}$ at $z$ tends to $p_{i}^{\pm};$
\item[{\bf(2)}] $\check{a}\circ \varphi= a+C,$ $\check{\tilde{u}}\circ \varphi=\tilde{u}$ for some constant $C$;
\item[{\bf(3)}] near every periodic orbit $x$, $\tilde{u}$ and $\check{\tilde{u}}\circ \varphi$ may differ by a canonical coordinate transformation \eqref{canonical}.
\end{itemize}
\end{definition}

\begin{definition}  A  J-holomorphic map $(u;((\Sigma,{\bf j}),{\bf y},{\bf p}))$ is said to be stable if for each $v$ one of the following conditions holds:
\begin{itemize}
\item[(1).] $\tilde{E}(u\circ \pi_{\Sigma_{v}})\ne 0$,
\item[(2).] Let $val_v$ be the number of special points on $\Sigma_{v}$ which are nodal points, marked points or puncture points. Then
$val_v + 2g_{v}\geq 3.$
\end{itemize}
\end{definition}

\v

\begin{definition}  Put \begin{align*}Aut(u;((\Sigma,{\bf j}),{\bf y},\check{\bf p}^{+},\check{\bf p}^{-}))=\{\phi:\Sigma\rightarrow \Sigma| &\phi \mbox{ is an automorphism satisfying (1), (2)}  \\ &\mbox{ and (3) in Definition \ref{equivalent_R}}  \} .\end{align*} We call it the automorphism group of $(u;((\Sigma,{\bf j}),{\bf y}, \check{\bf p}^{+},\check{\bf p}^{-}))$.
\end{definition}

For any $A\in H^{2}(\bR\times \widetilde M,\Re;\mathbb Z)$ we define $d\lambda (A)$ as following:
let $v: \mathbb R \times S^1\rightarrow \bR\times \widetilde M$ be a $C^{\infty}$ map such that $[v(\mathbb R \times S^1)]=A$, we define $d\lambda (A):= \int_{\mathbb R \times S^1} v^*(d\lambda).$

We fix $A\in H^{2}(\bR\times \widetilde M,\Re;\mathbb Z)$ and
${\bk}^\pm=(k_1^{\pm},...,k_{\nu^{\pm}}^{\pm})$, $\fe^{\pm}=(\fe_1^{\pm},...,\fe_{\nu^{\pm}}^{\pm})$
satisfying
\begin{equation}\label{energy_condtion_RM}
d\lambda(A)=\sum_{i=1}^{\nu^{+}} k_i^{+}\cdot T_{\fe_{i}^{+}}-\sum_{i=1}^{\nu^{-}} k_i^{-}\cdot T_{\fe_{i}^{-}} .
\end{equation}

We define
${\mathcal{M}}_{A}({\mathbb{R}}\times \widetilde{M},g,m+\nu^{+}+\nu^{-},({\bf k}^-,\fe^{-}), ({\bf k}^+,\fe^{+}))$ to be the space of equivalence classes
of all stable $J$-holomorphic maps in ${\mathbb{R}}\times \widetilde{M}$ representing $A$ and converging to a $k_i^{\pm}\cdot T_{\fe_{i}^{\pm}}$-periodic orbit $x_{k_{i}^{\pm}}\subset \mathcal{F}_{\mathfrak{e}_{i}^{\pm}}$ as $z$ tends to $p_i^{\pm}$. For any $(u,\Sigma,{\bf y},{\bf p}^{+},{\bf p}^{-})\in {\mathcal{M}}_{A}({\mathbb{R}}\times \widetilde{M},g,m+\nu^{+}+\nu^{-},({\bf k}^-,\fe^{-}), ({\bf k}^+,\fe^{+})),$ by Stoke's formula we have
$$d\lambda(A)=\int_{\Sigma}u^*d\lambda=\sum_{i=1}^{\nu^{+}} \lambda(k_{i}^{+}x_{i}^{+}) -\sum_{i=1}^{\nu^{-}}\lambda(k_{i}^{-}x_{i}^{-}) =\sum_{i=1}^{\nu^{+}} k_i^{+}\cdot T_{\fe_{i}^{+}}-\sum_{i=1}^{\nu^{-}} k_i^{-}\cdot T_{\fe_{i}^{-}} .$$

We call ${\mathcal{M}}_{A}({\mathbb{R}}\times \widetilde{M},C;g,m+\nu^{+}+\nu^{-},({\bf k}^-,\fe^{-}), ({\bf k}^+,\fe^{+}))$ a holomorphic rubber block in ${\mathbb{R}}\times \widetilde{M}$.
\begin{remark}
Let $\Sigma_v\subset \Sigma$ be a smooth connected component of $\Sigma$. Suppose that $u(\Sigma_v)$ lie in a compact set $K\subset {\mathbb{R}}\times \widetilde{M}$. Then every nodal point of $\Sigma_v$ is a removable singular point. As  $\omega_{\phi}=d(\phi \lambda)$ we conclude that $E_{\phi}(u)\mid_{\Sigma_v}=0$, then $u(\Sigma_v)$ is one point, then
$\widetilde{E}(u)\mid_{\Sigma_v}=0$.
\end{remark}

\subsection{Local coordinate system of holomorphic blocks}

Let  $b=(u,\Sigma,{\bf y}, {\bf p}^{+},{\bf p}^{-})\in {\mathcal{M}}_{A}({\mathbb{R}}\times \widetilde{M},C;g,m+\nu^{+}+\nu^{-},({\bf k}^-,\fe^{-}),({\bf k}^+,\fe^{+})) $, and
$h\in {\mathcal W}^{1,p,\alpha}(\Sigma, u^*T( {\mathbb{R}}\times \widetilde{M}))$ such that $D_u(h)=0$. Under the  transformation
\begin{equation}\label{a}
\hat a=a+C\end{equation}
and the canonical transformation of local pseudo-Daubaux coordinate system
\begin{equation}\label{a}
\hat \vartheta =\vartheta + \vartheta_0\end{equation}
$h$ is invariant ( see subsection \S\ref{Lower strata}, Case 2 below), so we can view $h$ as a element
in $T_b({\mathcal{M}}_{A}(\bR\times \widetilde{M},C;m+\nu^{+}+\nu^{-},({\bf k}^-,\fe^{-}),({\bf k}^+,\fe^{+}) )$. We assume that $D_u$ is surjective. Let $e_1,...,e_a$ be a base of $ker D_u$. Then for any $h\in ker D_u$ we have $h=\sum_i^a x_i e_i.$ Let $w_1,...,w_d$ be the local coordinates of ${\mathcal{M}}_{g,m+\nu^{+}+\nu^{-}}$ around $\Sigma$, where $d=6g-6+2(m+\nu^{+}+\nu^{-})$. Then the  $(w_1,...,w_d,x_1,...,x_a)$ is a local coordinate system of ${\mathcal{M}}_{A}(\bR\times \widetilde{M},C;m+\nu^{+}+\nu^{-},({\bf k}^-,\fe^{-}),({\bf k}^+,\fe^{+})) $.

Similarly, let $b=(u,\Sigma,{\bf y}, {\bf p}^{+})\in {\mathcal{M}}_{A}(M^{+},C;g,m+\nu,({\bf k},\fe))$,
$h\in {\mathcal W}^{1,p,\alpha}(\Sigma, u^*T(M^+))$ with $D_u(h)=0$, we view $h$ as a element
in $T_b({\mathcal{M}}_{A}(M^{+},C;g,m+\nu,({\bf k},\fe))$. Assume that $D_u$ is surjective. We can choose a local coordinate system in ${\mathcal{M}}_{g,m+\nu}$ around $\Sigma$ together with a coordinate system in $ker D_u$ as a local coordinate system of ${\mathcal{M}}_{A}(M^{+},C;g,m+\nu,({\bf k},\fe))$. \v

\section{Compactness theorems }\label{compact_theorem}
We will use Li-Ruan's Compactification to our setting.
Roughly speaking, the idea of Li-Ruan's Compactification is that firstly let the Riemann surfaces degenerate in Delingne-Mumford space and then let $M^+$ degenerate compatibly. We will give a detail description of the Li-Ruan's idea and give the definition of the stable compactification using language of graphs. This is one of the key parts of this paper.
\v

\subsection{Bubble phenomenon }\label{bubble_phenomenon}
\vskip 0.1in
\noindent
The proofs of the following  two theorems are standard (see \cite{BEHWZ}).
\begin{theorem}\label{energy_lower_bound-1}
There is a constant $\hbar>0$ such that
for every $J$-holomorphic map $u=(a,\widetilde{u}):
{\bC}\rightarrow \bR\times \widetilde{M}   $ with finite energy,
\begin{itemize}
\item[\bf(1)] if \;$ E_{\phi}(u)\neq 0 ,$ we have $E_{\phi}(u)\geq \hbar$;
\item[\bf(2)] if \;$\widetilde{E}(u)\neq 0 ,$ we have $\widetilde{E}(u)\geq \hbar$.
\end{itemize}
\end{theorem}

\begin{theorem}\label{lambda_lower_bound}
There is a constant $\hbar>0$ such that
for every $J$-holomorphic map $u\in  \mathcal{M}_A({\mathbb{R}}\times \widetilde{M},C;g,m+\nu^{\pm},({\bf k}^-,\fe^{-}), ({\bf k}^+,\fe^{+})) $ with finite energy,
  if \;$\widetilde{E}(u)\neq 0 ,$ we have $$\widetilde{E}(u)=\sum_{i=1}^{\nu^+} k_{i}^{+}-\sum_{i=1}^{\nu^{-}}k_{i}^{-} \geq \hbar.$$
\end{theorem}

\v

Let $\Gamma\upper=(u\upper,\Sigma\upper;{\bf y}\upper,{\bf p}\upper) \in {\mathcal M}_{A}(M^{+},C,g,m+\nu,A,({\bf k},\fe))$ be a sequence. Then there is a constant $C>0$ such that
$$\widetilde{E}(u\upper)+E_{\phi}(u\upper)<C\;\;\;\forall\;i.$$

Suppose that $(\Sigma\upper; {\bf y}\upper,{\bf p}\upper)$ is stable and
converges to $(\Sigma; {\bf y}, {\bf p})$ in ${\overline{\mathcal{M}}}_{g,m+\nu }$. On a neighborhood of $(\Sigma;{\bf y}, {\bf p})$ in
${\overline{\mathcal{M}}}_{g,m+\nu}$ we construct a smooth family of metrics on each
$({\Sigma}\upper;{\bf y}\upper, {\bf p}\upper)$
such that near each marked point the metric is the Euclidean metric on disc in $\mathbb C$, and near each nodal point the metric is
the standard cylinder metric.
\v
\subsubsection{Bound of the number of singular points }\label{number of singular points}

 Following McDuff and Salamon \cite{MS} we introduce the notion of
singular points for a sequence $u^{(i)}$ and the notion of mass of singular
points. We show that there is a constant $\hbar >0$ such that the mass of every singular
point is large than $\hbar$. Let $q$ be a singular point and $q^{(i)}\in \Sigma\upper$, $q^{(i)}\rightarrow q$. In case $u^{(i)}(q^{(i)})\in M^+_0$
the argument is standard (see \cite{MS}). We only consider $\widetilde{E}$ over the cylinder end. Since the metrics near nodal points are the standard cylinder metrics, we have, in the cylinder coordinates
$D_{1}(q^{(i)})\subset \Sigma^{(i)}-\{nodal\; points\},$ where $D_{1}(q^{(i)})=\{(s\upper, t\upper)\;|\; (s\upper -s\upper(q\upper))^2+(t\upper -t\upper(q\upper))^2 \leq 1\}$. We identify $q\upper$ with $0$ and consider $J$-holomorphic maps $u^{(i)}: D_1(0)\to N.$
\v
The proof of the following  lemma is similar to Theorem 4.6.1 in \cite{MS}. We give the proof here for the reader's convenience.
\v\n
\begin{lemma}\label{lower_bound_of_singular_points}
 Let $u^{(i)}:D_1(0)\to \mathbb{R}\times \widetilde{M}$ be a sequence of J-holomorphic maps with finite energy such that
  $$\sup_{i}\widetilde{E}(u\upper)<\infty,\;\;\;\;\; |du^{(i)}(0)|\longrightarrow \infty ,\;\; as \;\;i\to \infty.$$
Then there is a constant $\hbar >0$ independent of $u^{(i)}$ such that,  for every $\epsilon>0$
\begin{equation}
 \liminf\limits_{i\to\infty}\widetilde{E}(u^{(i)};D_{\epsilon}(0)) \geq \hbar.
\end{equation}
\end{lemma}
\noindent
{\bf Proof: } Consider the function $$F\upper(z)=|du^{(i)}|(z)d^{2}(z,\partial D_{1}(0))$$
where $d(z,\partial D_{1}(0))$ denotes the distance from $z$ to $\partial D_{1}(0)$ with respect to the standard Eucildean metrics. Obviously, $F\upper$ attains its maximum at some interior point $q^\star_{i}\in  D_{1}(0)$ and
$\lim\limits_{i\to\infty}F\upper(q_{i}^\star)=\infty.$ Set $\delta_{i}=\frac{1}{2}d(q^\star_{i},\p D_{1}(0)).$ Then for any $q\in D_{\delta_{i}}(q^\star_{i}),$
$$|du^{(i)}|(q)\leq 4 |du^{(i)}|(q^\star_{i}):= 4A _{i}.$$
 Consider the re-scaling sequence  $$v^{(i)}(z)=u^{(i)}\left(q^\star_{i}+\frac{z}{A_{i}}\right ) $$
  As $\lim\limits_{i\to\infty}F_{i}(q_{i}^\star)=\infty$, we have $\lim\limits_{i\to\infty}\delta _{i}A_{i}=\infty.$  Then in $D_{\delta_{i}A_{i}}(0)$
 $$\sup|dv^{(i)}|\leq 4,\;\;\; |dv^{(i)}|(0)=1,\;\;\;\;\widetilde{E}(v^{(i)},D_{\delta_{i}A_{i}}(0))=\widetilde{E}(u^{(i)},D_{\delta_{i}}(q^\star_{i})).$$
 By choosing a subsequence we conclude that $v_{i}$ locally uniformly converges to a nonconstant $J$-holomorphic map with finite energy
 $v:\mathbb C\longrightarrow \mathbb{R}\times \widetilde{M}$. Then the lemma follows from Lemma \ref{energy_lower_bound-1}.    $\;\;\;\;\Box$
\v


\v
Recall that $N$ denotes one of $M^+$ and $\mathbb R\times \tilde{M}$.
By Lemma \ref{lower_bound_of_singular_points} we conclude that the rigid singular points are isolated and the limit
$$m_{\epsilon}(q)=\lim_{i \rightarrow \infty} \widetilde{E}(u^{(i)};D_{q^{(i)}}
(\epsilon,h^{(i)}))$$ exists for every sufficiently small $\epsilon>0$.
The mass of the singular point q is defined to be
$$m(q)=\lim_{\epsilon\rightarrow 0}{m_{\epsilon}(q)}.$$

\v
Denote by $P\subset \Sigma $ the set of  singular points for $u^{(i)}$,
the double points and the puncture points. By Lemma \ref{lower_bound_of_singular_points} and \eqref{energy_bound}, $P$ is
a finite set. By definition, $|du^{(i)}|_{h^{(i)}}$ is uniformly bounded on every
compact subset of $\Sigma - P$.   By a possible translation along $\mathbb R$ and passing to a  subsequence we may
assume that $u^{(i)}$ converges uniformly  with all derivatives on every compact
subset of $\Sigma - P $ to a $J$-holomorphic map $u:\Sigma - P
\rightarrow N.$ Obviously, $u$ is a finite energy $J$-holomorphic map.
\v\n

We need to study the behaviour of the sequence $u\upper$ near each
singular point for $u\upper$. Let $q \in \Sigma $ be a rigid singular point for
$u\upper$. We have two cases.
\v\n
{\bf (a)} $q\in \Sigma -\{nodal\; points\}$. We consider $J$-holomorphic maps $u^{(i)}: D_1(0)\to N$.
\begin{enumerate}
\item[(a-1)] there are $\epsilon > 0$ and a compact set $K\subset N$ such that $u^{(i)}(D_{\epsilon}(q))\subset
K$.
\item[(a-2)] $q$ is a nonremovable singularity.
\end{enumerate}
\v\n
{\bf (b)} $q\in \{nodal\; points\}$. In this case a neighborhood of $q$ is two discs $D_1(0)$ joint at $0$, where $D_1(0)=\{|z|^2\leq 1\}$.
\v\n

For (a-1) we  construct bubbles as usual for a compact symplectic manifold (see
\cite{RT,PW,MS}). We call this type of bubbles (resp. bubble tree) the normal bubbles (resp. normal bubble tree).
\v\n
 \subsubsection{Construction of the bubble tree for (a-2)}\label{bubble_tree_a}
\v\n
We use cylindrical coordinates $z=e^{-s-2 \pi \sqrt{-1}t}$
and write
$$u\upper(s,t)=(a\upper(s,t), \widetilde{u}\upper(s,t))$$
$$u(s,t)=(a(s,t),\widetilde{u}(s,t)).$$

\v
\v
Note that the graduate $|du^{(i)}|$ depends not only on the metric $<, >$ on $N$ but also depends on the metric on $\Sigma^{(i)}$.  The energy don't depend on the metric on $\Sigma^{(i)}$. To construct bubbling in present case it is more
convenient to take the Euclidean metric $|dz|^2$ on the disk $D_1(0)$ and pullback to the
coordinate $(s,t)$ through $z=e^{-s-2 \pi \sqrt{-1}t}.$ Obviously, if $q$ is a singular point with respect to the cylinder metric, then it is also a singular point with respect to the disk metric.
\v

By Theorem \ref{exponential_estimates_theorem}, we have
$$\lim_{s \rightarrow \infty} \widetilde{u}(s,t)=x(kTt)$$
in $C^{\infty}(S^1)$, where $x(\;,\; )$ is a $kT$-periodic orbit on $\widetilde{M}$.
Choosing $\epsilon $ small enough we have $$|{m_{\epsilon}(q)}-m(q)|\leq \frac{1}{10}\min\{\hbar,T\}$$
For every $i$ there exists $\delta_i>0$ such that
\begin{equation}\label{choose_delta}
\widetilde{E} (u\upper; D_{\delta_i}(0))  =m(q) - \frac{1}{2}\min\{\hbar,T\},
 \end{equation}
 where $T$ is the minimal positive periodic.
By definition of the mass $m(q)$, the sequence $\delta_i$ converges to $0$.
Since $u^{(i)}$ converges uniformly with all derivatives to $u$ on any compact set
of $D_{\epsilon}(0)-\{0\}$, $\delta_i $ must converge to 0.
Put
\begin{equation}\hat{s}\upper= s + log \delta_i,\;\;\;\hat{t}\upper= t\;\;for\;|s_i\upper|>2R_{0}.\end{equation}
\begin{equation}\label{a_rescaling_DT}
\hat{a}\upper= a -kT(-\log\delta_i).
\end{equation}
 Define the $J$-holomorphic curve $v\upper(\hat{s},t)$ by
\begin{equation}\label{bubble construction}
v\upper(\hat{s},t)=(\hat{a}\upper(\hat{s},t),\widetilde{v}\upper(\hat{s},t))=
\left(a\upper(-\log\delta_{i}+ \hat{s},t) -kT(-\log\delta_i) ,
\widetilde{u}\upper(-\log\delta_{i}+ \hat{s},t)\right).
\end{equation}

\begin{lemma}\label{compact_lemma}
Suppose that 0 is a nonremovable singular point of $u$.
Define the $J$-holomorphic map $v\upper$ as above. Then there exists a
subsequence (still denoted by $v\upper$) such that
\begin{description}
\item[(1)] The set of singular points $\{Q_1,\cdot \cdot \cdot,Q_d\}$ for
$v\upper$ is
finite and tame, and is contained in the disc $D_{1}(0)=\{z \mid \mid z
\mid \leq 1\};$
\item[(2)] The subsequence $v^{(i)}$ converges with all derivatives uniformly on
every compact subset of  ${\mathbb C}\backslash\{Q_1,\cdot \cdot \cdot,Q_d\}$
to a nonconstant
J-holomorphic map $v:{\mathbb C}\backslash\{Q_1,\cdot \cdot \cdot,Q_d\}
\rightarrow {\mathbb R}\times \widetilde{M};$
\item[(3)]  $\widetilde{E}(v)>\frac{1}{3} \min\{\hbar ,T \}$;
     \item[(4)] $\widetilde{E}(v)+\sum\limits_{1}^{d}m(Q_i)=m(0). $
\item[(5)] $\lim\limits_{s \rightarrow \infty} \widetilde{u}(s,t)=
\lim\limits_{\hat{s} \rightarrow -\infty} \widetilde{v}(\hat{s},t).$
\end{description}
\end{lemma}
\vskip 0.1in
\noindent
{\bf Proof: } The proofs of
{\bf (1)}, {\bf (2)} and {\bf (4)} are standard (see \cite{MS}), we omit them here. The proof of {\bf (5)} will be given in our next paper (see \cite{LS1}). We only prove {\bf (3)}.
 \vskip
0.1in \noindent
 {\bf(3)} Note that
\begin{equation}
\left|\int_{\widetilde{u}\upper(-\log \delta_{i},S^{1} )}  \lambda - \int_{\widetilde{u}\upper(-\log\epsilon,S^{1} )}  \lambda \right|=\widetilde{E}(u; -\log\epsilon\leq s \leq - \log \delta_{i}) =\frac{1}{2}\min\{\hbar, T\}.
\end{equation}
Since $\lim\limits_{s\to \infty}\widetilde{u}(s,t)=x(kTt)$ we have
$$\left|\int_{\widetilde{u}\upper(-\log\epsilon,S^{1} )}\lambda-kT\right|\leq \frac{1}{8} \min\{\hbar, T\}. $$
 when $\epsilon$ small enough and $i$ big enough.
Then by the locally uniform convergence of $v\upper$
\begin{equation}
 \widetilde{E}( v)\geq \left|kT-\int_{\widetilde{v} (0,S^{1} )}  \lambda \right| \geq \frac{1}{3}\min\{\hbar, T\}.
\end{equation}
 $\Box$
\v

We can repeat this again to construct bubble tree.

\v\n
We introduce a terminology. Let $\Sigma_1$ and $\Sigma_2$
join at $p$, and $u_i: \Sigma_i \rightarrow
\bR\times \widetilde{M}$ a map. Choose holomorphic cylindrical
coordinates $z_1=(s_1,t_1)$ on $\Sigma_1$ and $z_2=(s_2,t_2)$
on $\Sigma_2$ near $p$ respectively. Suppose that
$$\lim_{s_1 \rightarrow -\infty} \widetilde{u_1}(s_1,t_1)= x_1(k_1T_1t_1)$$
$$\lim_{s_2 \rightarrow +\infty} \widetilde{u_2}(s_2,t_2)= x_2(k_2T_2t_2).$$
We say $u_1$ and $u_2$ converge to a same periodic orbit, if $k_1=k_2$,
and in the pseudo-Daubaux coordinates $(a,\vartheta,{\bf w}),$ (see Lemma \ref{Darboux}) ${\bf w}(x_1)={\bf w}(x_2)$ (in this case $T_1=T_2$ holds naturally).
\v\n

\subsubsection{Construction of bubble tree for (b)}\label{bubble_tree_b}
\v\n

Let $\Sigma=\Sigma_1 \wedge \Sigma_2$,
where $\Sigma_1$ and $\Sigma_2$ are smooth Riemann surfaces of genus $g_1$ and $g_2$ joining at $q$.
Let $z_1$ , $z_2$ be the local complex coordinates of $\Sigma_1$ and $\Sigma_2$ with
$z_1(q)=z_2(q)=0$. Recall that a neighborhood of $\Sigma_1 \wedge \Sigma_2\in \overline{{\mathcal M}}_{g,m+\nu}$ is given by
$$z_1z_2=w=e^{-R - 2\pi\sqrt{-1} \tau},\;\;w \in \mathbb C,$$
where $R=2l r$, $l\in \mathbb Z^{+}$.
We will use $(r,\tau)$ as the local coordinates in the neighborhood of
$\Sigma\in \overline{{\mathcal M}}_{g,m+\nu}$. Let
$$z_1=e^{-s_1 - 2\pi\sqrt{-1} t_1},\;\;\;z_2=e^{s_2 +2\pi\sqrt{-1} t_2}.$$
$(s_i, t_i)$ are called the holomorphic cylindrical coordinates near $p$.  In terms of the holomorphic cylindrical
coordinates we write
$$\Sigma_1-\{p\}=\Sigma_{10}\bigcup\{[0,\infty)\times S^1\},$$
$$\Sigma_2-\{p\}=\Sigma_{20}\bigcup\{(-\infty,0]\times S^1\}.$$
\begin{eqnarray}
&& s_1=s_2 + 2l r \\
&&t_1=t_2 + \tau.
\end{eqnarray}
Then $\Sigma_{(r)}\rightarrow  \Sigma_1 \wedge \Sigma_2$ as $r\rightarrow \infty$ in the $\overline{{\mathcal M}}_{g,m+\nu}$.
\v

We consider the case that $u\upper:\Sigma\upper\rightarrow N$  is a sequence of $J$-holomorphic maps, where $\Sigma\upper:=\Sigma_{(r\upper)}$, $q^{(i)}\in \Sigma\upper$,
$q^{(i)}\rightarrow q$.

Without loss of generality we assume that there is a subsequence, still denoted by $i$, such that $|s_{1}\upper(q^{(i)})|\leq l r\upper$. In this case we may identify a neighborhood of $q^{(i)}$ in
$\Sigma\upper$ with $D_1(0)\setminus D_{\epsilon\upper}(0)\subset \Sigma_1$, where $\lim\limits_{i}\epsilon\upper=0.$  The sequence
$u^{(i)}$ is considered to be a sequence of $J$-holomorphic maps from $D_1(0)\setminus D_{\epsilon\upper}(0)$
into $N$, and $q^{(i)}\in D_1(0)\setminus D_{\epsilon\upper}(0)$, $q^{(i)}\rightarrow 0$. In terms of the cylinder coordinates we can identify the coordinates $(s\upper,t\upper)$ in $\Sigma\upper$ with the coordinates $(s_1,t_1)$ of $\Sigma_1$.

By Theorem \ref{exponential_estimates_theorem}, we have
$$\lim_{s_1 \rightarrow \infty} \widetilde{u}(s_1,t_1)=x(kTt_1)$$
in $C^{\infty}(S^1)$, where $x(\;,\; )$ is a $kT$-periodic orbit on $\widetilde{M}$.
Choosing $\epsilon $ small enough we have $$|{m_{\epsilon}(q)}-m(q)|\leq \frac{1}{10}\min\{\hbar,T\}$$
For every $i$ there exists $\delta_i>0$ such that
\begin{equation}\label{choose_delta}
\widetilde{E} (u\upper; D_{\delta_i}(0))  =m(q) - \frac{1}{2}\min\{\hbar,T\},
 \end{equation}
 where $T$ is the minimal positive periodic.
We choose a smooth conformal transformation $\psi\upper:\Sigma\upper \rightarrow   \Sigma\upper$ such that
\begin{equation}\label{D-rescaling 1}
\psi\upper|_{\mathfrak{D}(R_{0})}=I,\end{equation}
\begin{equation}\label{D-rescaling 2}
z\upper=\delta_i \hat{z}\upper = \delta_i e^{-\hat{s}\upper-2\pi \sqrt{-1} \hat{t}\upper},\;\;for\;|s\upper|>2R_{0}
\end{equation}
where $(\hat s\upper,\hat t\upper)=\psi_{i}(s\upper,t\upper)$, i.e.,
\begin{equation}\hat{s}\upper= s\upper + log \delta_i=s_1 + log \delta_i ,\;\;\;\;\;\;\; \hat{t}\upper= t\upper= t_1\;\;for\;|s_i\upper|>2R_{0}.\end{equation}
We call coordinates rescalling of Riemann surface a  $\mathcal{D}$-rescaling.
In our present case we need not only a $\mathcal{D}$-rescaling, but also a translation along $\mathbb R$, called a $\mathcal{T}$-rescaling. We call such  composition a $\mathcal{D}$$\mathcal{T}$-rescaling. Put
\begin{equation}\label{a_rescaling_DT}
\hat{a}\upper=a -kT(-\log\delta_i).
\end{equation}
 Define the $J$-holomorphic curve $v\upper(\hat{s},t)$ by
\begin{equation}\label{bubble construction}
v\upper(\hat{s},t)=(\hat{a}\upper(\hat{s},t),\widetilde{v}\upper(\hat{s},t))=
\end{equation}
$$\left(a\upper(-\log\delta_{i}+ \hat{s},t) -kT(-\log\delta_i) ,
\widetilde{u}\upper(-\log\delta_{i}+ \hat{s},t)\right).$$
By the same argument as in subsection \S\ref{bubble_tree_a} we construct a bubble $S^2$ with
$E(v)\mid_{S^2}>\frac{1}{3} \min\{\hbar ,T \}$ or $\widetilde{E}(v)\mid_{S^2}> \frac{1}{3} \min\{\hbar ,T \}$, inserted between $\Sigma_1$ and $\Sigma_2$. The same results as Lemma \ref{compact_lemma}  still hold. We can repeat this again to construct bubble tree.
\vskip 0.1in
\noindent

\v
\subsubsection{For the case of genus $0$}\label{genus 0}

Let $\Gamma\upper=(u\upper,\Sigma\upper;{\bf y}\upper,{\bf p}\upper) \in \mathcal{M}_{A}(M^{+},C,g,m+\nu,{\bf y},{\bf p},(\mathbf{k},\mathfrak{e}))$ be a sequence.
Let $\Sigma\upper = \bigcup_{v=1}^N \Sigma\upper_v$. Assume that there is one component of $\Sigma\upper$ that has genus $0$ and is unstable. Let $\Sigma_1\upper$ is such a component.
We identify $\Sigma_1\upper$ with a sphere $S^2$, and consider $u\upper: S^2\rightarrow N$. We discuss several cases:
\vskip 0.1in
\noindent
{\bf 1).} $u\upper\mid_{S^2}$ has no singular point. Then $\|\nabla u\upper\mid_{S^2}\|$ are uniformly bounded above. As
$(u\upper,S^2;{\bf y}\upper,{\bf p}\upper)$ is stable, $E(u\upper)\mid_{S^2}\geq \hbar$ or $\widetilde E(u\upper)\mid_{S^2}\geq \hbar$. Then $u\upper\mid_{S^2}$ locally uniformly converges to $u\mid_{S^2}$ with $E(u)\mid_{S^2}\geq \hbar$ or $\widetilde E(u)\mid_{S^2}\geq \hbar$, so $(u; S^2)$ is stable.
\v\n
{\bf 2).} $S^2$ joints with $\Sigma\upper_2$ at a nodal point $q$. Then the number of the special points, including the marked points and the puncture points, is $\leq 1$.
\v\n
{\bf 2a)} There is one singular point $p$ and one special point $\mathcal{X}$.
\v\n
{\bf (2a-1)} $p\neq q$ and $p\neq \mathcal{X},$
\v\n
{\bf (2a-2)} $p= \mathcal{X}$ or $p=q$
\v \n
{\bf 2b)} There are two singular points $p_{1},p_{2}$.
\v\n
{\bf (2b-1)} $p_1\neq q$ and $p_2\neq q,$
\v\n
{\bf (2b-2)} $p_1=q $ or $p_2=q$.
\v\n
{\bf 2c)} There is only one singular point $p$ and no special point.
\v\n
{\bf 3).}  $S^{2}$ joint with $\Sigma\upper_2$ and $\Sigma\upper_3$ at $q_1$ and $q_2$ respectively. There is only one singular point $p$.
\v\n
{\bf (3a)} $p\neq q_1$ and $p\neq q_2,$
\v\n
{\bf (3b)} $p=q_1 $ or $p=q_2$,
\v\n
For the cases {\bf (2a-1)}, {\bf (2b-1)}, {\bf 3a)} we construct  bubble tree at singular points as in subsection \S\ref{bubble_tree_a}  to get a stable map $(u,S^{2})$.  For the cases {\bf (2a-2)}, {\bf (2b-2)}, {\bf 2c)} and {\bf (3b)} we  forget the map $u\upper\mid_{S^{2}}$ and contract $S^2$ as a point, then we construct a bubble as in subsection subsection \S\ref{bubble_tree_b}. By $({\bf 3})$ of Lemma \ref{compact_lemma} we get a stable map $(v,S^{2})$. We can repeat the procedure to construct bubble tree.

\subsection{$\mathcal{T}$- rescaling }\label{T_rescaling}

Let $\Gamma\upper=(u\upper,\Sigma\upper;{\bf y}\upper,{\bf p}\upper) \in \mathcal{M}_{A}(M^{+},C,g,m+\nu,{\bf y},{\bf p},(\mathbf{k},\mathfrak{e}))$ be a sequence. If there is some $\Sigma_v$ of genus 0, we treat it as in subsection \ref{genus 0}. In the following we
assume that
$(\Sigma\upper;{\bf y}\upper,{\bf p}\upper)$ is stable and converges to $(\Sigma;{\bf y}, {\bf p})$ in
${\overline{\mathcal{M}}}_{g,m+\nu}$. \vskip 0.1in
\noindent
Note that we have fixed the degeneration of Riemann surfaces now, we are concern with the $\mathcal{T}$-rescaling. We explain our procedure of $\mathcal{T}$-rescaling. We discuss maps into $M^+$, for maps into $\mathbb R \times \widetilde M$ the situation are the same.

Let $\Sigma$ be a Riemann surface with $\nu$ puncture points ${\bf p}=(p_1,..., p_\nu)$. We choose holomorphic coordinates $(s_i,t_i)$ near $p_i$.
Let $u:\overset{\circ}{\Sigma} \rightarrow M^+$ be a $J$-holomorphic map. Denote $u=(a,\tilde u).$ Suppose that
$$\lim_{|s_{j}|\rightarrow \infty}\widetilde
u(s_{j},t_j)=x(k_{j}T_{\fe_{j}}t_j)$$ in $C^{\infty}(S^1)$ for some $k_{j}T_{\fe_{j}}$-periodic orbit
$x(k_{j}T_{\fe_{j}}t)$.
For each $j$ we choose a local pseudo-Darboux coordinates $(a_j,\vartheta_j, {\bf w_j})$ near $p_j$  such that $x(k_{j}T_{\fe_{j}}t_j)=(k_jT_{\fe_{j}}t_{j}+\vartheta_{j0},0)$ for some constant $\vartheta_{j}.$ Note that for any $0< l,j\leq \nu$,
$$a_l = a_j + C_{lj},$$
where $C_{lj}$ are constants.

\v
 For simplicity, we consider the case as in Figure 1, the other cases is similar.
Let $ \Sigma=\bigcup_{j=1}^{6} \Sigma_{j},$ where each $\Sigma_{j}$ is a smooth connected component. Let
$T_{j}\subset \Sigma_{j}- {P}$ be compact sets, where $P\subset \Sigma $ denotes the set of  singular points for $u^{(i)}$, the double points and the puncture points. Denote $I=\{1,\cdots,6\}.$ We write
$$u\upper(z) =(a\upper(z) , \widetilde{u}\upper(z) ).$$
\begin{figure}[H] \label{Figure2}
\centering\includegraphics[height=5cm]{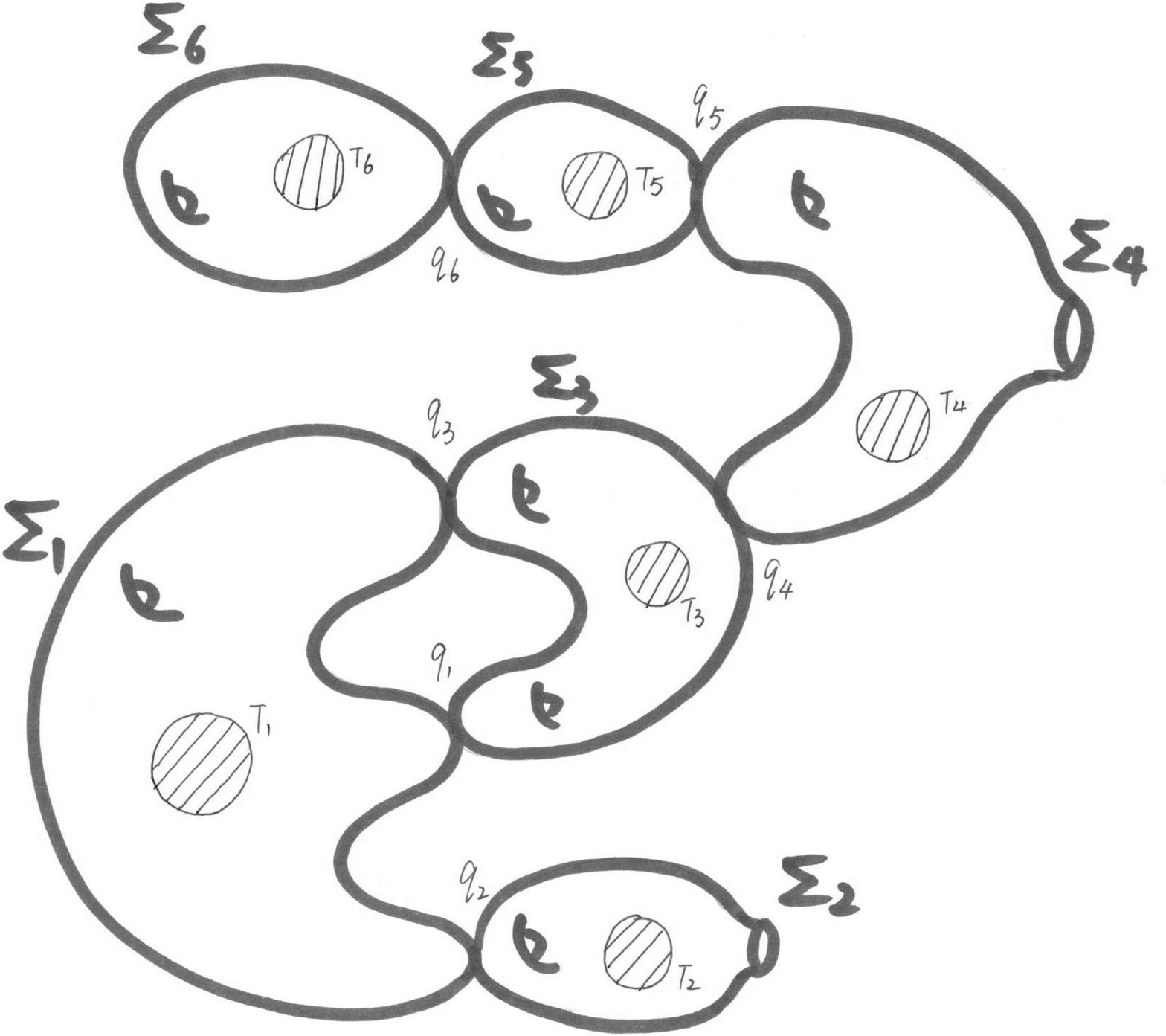}
\caption{}
\end{figure}
For each $j\in I$, $T_{j}\upper\subset\Sigma\upper$  converges to $T_{j}.$
We fix points $z_{j}\in T_{j}$. Suppose that $z_{j}\upper\in T_{j}\upper$ is a sequence of points which converges to $z_{j}$.
Without loss of generality, we assume that
\begin{itemize}
\item[(1)]  $\left|a\upper(z_{1}\upper)\right|\leq \min_{j\in I\setminus\{1\}}\left|a\upper(z_{j}\upper)\right|,\;\;\forall i,$
\item[(2)]   $\sup\limits_{i}\left|a\upper(z_{1}\upper)\right|<\infty$,
\end{itemize}
\v
Restricting to $T_1,$ $u\upper$ uniformly converges to a map $u_{1}.$  For any compact $K\subset \Sigma_1-P,$ $u\upper$
  also  converges to $u_{1}.$ Then $u$ is naturally defined on $\Sigma_1-P.$
\v
For $\Sigma_2$ we consider two different cases:
\v
{\bf Case 1.}  $\sup_{i}\left|a\upper(z_{2}\upper)\right|<\infty$.   Restricting to $T_2,$ $u\upper$ uniformly converges to a map $u_{2}.$ Then $u_{2}$ is naturally defined on $\Sigma_2-{P}.$ Then $(u_{1};\Sigma_1,{\bf y}_{1},{\bf p}_{1})$ and $(u_{2};\Sigma_2,{\bf y}_{2},{\bf p}_{2})$ are belong to the same holomorphic block in $M^{+}.$
\v
{\bf Case 2.} $\lim\limits_{k\to\infty}a\upper(z_{2}\upper)=\infty.$ In this case we take a coordinate transformation:
\begin{equation}\label{T_rescaling-1}
 (a^{*})\upper =a  - a\upper(z_{2}\upper).
\end{equation}
Define  \begin{equation}\label{T_rescaling-2}
 (u^{*})\upper(z)=((a^{*})\upper(z),  ({\tilde u}^{*})\upper(z)),\;\;\;\;\;\;\;({\tilde u}^{*})\upper(z)=\tilde u\upper(z).
\end{equation}
Note that $|du\upper|$ is invariant under translation along ${\mathbb{R}}$.
As above, restricting to $T\upper_2$, $(u^{*})\upper$ locally uniformly converges to a  map
$u^{*}=(a^{*},\tilde{u}^{*}): T_2\longrightarrow \mathbb R \times \widetilde M$ with $a^{*}(z_2)=0$,
which extends to $u^{*}: \Sigma_2-{P}\longrightarrow \mathbb R \times \widetilde M$.
If $\lim\limits_{\epsilon\to 0}m_{\epsilon}(q_2)\geq\hbar,$ then we can construct bubbles as in subsection \S\ref{bubble_tree_b}.  Otherwise,  $u$ and
$u^{*}$ converges to the same periodic orbit as $z\to q_2$. Note that, near $q_2$, in terms of the local cylinder coordinates the degeneration is given by
\begin{eqnarray}\label{T_rescaling-3}
 (s^*)\upper = s- 2l r\upper,\;\;\;\;\;\;(t^*)\upper= t - \tau\upper.
\end{eqnarray}
In local pseudo-Dauboux coordinate system $(a,\vartheta,{\bf w})$ and the cylinder coordinates of Riemann  surface near $q_{2}$ we write
\begin{eqnarray}
u(s,t)=(a(s,t),\widetilde u(s,t)),\;\;\;\; u^{*}(s^{*},t^{*})=(a^{*}(s^{*},t^{*}),\widetilde u^{*}(s^{*},t^{*})).\end{eqnarray}
Then
\begin{eqnarray}
\lim_{|s |\rightarrow \infty}\widetilde
u(s ,t )=x(kTt),\;\;\;\;\; \lim _{|s ^{*}|\rightarrow \infty}\widetilde
u^*(s^{*},t^{*})=x(kTt^{*}+\vartheta_{0}^*)\end{eqnarray}
 in $C^{\infty}(S^1)$ for some $kT$-periodic orbit
$x(kTt)$.
\v
We call $u^{*}: \Sigma_2\longrightarrow
\mathbb R \times \widetilde M$ a rubber component.

  \v
By the same argument we can construct maps on $\Sigma_{j}-P,j\geq 3.$
\v

\begin{remark}
We use notations above. Suppose that $\sup\limits_{i}\left|a_{1}\upper(z_{1}\upper)\right|\rightarrow\infty$ and there is only one singular point $p\in \Sigma_1 $ with
$$p\upper \rightarrow p,\;\;\sup_{i}\left|u\upper(p\upper)\right|<\infty.$$
We re-scale as above and construct bubble tree as usual. Then all main part of $u(\Sigma)$ are rubber component and there is a bubble tree with $u(p)\in M^+$.
\end{remark}

 \begin{remark}
For Figure 1, it is possible that some $u\upper(T_{j}\upper)\subset M^{+}$ for some $j\geq 2.$ We assume that $$u\upper(T_{1}\upper)\subset M^{+},\;\;\;\;u\upper(T_{6}\upper)\subset M^{+},\;\;\;$$
and
\begin{equation}\label{relation_1_6}
\sup\limits_{i\to\infty} |a\upper(z_{1}\upper)|<\infty,\;\;\;\sup\limits_{i\to\infty}  |a\upper(z_{6}\upper)|<\infty.
\end{equation}
In this case there is a relation between $a\upper(z_j\upper)$, $1\leq j \leq 6$. To see this and to simplify notations, we omit the index $(i)$ and let
 $$l_{1}=a(z_{3})-a(z_{1}),\;\;\; l_{2}=a(z_{4})-a(z_{3}),\;\;$$
 $$l_{3}=a(z_{5})-a(z_{4}),\;\;\;l_{4}=a(z_{6})-a(z_5).\;\;\; $$
 Then
\begin{equation}\label{relation_6_1}
l_{1}+l_{2}+l_{3}+l_{4}=a(z_6)-a(z_1).
\end{equation}
This means that $l_{1}+l_{2}+l_{3}+l_{4}$ is a $\mathcal{T}$ re-scaling invariant.
 By \eqref{relation_1_6},  there is only three independent parameter.

\v
If we start from $\Sigma_6$ to re-scale by the above procedure, we get the same result. Similarly, we start from any component
 $\Sigma_i$, $i\geq 3$, we get also the same result.
This suggests a alternative way to do
$\mathcal{T}$-re-scaling for $\Sigma$: we first re-scale for both $\Sigma_1$
and $\Sigma_6$ by
$$\check a\upper(z)=a\upper(z) - a\upper(z_3)= a\upper(z)- l_1\upper - a\upper(z_1),$$
$$\acute{a}\upper(z)=a\upper(z)- a\upper(z_5)=a\upper(z)+ l_4\upper - a\upper(z_6),$$
then we re-scale for both $\Sigma_3$ and $\Sigma_5$ by
$$a^{\bullet\upper}(z)=\check a\upper(z)-\check a\upper(z_4)=\check a\upper(z)-l_2\upper,$$
$$a^{\bullet\upper(z)}=\acute{a}\upper(z) -\acute{a}\upper(z_4) = \acute{a}\upper(z)+ l_3\upper.$$
In view of \eqref{relation_6_1} we find that the above re-scalings are equivalent.
\end{remark}

In the following we consider the singular points. Let $\tilde{q}$ be a singular point of $u\upper$. We discuss two cases:
\v
{\bf Case a.} $\tilde{q}$ is not a node of $\Sigma$. In this case we construct bubble tree as usual.\v
{\bf Case b.} $\tilde{q}$ is a node of $\Sigma$. Suppose that $\tilde{q}=q_1$. We construct bubble as in subsection \S\ref{bubble_tree_b} to get $S^2$, inserted between $\Sigma_1$ and $\Sigma_3$, with $\widetilde{E}(v)\mid_{S^2}> \frac{1}{3} \min\{\hbar ,T \}$. Denote $q_{1}'=\Sigma_{1}\cap S^{2}$, $q_{1} =\Sigma_{3}\cap S^{2}$.

We fix a point $\fkz$ in the compact set of $S^{2}-\{q_1,q_1'\}$. Let $\fkz\upper\in \Sigma\upper$ such that
$\fkz\upper\rightarrow \fkz.$ Then $\Sigma'=\Sigma\cup S^2$ is a semi-stable curve and there exists  a fixed degeneration $ \Sigma\upper\rightarrow\Sigma'.$

 Consider the coordinate transformations
$$\hat{a}\upper=a -kT(-\log\delta_{i}),$$
and
$$(a^\circ)\upper=\hat{a}\upper-\hat{a}\upper(z_{3}\upper). $$
Note that for any $z\in \Sigma_{3}$
$$ \hat{a}\upper(z)-\hat{a}\upper(z_{3}\upper) =  {a}\upper(z)- {a}\upper(z_{3}\upper).$$
Then we have
$$ (a^\circ)\upper(z)=(\check a)\upper(z) .$$
To simplify notations, we omit the index $(i)$ and let
 $$l_{1}'=a(\fkz)-a(z_{1}),\;\;\;l_{1}^{\star}=a(z_{3})-a(\fkz),\;\;\; l_{2}=a(z_{4})-a(z_{3}),\;\;$$
 $$l_{3}=a(z_{5})-a(z_{4}),\;\;\;l_{4}=a(z_{6})-a(z_5).\;\;\; $$  Then
\begin{equation}\label{relation_6_1}
l_{1}'+l_1^{\star}+l_{2}+l_{3}+l_{4}=a(z_6)-a(z_1).
\end{equation}
We conclude that the $\mathcal{T}$-rescalling based on $\Sigma\upper\rightarrow\Sigma$ and
the $\mathcal{T}$-rescalling based on $\Sigma\upper\rightarrow\Sigma'$ are equivalent.

\subsection{Equivalent $\mathcal{D}\mathcal{T}$- rescaling}\label{DT_rescaling}
 We consider the example in subsection  \ref{T_rescaling} (see Figure 1). For simplicity we only consider the degeneration of $\Sigma_{(r)}$ at $q_{2}$ and consider the {\bf Case 2} here. Let $q^{(i)}\in \Sigma\upper$,
$q^{(i)}\rightarrow q_2$. Suppose that there exists a constant $N>0$ such that
\begin{equation}\label{bound_energy_rubber1}
\lim_{i\rightarrow \infty} \sup \tilde {E}(u\upper; N-1\leq s\leq 2lr \upper-N+1) \leq
\frac{1}{2}\min\{\hbar,T\}.
\end{equation}
\v
We have two $\mathcal{D}\mathcal{T}$-rescaling:
\v
{\bf A.} The degeneration of Riemann surfaces in the Delingne-Mumford space, together with the $\mathcal{T}$-re-scaling in Section \S\ref{T_rescaling}, we have a $\mathcal{D}\mathcal{T}$-rescaling given by \eqref{T_rescaling-1}, \eqref{T_rescaling-2} and \eqref{T_rescaling-3}
\v
{\bf B.} There is  another re-scaling as following.  Assume that $\Sigma\upper$ degenerate at $q_2$ with the formula \eqref{T_rescaling-3}. We choose a new coordinate system $\check{a}\upper$ by
\begin{equation}\label{new_a_check}
\check{a}\upper= a -2kTlr\upper.
\end{equation}
Define
\begin{equation}\label{new_cord_check}
\check u \upper(z)=(\check{a}\upper(z),\check { \tilde u}  \upper(z))=( a\upper(z)-2kTlr\upper, \tilde u\upper(z)).
\end{equation}
Then the two coordinate systems $\check{a}\upper$ and $(a^*)\upper$ satisfies
\begin{equation}\label{equivalent-1}
\check{a}\upper= (a^*)\upper+a\upper(z_{2}\upper) -2kTlr\upper.\end{equation}
Choose two sequences of points $\fkz_1\upper,\fkz_2\upper$ such that $(s,t)(\fkz_{1}\upper)=(N,0)$ and  $((s^{*})\upper,(t^{*})\upper)(\fkz_{2}\upper)=(-N,0).$ Obviously, $\fkz_j\upper \rightarrow \fkz_{j} \in \Sigma_{j}-P,$ $j=1,2$.
It follows from the convergence of $u\upper$ and $(u^{*})\upper$ on compact sets that
\begin{equation}\label{eqn_a_fkz}
|a\upper(\fkz_{1}\upper)-kTs(\fkz_{1}\upper)|\leq C_2,\;\;\;\; |a\upper(\fkz_{2}\upper)-a\upper(z_{2}\upper)|=|(a^{*})\upper(\fkz_{2}\upper)-(a^{*})\upper(z_{2}\upper)|\leq C_{3}.
\end{equation}
By the same argument of ({\bf 5}) in Lemma \ref{compact_lemma} we have
 \begin{equation*}
|a\upper (\fkz_{2}\upper)-kTs(\fkz_{2}\upper)|\leq C_3'.
\end{equation*}
By \eqref{T_rescaling-3} we have
$$s(\fkz_{2}\upper)-2lr\upper= s^*(\fkz_{2}\upper)= -N.$$
Then
\begin{equation}\label{eqn_bound_fkz}
|a\upper(\fkz_{2}\upper) -2kTlr\upper|\leq |a\upper(\fkz_{2}\upper)-kTs(\fkz_{2}\upper)| +|kTs(\fkz_{2}\upper)-2kTlr\upper|\leq C_{3}' +kTN\leq  C_4.
\end{equation}
It follows from \eqref{equivalent-1}, \eqref{eqn_a_fkz} and \eqref{eqn_bound_fkz} that for any point $z\in \Sigma_{2}-P$ we have
\begin{equation} \label{two_re-scalings_A}
\;\;\;\;\;\;\;|\check a\upper(z)-(a^{*})\upper(z)|\leq C
\end{equation}
for some positive constant $C.$

\v
We call two $\mathcal{D}\mathcal{T}$-rescalings are equivalent if there exists a constant $C>0$ independent of $i$ such that \eqref{two_re-scalings_A} holds. We have proved
\begin{lemma}\label{one_nodal_point}
Let $\Sigma=\Sigma_1 \wedge \Sigma_2$,
where $\Sigma_1$ and $\Sigma_2$ are smooth Riemann surfaces of genus $g_1$ and $g_2$ joining at $q$. Assume that
$\Sigma\upper$ is a sequence of smooth Riemann surface which converges to $\Sigma$ in the
Delingne-Mumford space  as $i\to \infty.$ Suppose that
\begin{equation}
\lim_{i\rightarrow \infty} \sup {E}(u\upper; D_{\epsilon}(q)) \leq
\frac{1}{2}\min\{\hbar,T\},
\end{equation}
and restricting on $\Sigma_{2}\setminus D_{\epsilon}(q)$ there is no singular points of $u\upper$. Then the
two $\mathcal{D}\mathcal{T}$-rescalings {\bf A} and {\bf B} are equivalent.
\end{lemma}
\v

The above discussion can be immediately generalized to the case of several nodal points.
\begin{lemma}\label{several_nodal_points}
Let $\Sigma=\Sigma_1 \wedge \Sigma_2$,
where $\Sigma_1$ and $\Sigma_2$ are smooth Riemann surfaces of genus $g_1$ and $g_2$ joining at $q_{1},\cdots,q_{\nu}$. Assume that
$\Sigma\upper=\Sigma_{1}\#_{({\bf r}\upper)}\Sigma_{2}$ be a sequence of smooth Riemann surface which converges to $\Sigma $ in the
Delingne-Mumford moduli space  as $i\to \infty.$ Let $u\upper:\Sigma\upper\rightarrow M^{+}$ be a sequence $J$-holomorphic maps.
Suppose that for any $j\in \{1,\cdots,\nu\},$
\begin{equation}
\lim_{i \rightarrow \infty} \sup {E}(u\upper; D_{\epsilon}(q_{j})) \leq
\frac{1}{2}\min\{\hbar,T\},
\end{equation}
 and
restricting on $\Sigma_{2}-\bigcup_{j=1}^{\nu} D_{\epsilon}(q_{j})$ there is no singular points of $u\upper$. Then the
$\mathcal{D}\mathcal{T}$-rescaling of type {\bf A} and any
$\mathcal{DT}$ re-scalings of type {\bf B} are
equivalent.
\end{lemma}

\subsection{Procedure of re-scaling}\label{procedure_of_re-scaling}

To sum up, for any sequence $\Gamma\upper=(u\upper,\Sigma\upper;{\bf y}\upper,{\bf p}\upper) \in \mathcal{M}_{A}(M^{+},C,g,m+\nu,({\bf k},\mathfrak{e}))$ our procedure
is following:
\v
\begin{enumerate}
\item[(\bf 1)] If there is some $\Sigma_v$ of genus 0, we treat it as in subsection \ref{genus 0}. In the following we
assume that
$(\Sigma\upper;{\bf y}\upper,{\bf p}\upper)$ is stable and converges to $(\Sigma;{\bf y}, {\bf p})$ in
${\overline{\mathcal{M}}}_{g,m+\nu}$.
\item[(\bf 2)] By Lemma \ref{lower_bound_of_singular_points} the number of singular points of $\Sigma$ is finite. Denote by $P\subset \Sigma $ the set of  singular points for $u\upper$,
the nodal points and the puncture points.
\item[(\bf 3)] We first find a component $\Sigma_k$ of $\Sigma$, for example $\Sigma_1$ in Figure 1, such that
$$\left|a_{1}\upper(z_{1}\upper)\right|\leq \min_{j\in I\setminus\{1\}}\left|a_{j}\upper(z_{j}\upper)\right|,\;\;\forall i.$$
 Without loss of generality we suppose that $\sup\limits_{i}\left|a_{1}\upper(z_{1}\upper)\right|<\infty$, that is, $u\upper(z\upper_1)\subset M^+$.
  Find a set $J\subset I$ such that $j\in J$ if and only if 
  $u(\Sigma_j)\subset M^+$, for example $ \Sigma_1$ and $\Sigma_6$ in Figure 1, i.e., $J=\{1,6\}.$
  Let $\Sigma^{1\bigstar}=\Sigma - \bigcup_{j\in J} \Sigma_j$. $\Sigma^{1\bigstar}$ may have several connected components. For every connected component
  of $\Sigma^{1\bigstar}$ we do $\mathcal{T}$-rescaling independently as in \S\ref{T_rescaling}.

We repeat the procedure. Finally we will stop after finite steps.
\item[(\bf 4)] Then we construct bubble tree for every singular point independently to get $\Sigma^{\prime}$, where
${\Sigma}^{\prime}$ is obtained by joining chains of ${\bf P^1}s$ at some double points of $\Sigma$ to separate the
two components, and then attaching some trees of ${\bf P^1}$'s. Then we have
\v
(a) for every nodal point $q \in \mathfrak{N}$ there is a neighborhood $D_{\epsilon}(q)$ so that
$$\lim_{i \rightarrow \infty} \sup {E}(u^{(k)}; D_{\epsilon}(q)) \leq
\frac{\hbar}{2},$$
where $\mathfrak{N}$ denotes the set of all nodal points of $\Sigma^{\prime}$;
\v
(b) restricting on $\Sigma^{\prime}\setminus \bigcup _{q\in \mathfrak{N}}D_{\epsilon}(q)$ there is no singular points of $u^{(i)}$. $\Box$
\end{enumerate}
\v

For every sequence $\Gamma\upper=(u\upper,\Sigma\upper;{\bf y}\upper,{\bf p}\upper) \in \mathcal{M}_{A}(M^{+},C,g,m+\nu,(\mathbf{k},\mathfrak{e}))$, using our procedure we get $\Gamma=(u, \Sigma',{\bf y},{\bf p})$, where
\begin{itemize}
\item[{\bf(A-1)}] ${\Sigma}^{\prime}$ is obtained by joining chains
of ${\bf P^1}s$ at some double points of $\Sigma$ to separate the
two components, and then attaching some trees of ${\bf P^1}$'s. ${\Sigma'} $ is a connected curve with
normal crossings. We call components of ${\Sigma}$ principal components and others bubble components.
\item[{\bf(A-2)}] $u:\Sigma'\rightarrow (M^{+})^{\prime}$ is a continuous map, where $(M^{+})^{\prime}$ is obtained by attaching some ${\mathbb{R}}\times
\widetilde{M}$ to $M^{+}$.  Let $\Sigma_1$ be a connected component of
$\Sigma'$ and $u|_{\Sigma_1}:\Sigma_1\rightarrow {\mathbb{R}}\times \widetilde{M}$. We call $(u; \Sigma_1)$, modulo the translations  ${\mathbb{R}}\times
\widetilde{M}$ and the $S^1$-action on periodic orbits, a rubber component.
\item[{\bf(A-3)}] If we attach a tree of ${\bf P^1}$ at a marked point $y_i$ or a
puncture point $p_i$, then $y_i$ or $p_i$ will be replaced by a point
different from the intersection points on a component of the tree. Otherwise,
the marked points or puncture points do not change;
\item[{\bf(A-4)}] Let $val_v$ be the number of points on ${\Sigma}_{v}$ which are
nodal points or  marked points or puncture points.  In case $u(\Sigma_{v})\subset M^{+}$, if
$u|_{{\Sigma}_{v}}$ is constant  then $val_v+2g_{v}\geq3$; in case $u:\Sigma_{v}\rightarrow \bR\times \widetilde M,$
 if $\tilde{u}|_{\Sigma_{v}}$ is  constant  then $val_v+2g_{v}\geq3$;
\item[{\bf(A-5)}]  $u$ converges exponentially to $(k_{1}T_{\fe_{1}},\cdots,k_{\nu}T_{\fe_{\nu}})$ periodic orbits
$(x_{k_1},...,x_{k_{\nu}})$ as the variable tends to the puncture
$(p_1,...,p_{\nu})$; more precisely, $u$ satisfies \eqref{exponential_decay_a}-\eqref{exponential_decay_y};
\item[{\bf(A-6)}] The restriction of $u$ to each component of $\Sigma' $ is J-holomorphic. Let $q$ be a nodal point of $\Sigma^{\prime}$ . Suppose $q$ is
the intersection point of ${\Sigma}_{v}$ and ${\Sigma}_{w}$. If $q$ is a
removable singular point of $u$, then u is continuous at $q$;
If $q$ is a nonremovable singular point of $u$, then
$u|_{\Sigma_{v}}$  and $u|_{\Sigma_{w}}$ converge exponentially to the same
periodic orbit on $\widetilde{M}$ as the variables tend to the nodal point $q$.
\end{itemize}

\v\v

\subsection{Weighted dual graph with a oriented decomposition}
\v
It is well-known that the moduli space of stable maps in a compact symplectic manifold has a stratification indexed by the combinatorial type of its decorated dual graph. In this section we generalizes this construction to
our setting and in the next section we state Li-Ruan's compatification by using weighted dual graphs.

Let $G$ be a graph. Denote $G=(V(G), E(G))$, where $V (G)$ is a finite nonempty set of vertices and $E(G)$ is a finite set of edges. Suppose that $V=\{v_1,...,v_N\}.$ Given a partition  of $\{1,2,...,N\}$
$$\mathfrak{d}: \{1,\cdots,N\}=(\cup^c_{i=1}I_{i} )\bigcup(\cup_{\alpha=1}^{d} J_{\alpha} ),$$
it induces a decomposition of $V$, still denoted by $\mathfrak{d}$,
$$\mathfrak{d}: V= \mathfrak{A}\bigcup \mathfrak{B},$$
where
$$\mathfrak{A}= \bigcup_{i=1}^{c} W_{I_{i}},\;\;\;\;\;\mathfrak{B}= \bigcup_{\alpha=1}^{d} W_{J_{\alpha}},\;\;\;\;\;W_{I_{i}}=\{v_{k}|k\in I_{i}\},\;\;\;\;\;W_{J_{\alpha}}=\{v_{k}|k\in J_{\alpha}\}. $$
Obviously, $\mathfrak{A}\bigcap \mathfrak{B}=\emptyset.$
Every subset $W_{I_{i}}$ ( resp. $W_{J_{\alpha}}$) determines a induced subgraph $G_{I_{i}}$ ( resp. $G_{J_{\alpha}}$) of $G$.
\v {\bf Assumption 1.} For any $W_{I_{i}},W_{I_{j}}\subset \mathfrak{A}, i\neq j$, there is no edge in $E(G)$ connecting $G_{I_{i}}$ and $G_{I_{j}}$.

 \v
Denote by $R(G)$ all the edges which connect two subgraphs above. Obviously,
$$R(G)=E(G)-(\cup_{i=1}^c E(G_{I_{i}}))\bigcup(\cup_{\alpha=1}^d E(G_{J_{\alpha}})).$$
Let $\ell\in R(G)$ be an edge connecting $v_i$ and $v_j$. We give an orientation to $\ell$, denoted by
\begin{equation}
\overrightarrow{\ell}:v_{i}\xrightarrow{\ell} v_{j}.\end{equation}
Sometimes we denote simply by $v_{i}\xrightarrow{\ell} v_{j}$.

We call $\overrightarrow{\ell}\in R(G)$ an oriented edge, the edges in $E(G)-R(G)$ are called normal edges, or simply edges. If for all $\ell\in R(G)$ we have an orientation, we say that $\mathfrak{d}$ is an oriented decomposition.
\v {\bf Assumption 2.} There is an orientation of $\mathfrak{d}$ such that
  \begin{description}
 \item[(1)] For any edge $\ell\in R(G)$ connecting $v_i\in G_{I_{i}}$ and $v_j\in G_{J_{\alpha}},$ we have
   $$\overrightarrow{\ell}:v_{i}\xrightarrow{\ell} v_{j},$$
   and we denote $G_{I_{i}}\rightarrow G_{J_{\alpha}}.$
 \item[(2)] For any connected subgraph $G_{J_{\alpha}}$ and $G_{J_{\beta}}$, all edges between $G_{J_{\alpha}}$ and $G_{J_{\beta}}$ have the same orientation, that is, if there is an edge $\ell\in R(G)$ connecting $v_{i}\in G_{J_{\alpha}}$ and $v_{j}\in G_{J_{\beta}}$ satisfying $v_{i}\xrightarrow{\ell} v_{j},$ then any edge $\ell'\in R(G)$ connecting $v_{i}'\in G_{J_{\alpha}}$ and $v_{j}'\in G_{J_{\beta}}$ satisfies  $v_{i}'\xrightarrow{\ell} v_{j}'.$ We denote $G_{J_{\alpha}}\rightarrow G_{J_{\beta}}.$
 \item[(3)] For any subgraph sequence $ G_{J_{\alpha_{1}}},\cdots, G_{J_{\alpha_{l}}}$ satisfying
  $$G_{J_{\alpha_{1}}}\rightarrow G_{J_{\alpha_{2}}}\rightarrow \cdots \rightarrow G_{J_{\alpha_{l}}} ,$$
  if there is an edge $\ell\in R(G)$ connecting $v_{i}\in G_{J_{\alpha_{i}}}$ and   $v_{j}\in G_{J_{\alpha_{j}}}$ with $j>i$ then we have
  $$v_{i}\xrightarrow{\ell} v_{j},\;\;\;\;\;\;G_{J_{\alpha_{i}}}\rightarrow G_{J_{\alpha_{j}}}.$$
 \item[(4)] For any vertices $v_{j_{0}},v_{j_{i}}\in G,$ if there exist  a walk of the form
 $$v_{j_{0}},e_{j_{1}},v_{j_{1}},e_{j_{2}},\cdots,e_{j_{i}},v_{j_{i}},$$
 where each edge $e_{j_{l}}$, $1\leq l\leq i,$ is the normal edge,
 then $v_{j_{0}},v_{j_{i}}$ belong to the same subgraph; otherwise $v_{j_{0}},v_{j_{i}}$ belong to different subgraphs.
 \end{description}
Let $v_{k}\xrightarrow{\ell} v_{j},$ $v_k\in G_{I_i}$ (or $G_{J_{\alpha}}$),  $v_j\in G_{J_{\beta}}$. When we consider the subgraphs $ G_{I_i}$, $G_{J_{\beta}}$ we attach a half edge  $\ell^{+}$ to $v_{k}$ and attach $\ell^{-}$ to $v_{j}$.
\v
Let $(V(G),E(G))$ be a graph and $\mathfrak{d}$ be an oriented decomposition satisfying {\bf Assumption 1} and {\bf Assumption 2}. We call the graph $G$ a graph with a oriented decomposition $\mathfrak{d}$, denoted by
$(V(G),E(G),\mathfrak{d})$.

\begin{definition} Let $g$, $m$ and $\nu$ be nonnegative integers.
A $(g,m+\nu,\mathfrak{d})$-weighted dual graph $G$
consists of $(V(G),E(G),\mathfrak{d})$
together with three weights, where
\begin{itemize}
\item[(1)] $(V(G),E(G)$ is a graph, and $\mathfrak{d}$ is an oriented decomposition of $V(G)$,
\item[(2)] $\mathfrak{g}:V(G) \rightarrow \mathbb Z_{\geq 0}$ assigning a nonnegative integer $g_{v}$ to each vertex $v$ such that $$g=\sum_{v\in V(G)}g_{v}+b_{1}(G)$$
where $b_1(G)$ is the first Betti number of the graph $G$;
\item[(3)] assign $m$ ordered tails  $\mathfrak{m}=(t_{1},\cdots,t_{m})$ to $V(G):$  attach $m_v$  tails to $v$ for each $v\in V(G)$,

\item[(4)] assign  $\nu$ ordered half edges  $\mathfrak{l}=(e_{1},\cdots,e_{\nu})$ to $V(G):$ attach $l_v$ half edges to  $v$ for each $v\in V(G)$.
\end{itemize}

\v

\end{definition}

We denote the $(g,m+\nu,\mathfrak{d})$-weighted dual graph $G$ by $(V(G),E(G), \mathfrak{g}, \mathfrak{m}, \mathfrak{l},\mathfrak{d})$.
\v
We introduce a terminology: a vertex $v$ is called an interior vertex if there is no half edge attached it, it is called boundary vertex if there are some half edges attached it.

\begin{definition} Let $g$, $m$ and $\nu$ be nonnegative integers.
A $H$-$(g,m+\nu,\mathfrak{d})$ weighted dual graph $G$
 consists of $(V(G),E(G), \mathfrak{g}, \mathfrak{m}, \mathfrak{l},\mathfrak{d})$ together with three weights:
\begin{itemize}
\item[(1)]  $\mathfrak{h}:V(G) \rightarrow   H_{2}(M^{+},\Re,\mathbb Z)\bigcup H_{2}(M^{+},\mathbb Z)$ assigning a $A_{v}\in  H_{2}(M^{+},\Re,\mathbb Z)$ to each boundary vertex, assigning a $A_{v}\in  H_{2}(M^{+},\mathbb Z)$ to each interior vertex $v\in \mathfrak{A} ,$   assigning a $0$ to each interior vertex $v\in \mathfrak{B} .$  For any  $G_{I_{i}}$ and $G_{J_{\alpha}},$ denote  $$A_{I_{i}}=\sum_{v\in V(G_{I_{i}})}A_{v}\in H_{2}(M^{+},\Re,\mathbb Z),\;\;\;\;\;\;\;\;\; A_{J_{\alpha}}=\sum_{v\in V(G_{J_{\alpha}})}A_{v}\in H_{2}(M^{+},\Re,\mathbb Z),$$
and
$$A=\sum_{i=1}^{c}A_{I_{i}}+\sum_{\alpha=1}^{d}A_{J_{\alpha}}  .$$
\item[(2)] assign  $\nu$ ordered weights ${\bf k}=(k_1,...,k_\nu)$ and weights $\mathfrak{e}=(\mathfrak{e}_{1},\cdots,\mathfrak{e}_{\nu})$ to the half edges $\mathfrak{l}=(e_{1},\cdots,e_{\nu})$ such that $\mathfrak{l}$ becomes weighted half edges $\mathfrak{l}^{(\bf k,\mathfrak{e})}=((k_{1},\mathfrak{e}_{1})e_{1},\cdots,(k_{\nu},\mathfrak{e}_{\nu})e_{\nu})$.
\item[(3)]  $\mathfrak{k}: R(G)\rightarrow   \mathbb Z^{+},$
for each $\ell \in R(G)$ with $v_{k}\xrightarrow{\ell} v_{j},$ $v_k\in G_{I_i}$ (or $G_{J_{\alpha}}$),  $v_j\in G_{J_{\beta}},$  assigning  $\mathfrak{k}(\ell)=  k_{\ell}>0 $ such that
  $$\sum_{e_{j}\in G_{I_{i}}}k_{j}+\sum_{\ell^{+}\in G_{I_{i}}}k_{\ell}>0, \mbox{ for any } G_{I_{i}}  $$
   and
   $$ d\lambda(A_{J_{\alpha}})=  \sum_{e_{j}\in G_{J_{\alpha}}}k_{j}\cdot T_{\fe_{j}}+\sum_{\ell^{+}\in G_{J_{\alpha}}}k_{\ell^{+}}\cdot T_{\ell^{+}}-\sum_{\ell^{-}\in G_{J_{\alpha}}}k_{\ell^{-}}\cdot T_{\ell^{-}}, \mbox{ for any } G_{J_{\alpha}}. $$
   and for each boundary vertex $v\in G_{J_{\alpha}}$        \begin{equation}
  d\lambda(A_{v})= \sum_{ e_{i}^{+}\in \fl_{v}^{+}} k^{+}_{i}\cdot T_{e_{i}^{+}}-\sum_{e_{i}^{-}\in \fl_{v}^{-}} k^{-}_{i}\cdot T_{e_{i}^{-}},   \end{equation}
  where $\fl_{v}^{\pm}$ is the subset of $\fl^{\pm},$ the half edges attached to  $v$
 \end{itemize}

\v

\end{definition}

We denote the $H$-$(g,m+\nu,\mathfrak{d})$ weighted dual graph $G$ by $(V(G),E(G),\mathfrak{g}, \mathfrak{m}, \mathfrak{l}^{(\bf k,\mathfrak{e})},\mathfrak{k},\mathfrak{h},\mathfrak{d}  )$.

\v

By a leg of $G$ we mean either a tail or a
half-edge.

\begin{definition}
Let $G$ be a $(V(G),E(G),\mathfrak{g}, \mathfrak{m}, \mathfrak{l}^{(\bf k,\mathfrak{e})},\mathfrak{k},\mathfrak{h},\mathfrak{d}  )$ graph. A vertex $v$ is called stable if one of the following holds:
\begin{itemize}
\item[(1)]  $2g_{v}+ val(v)\geq 3,$ where $val(v)$ denotes the sum of the number of legs attached to $v$;
\item[(2)] $A_v\not=0$,  when $v\in \mathfrak{A} ;$
\item[(3)]$d\lambda(A_v)\not=0$ when $v\in \mathfrak{B}.$
\end{itemize}
  $G$ is called stable if all vertices are stable.
 \end{definition}

Let $G$ be a stable $(V(G),E(G),\mathfrak{g}, \mathfrak{m}, \mathfrak{l}^{(\bf k,\mathfrak{e})},\mathfrak{k},\mathfrak{h},\mathfrak{d}  )$ graph. Then all subgraphs $G_{I_{i}},$ $G_{J_{\alpha}}$ are stable.
\v
   Two $H$-$(g,m+\nu,\mathfrak{d})$ weighted dual  graphs  $G_1$ and $G_2$ are called isomorphic
if there exists a bijection $T$ between  their vertices and edges keeping  oriented decomposition and all weights.
\v
 Let $S_{\mathfrak{g},\mathfrak{m},\mathfrak{l}^{(\bf k,\mathfrak{e})},\mathfrak{k},\mathfrak{h},\mathfrak{d}}$
 be the set of isomorphic classes of  $H$-$(g,m+\nu,\mathfrak{d})$ weighted dual graphs. Given $g$, $m$, $\nu$, $A\in H_{2}(M^{+},\Re,\mathbb Z)$ and two weights  $\mathbf{k}=(k_1,...,k_{\nu})$, $\mathfrak{e}=(\mathfrak{e}_{1},\cdots,\mathfrak{e}_{\nu})$, denote by $S_{g,m+\nu,A,\mathfrak{l}^{(\bf k,\mathfrak{e})}}$ the union of all possible $S_{\mathfrak{g},\mathfrak{m},\mathfrak{l}^{(\bf k,\mathfrak{e})},\mathfrak{k},\mathfrak{h},\mathfrak{d}}$.
\v

\v
\subsection{Li-Ruan's Compactification
}\label{sect_1.2}

In this section we state Li-Ruan's compactification on the moduli space of  maps in terms of language of graphs.
\v
Let $G$ be a stable $(V(G),E(G),\mathfrak{g}, \mathfrak{m}, \mathfrak{l}^{(\bf k,\mathfrak{e})},\mathfrak{k},\mathfrak{h},\mathfrak{d}  )$  graph with $N$ vertices $(v_1,...,v_N)$, $m$ tails and $\nu $ half edges, and $(\Sigma,{\bf y}, {\bf p})$ be a semi-stable curve with $m$ marked points and $\nu$ puncture points. Let $A\in H_{2}(M^{+},\Re,\mathbb Z)$.  A  stable $J$-holomorphic map of type $G$ is a quadruple $$(u;\Sigma, {\bf y},{\bf p})$$ where  $u:\Sigma\rightarrow (M^{+})^{\prime}$ is a continuous map, $(M^{+})^{\prime}$ is obtained by attaching some ${\mathbb{R}}\times
\widetilde{M}$ to $M^{+}$,
satisfying the following conditions:
\begin{itemize}
\item [{\bf [A-1]}] $\Sigma = \bigcup_{v=1}^N \Sigma_v$, where  each $v\in V(G)$  represents a
 smooth component $\Sigma_{v}$ of $\Sigma$.
\item [{\bf [A-2]}]for the $i$-th tail attached to the vertex $v$ there exists the $i$-th marked  point $y_i\in \Sigma_{v}$, $m_v$
is equal to the number of the marked points on $\Sigma_v$,
\item [{\bf [A-3]}]for the $j$-th half edge attached to the
vertex $v$ there  exists $j$-th puncture point $p_{j}\in \Sigma_{v}$, $l_v$
is equal to the number of puncture points on $\Sigma_{v}$;
\item [{\bf [A-4]}] if there is an edge connected the vertices $v$ and $w$, then there exists a node between $\Sigma_{v}$ and $\Sigma_{w},$  the number  of edges between $v$ and $w$ is equal to the number of node points between $\Sigma_{v}$ and $\Sigma_{w}$;
\item[{\bf[A-6]}] the restriction of $u$ to each component $\Sigma_v$ is $J$-holomorphic.
\item[{\bf[A-7]}]  $u$ converges exponentially to $(k_{1}\cdot T_{\fe_{1}},\cdots, k_{\nu}\cdot T_{\fe_{\nu}})$ periodic orbits
$(x_{k_1},...,x_{k_{\nu}})$ satisfying $x_{k_{i}}\subset \mathcal{F}_{\mathfrak{e}_{i}}$ as the variable tends to the puncture
$(p_1,...,p_{\nu})$; more precisely, $u$ satisfies \eqref{exponential_decay_a}-\eqref{exponential_decay_y};
\item[{\bf[A-8]}] let $q$ be a nodal point of $\Sigma$. Suppose $q$ is
the intersection point of ${\Sigma}_{v}$ and ${\Sigma}_{w}$ associated to the edge $\ell\in R(G).$ Then  $k_{\ell}>0$,  $u|_{\Sigma_{v}}$  and $u|_{\Sigma_{w}}$ converge exponentially to the same $k_{\ell}$
periodic orbit $x_{k_{\ell}}$ on $\widetilde{M}$ as the variables tend to the nodal point $q$.
\item[{\bf[A-9]}] For any $v\in V(G),$ $[ u (\Sigma_{v})]=A_{v}\in H_{2}( M ^{+},\Re,\mathbb Z) $ when $v$ is a boundary vertex; $[ u (\Sigma_{v})]=A_{v}\in H_{2}( M ^{+},\mathbb Z) $ when $v\in \mathfrak{A}$ is an interior vertex, $[ u (\Sigma_{v})]=0$ when $v\in \mathfrak{B}$ is an interior vertex; $A=\sum_{i=1}^vA_v.$

\end{itemize}
\v

\v
\begin{definition}\label{equivalent relation}
Two  stable $J$-holomorphic maps $\Gamma =(u, (\Sigma,{\bf j}),{\bf y},{\bf p})$ and $\check{\Gamma} =(\check{u}, (\check{\Sigma},
\check{\bf j}),{\bf \check{y}},{\bf \check{p}})$ of type $G$
 are called equivalent if there exists a diffeomorphism $\varphi:\Sigma\rightarrow \check{\Sigma}$  such that  it can be lifted to bi-holomorphic isomorphisms $\varphi_{v w}:(\Sigma_{v},j_v)\rightarrow (\check\Sigma_{w},\check j_w)$ for
 each component $\Sigma_{v}$ of $\Sigma$, and
\begin{itemize}
\item[{\bf(1)}] $\varphi(y_i)= \check{y_i}$,  $\varphi(p_j)= \check{p_j}$ for any $1\leq i\leq l$, $1\leq j\leq \nu$,
\item[{\bf(2)}] $u(\Sigma_v)$ and $\check{u}\circ \varphi(\Sigma_v)$ lie in the same holomorphic block for any $\Sigma_v$,
\item[{\bf(3)}]    In case $v\in \mathfrak{A}$ we have $\check{u}\circ \varphi= u$ on $\Sigma_{v}$; In case $v\in \mathfrak{B}$, we have
  $$\check{\tilde{u}}\circ \varphi= \tilde{u},\;\;\check{a}\circ \varphi= a+C\;\;on \;\Sigma_{v}.$$
  Moreover, near every periodic orbit $x$, $\tilde{u}$ and $\check{\tilde{u}}\circ \varphi$ may differ by a canonical coordinate transformation \eqref{canonical}.
  \end{itemize}
Denote by $\mathcal{M}_{G}$ the space of  the equivalence class of  stable $J$-holomorphic maps of type $G$.
\end{definition}
 \begin{remark} Every subgraph $G_{I_{i}}$ determines a holomorphic block of type $G_{I_{i}}$ in $M^+$, every subgraph $G_{J_{\alpha}}$ determines a holomorphic rubber block of type $G_{J_{\alpha}}$ in ${\mathbb{R}}\times \widetilde{M}$.
Roughly speaking, $\mathcal{M}_{G}$ is the gluing of several holomorphic blocks.
\end{remark}

\v
We define the automorphism group of $u$:
\begin{align*}Aut(u)=\{& \varphi \mid \varphi:\Sigma \rightarrow \Sigma \mbox{ is a holomorphic isomorphic such that } {\bf (1)}, {\bf (2)}  \\
&\mbox{ and {\bf (3)} hold
in {\bf Definition} \ref{equivalent relation}}.\}
\end{align*}
It is easy to see
\begin{lemma}  For any  stable holomorphic map $u$ of type $G$ the automorphism group $Aut(u)$ is finite .
\end{lemma}

\v
Given $g$, $m$, $\nu$, $A\in H_{2}(M^{+},\Re,\mathbb Z)$, two weights  $\mathbf{k}=(k_1,...,k_{\nu})$ and $\mathfrak{e}=(\mathfrak{e}_{1},\cdots,\mathfrak{e}_{\nu})$,
we define
 $${\overline{\mathcal{M}}}_{A}(M^{+},C;g,m+\nu,({\bf k},\mathfrak{e}))=\bigcup_{G\in S_{g,m+\nu,A,\mathfrak{l}^{(\bf k,\mathfrak{e})}}} \mathcal M_{G} .$$

\v\v

This gives a stratification of ${\overline{\mathcal{M}}}_{A}(M^{+},C;g,m+\nu,({\bf k},\mathfrak{e}))$. Denote by ${\mathcal D}^{J,A}_{g,m+\nu,{\bf k},\mathfrak{e}}$ the number of all possible $S_{\mathfrak{g},\mathfrak{m},\mathfrak{l}^{({\bf k},\mathfrak{e})},\mathfrak{k},\mathfrak{h},\mathfrak{d}}$. The following Lemma is obvious.
\v
\begin{lemma} ${\mathcal
D}^{J,A}_{g,m+\nu,{\bf k},\mathfrak{e}}$ is finite.
\end{lemma}
\v
We immediately obtains \vskip 0.1in
\noindent
\begin{theorem}\label{compact_moduli_space}
 ${\overline{\mathcal{M}}}_{A}(M^{+},C;g,m+\nu,({\bf k},\mathfrak{e}))$ is compact.
\end{theorem}
\v
\v

\begin{remark} Let $u:  S^2\longrightarrow {\mathbb R} \times \tilde{M}$ be a J-holomorphic map with $\tilde {E}(u)=0.$ Then $u$ must be $u=(kTs+d,x(kt+\vartheta_{0}))$. If the number
of the special points (including nodal points, puncture points and marked points) $\leq 2$, $u$ is called a unstable second class ghost bubble.  In our compactification there is no unstable second class ghost bubble.
\end{remark}

\section{Gluing theory-Pregluing}\label{gluing_pregluing}

The gluing is the inverse of the degeneration. In section \S\ref{compact_theorem} we see that for a sequence $\Gamma\upper=(u\upper,\Sigma\upper;{\bf y}\upper,{\bf p}\upper) \in {\mathcal M}_{A}(M^{+},C;g,m+\nu,({\bf k},\mathfrak{e}))$, each component degenerates independently, and the bubble trees are also constructed independently,
furthermore each node of a connected component degenerates also independently. So we also glue independently for each nodal points. Recall that
there are some freedoms of choosing the coordinates $a$, $\vartheta$ for $M^{+}$ and $\bR\times \widetilde M$ (see subsection \S\ref{subsection:7.2}).
\v
For any $r$ we glue $M^+$ and ${\mathbb{R}}\times \widetilde{M}$ with
parameter $r$ to get again $M^+$. We cut off the
part of $M^+$ and $\bR\times \widetilde M$ with cylindrical
coordinate $|a^{\pm}|>\frac{3lr}{2}$ and glue the remainders along
the collars of length $lr$ of the cylinders with the gluing formulas:
\begin{eqnarray} \label{gluing_local_relative_node_a}
a_1=a_2 + 2l r.
\end{eqnarray}

\subsection{Gluing one relative node for $M^+\cup( {\mathbb{R}}\times \widetilde{M})$}\label{relative_node_M_R_M}

Let $b=(u_1,u_2; \Sigma_1 \wedge \Sigma_2,j_1,j_2)\in {\mathcal{M}}_{A_1}(M^{+},C;g_1,m_1+1,(k,\fe_{1}))\times{\mathcal{M}}_{A_2}({\mathbb{R}}\times \widetilde{M},g_2,m_2+\nu + 1,(k,\fe_{1}), ({\bf k}^+,\fe^{+}))$,
where $(\Sigma_1,j_1)$ and $(\Sigma_2,j_2)$ are smooth Riemann surfaces of genus $g_1$ and $g_2$ joining at $q$ and
$u_1: \Sigma_1 \rightarrow M^+$, $u_2: \Sigma_2 \rightarrow {\mathbb{R}}\times \widetilde{M}$ are $J$-holomorphic maps such that
$u_i(z)$ converge to the same $kT$-periodic orbit $x$ as $z\rightarrow q$. To describe the maps $u_i$ we choose a local pseudo-Darboux coordinate system $(a_i,\vartheta_i,{\bf w})$, $i=1,2$, near $x$.
Suppose that
$$a_i(s_i,t_i)-kTs_i-\ell_i\rightarrow 0\;\;\;\;
\vartheta_i(s_i,t_i)-kt_i-\vartheta_{i0}\rightarrow 0.$$
If we choose a different origin in the periodic orbit $x$, we have a different coordinate system $\vartheta^*_i$. Suppose that
$$\vartheta^*_i=\vartheta_i +\tau_i.$$
Obviously,
$$\vartheta^*_i(s_i,t_i)-kt_i-(\vartheta_{i0}+\tau_i)\rightarrow 0.$$
Without loss of generality we assume that $\vartheta_{10}=0$. Set $\tau=-\vartheta_{20}$. We consider $\tau$ to be a parameter satisfying $0\leq \tau<1$. Then
\begin{equation}\label{gluing_a_r-1}
\vartheta_1=\vartheta_2 + \tau
\end{equation}
Given gluing parameters $(r)=(r,\tau)$ we construct a surface
$\Sigma_{(r)} =\Sigma_1\#_{(r)} \Sigma_2 $ with gluing formulas:
\begin{eqnarray}
&& s_1=s_2 + \tfrac{2lr}{Tk} \\
&&t_1=t_2 + \tfrac{\tau + n}{k} \label{glu_t}
\end{eqnarray}
for some $n \in Z_k$.
\v
To get a pregluing map $u_{(r)}$ from $\Sigma_{(r)}$ we set\\
\[
u_{(r)}=\left\{
\begin{array}{ll}
u_1 \;\;\;\;\; on \;\;\Sigma_{10}\bigcup\left\{(s_1,t_1)|0\leq s_1 \leq
\frac{lr}{2Tk}, t_1 \in S^1 \right\}    \\  \\
\left(kTs_1, x(kTt_1)\right) \;\;\;\;\;
on \;\;\left\{(s_1,t_1)| \frac{3lr}{4Tk}\leq s_1 \leq
\frac{5lr}{4Tk}, t_1 \in S^1 \right\}  \\   \\
u_2 \;\;\;\;\; on \;\;\Sigma_{20}\bigcup\left\{(s_2,t_2)|0\geq s_2
\geq - \frac{lr}{2Tk}, t_2 \in S^1 \right\},
\end{array}
\right.
\]
where $x(kTt_1)=(kt_1,0)$ as $\vartheta_{10}=0$.
\v
To define the map $u_{(r)}$ in the remaining part we fix a smooth cutoff
function $\beta : {\mathbb{R}}\rightarrow [0,1]$ such that
\[
\beta (s)=\left\{
\begin{array}{ll}
1 & if\;\; s \geq 1 \\
0 & if\;\; s \leq 0
\end{array}
\right.
\]
and $\sqrt{1-\beta}$ is a smooth function, $|\beta^{\prime}(s)|\leq 2.$
We assume that $r$ is large enough such that $u_i$ maps the tube
$\{(s_i,t_i)||s_i|\geq \frac{lr}{4Tk},t_i \in S^1 \}$ into a domain
with pseudo-Darboux coordinates $(a_i,\vartheta_i, {\bf w})$. We write
 $u_{(r)} =\left(a_{(r)} , \widetilde{u}_{(r)} \right) $ and define\\
\begin{equation}\label{gluing_a_r}
a_{(r)}= kTs_1  +  \left(\beta(3-\tfrac{4Tks_1}{lr}) (a_1(s_1,t_1)-
kTs_1) + \beta(\tfrac{4Tks_1}{lr}-5) (a_2(s_2,t_2)- kTs_2
) \right),
\end{equation}
\begin{equation}
\label{gluing_t_u_r}
\widetilde{u}_{(r)}=x(kTt_1)+
 \left(\beta(3-\tfrac{4Tks_1}{lr}) (\widetilde{u}_1(s_1,t_1)-x(kTt_1))+
\beta(\tfrac{4Tks_1}{lr}-5)
(\widetilde{u}_2(s_2,t_2)-x(kTt_2)) \right),
\end{equation}
where $x(kTt_2)=(kt_2,0)$ as \eqref{glu_t} and \eqref{gluing_a_r-1}.
It is easy to check that   $u_{(r)}$ is a
smooth function.
\v
\subsection{Gluing two relative nodal points for $M^+\cup( {\mathbb{R}}\times \widetilde{M})$}\label{relative_node_M_R_M-1}

Recall Lemma \ref{several_nodal_points}, we glue   two relative nodal points independently.

Let $b=(u_1,u_2; \Sigma_1 \wedge \Sigma_2,j_1,j_2)\in {\mathcal{M}}_{A_1}( M^{+},C;g_1,m_1+2,({\bk},\fe))\times {\mathcal{M}}_{A_2}(\bR\times\widetilde{M},C;g_2,m_2+2,({\bk}^{-},\fe^{-}),(k',\fe'))$,
where ${\bk}={\bk}^{-}=(k_{1},k_2), {\fe}={\fe}^{-}=(\fe_{1},\fe_{2}),$ and  $k'$ satisfying $\sum_{j=1}^{2}k_{j}<k'.$
 Here $(\Sigma_1,j_1)$ and $(\Sigma_2,j_2)$ are smooth Riemann surfaces of genus $g_1$ and $g_2$ joining at $q_{1},q_{2},$ and
$u_1: \Sigma_1 \rightarrow M^+$, $u_2: \Sigma_2 \rightarrow \bR\times\widetilde{M}$ are J-holomorphic maps such that
$u_{j}(z)$ converge to the same $k_{j}T_{\fe_{j}}$-periodic orbit $x_{j}(k_{j}T_{\fe_{j}})$ as $z\rightarrow q_{j},j=1,2$ (see Figure 2).  Choose cylinder coordinates  $(s_{1j},t_{1j})$ and  $(s_{2j},t_{2j})$ on $\Sigma_{1}$ and $\Sigma_{2}$ near node $q_{j}.$  We choose local pseudo- Darboux coordinate systems $(a_{1},\vartheta_{1j},{\bf w}_{j})$ on the cylinder end of $M^+$, $(a_{2},\vartheta_{2j},{\bf w}_{j})$ on ${\mathbb{R}}\times \widetilde{M}$
near $x_{j}$, where ${\bf w}_{j}$ is a local coordinates near $x_{j}$.
Suppose that
$$a_{1}(s_{1j},t_{1j})-kTs_{1j}-\ell_{1j}\rightarrow 0,\;\;\;\;
\vartheta_{1j}(s_{1j},t_{1j})-kt_{1j}-\vartheta_{1j0}\rightarrow 0,\;\;\;\;\;j=1,2,$$
$$a_{2}(s_{2j},t_{2j})-kTs_{2j}-\ell_{2j}\rightarrow 0,\;\;\;\;
\vartheta_{2j}(s_{2j},t_{2j})-kt_{2j}-\vartheta_{2j0}\rightarrow 0,\;\;\;\;\;j=1,2.$$
\begin{figure}[H] \label{Figure3}
\centering\includegraphics[height=5cm]{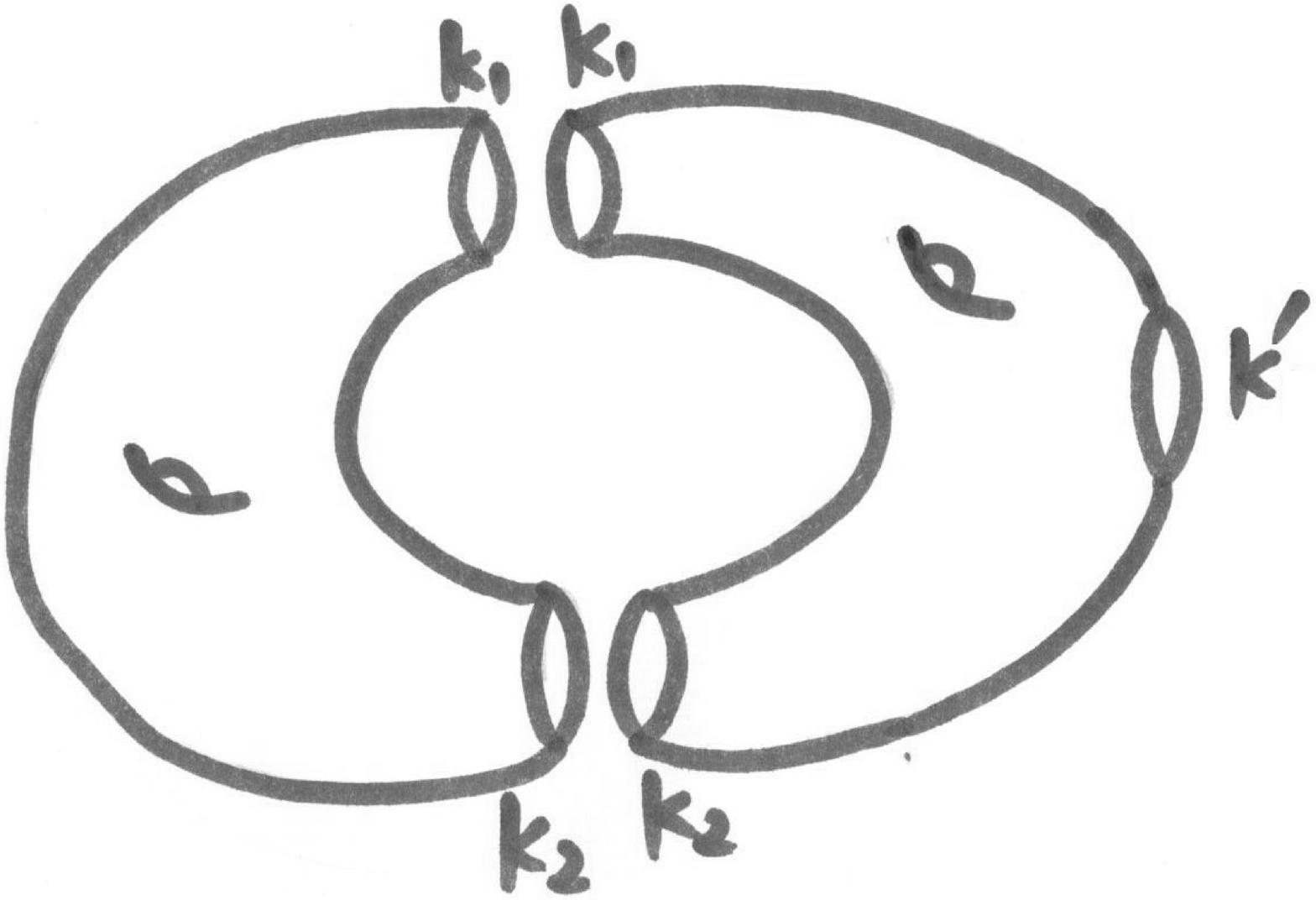}
\caption{}
\end{figure}
Without loss of generality we assume that $\ell_{1j}=\ell_{2j}=0$, for $j=1,2$. Denote $T_{j}=T_{\fe_{j}}$, $j=1,2.$

\v
For any parameter $\varrho>0$, we can gluing $M^+$ and $\bR\times\widetilde{M}$ to get $M^{+}_{\varrho}$  with gluing formula:
\begin{eqnarray}\label{gluing_relative_node_a}
&&a_{1}=a_{2} + 2l\varrho.
\end{eqnarray}
\v
For each $x_j$ we take $\mathcal{T}$ translation:
\begin{eqnarray}\label{freedoms_1}
a_{1j}= a_{1}-c_{1j},\;\;a_{2j}=a_{2}-c_{2j}
\end{eqnarray}
for some constants $c_{1j}\geq 0$.
Then $(a_{1j}, \vartheta_{1j}, {\bf w}_{j})$ and $(a_{2j}, \vartheta_{2j}, {\bf w}_{j})$ are local coordinate systems near $x_{j}$ over cylindrical end of $M^+$ and over $\bR\times\widetilde{M}$.

 Choose $R_{0}$ such that
\begin{equation}\label{condition_R_0}
 \sum_{i=1}^{2}\sum_{j=1}^{2}\widetilde{E}(u_i; |s_{ij}|\geq \frac{R_{0}}{2})\leq \frac{\min\{\hbar,T\}}{8},\;\;\;\;\;
 2\mathcal{C}_{1}\fc_{1}^{-1}e^{-\frac{\fc R_{0}}{4}}\leq \frac{\hbar_{1}}{8},\;\;\;\;
 \end{equation}
 where $\hbar,\mathcal{C}_{1},\fc$ and $\hbar_{1}$ are the constants in Theorem \ref{tube_exponential_decay} and Theorem \ref{details_implicit_function_theorem}.

For any $r_{j}\geq R_0$ we glue $M^+$ and $\bR\times\widetilde{M}$ in coordinates $a_{1j}$  and $a_{2j}$ with the gluing formula:
\begin{eqnarray} \label{gluing_local_relative_node_a}
&&a_{1j}=a_{2j} + 2lr_{j}.
\end{eqnarray}

\v
Now we return to the coordinates $a_{i},i=1,2$.
By relation \eqref{freedoms_1} the gluing formula can be re-written as
  \begin{eqnarray}
a_{1}=a_{2}+c_{1j}-c_{2j} + 2l r_{j}.
\end{eqnarray}
Choose
  \begin{equation}\label{choose parameters}
c_{1j}=-c_{2j}= l(\rho- r_{j}). \end{equation}
By \eqref{choose parameters}, in order to get $M_{\varrho}^{+}$ we need
\begin{equation}\label{parameters bound}
 R_0\leq r_{j} \leq \varrho,\;\;\;\;\;j=1,2.
\end{equation}

Now we glue $J$-holomorphic maps. We express $(u_1,u_2)$ in terms of the coordinates $(a_{1j},\vartheta_{1j})$ and $(a_{2j},\vartheta_{2j})$:
$$a_{1j}(s_{1j},t_{1j})-k_{j}T_{j}s_{1j}+c_{1j}\rightarrow 0,\;\;\;\;
\vartheta_{1j}(s_{1j},t_{1j})-k_{j}t_{1j}-\vartheta_{1j0}\rightarrow 0,\;\;\;\;\;j=1,2,$$
$$a_{2j}(s_{2j},t_{2j})-k_{j}T_{j}s_{2j}+c_{2j}\rightarrow 0,\;\;\;\;
\vartheta_{2j}(s_{2j},t_{2j})-k_{j}t_{2j}-\vartheta_{2j0}\rightarrow 0,\;\;\;\;\;j=1,2.$$
As in subsection \S \ref{relative_node_M_R_M} we assume $\vartheta_{1j0}=0$ and consider $\tau_j=-\vartheta_{2j0}$ as parameters, and construct a surface
$\Sigma_{(r_{j})} =\Sigma_1\#_{(r_{j})} \Sigma_2 $ with gluing formulas:
\begin{eqnarray}
\label{gluing_s}&& s_{1j}=s_{2j} + \tfrac{2lr_{j}}{k_{j}T_{j}},\;\;\;\;\;\;j=1,2, \\
\label{gluing_t}&&t_{1j}=t_{2j} + \tfrac{\tau_{j} + n_{j}}{k_{j}},\;\;\;\;\;\;j=1,2,
\end{eqnarray}
for some $n_{j} \in Z_{k_{j}},j=1,2$. By \eqref{choose parameters}, in terms of $\varrho$ the gluing formulas \eqref{gluing_s} and \eqref{gluing_t} can be written as
\begin{eqnarray}
&& s_{1j}=s_{2j} + \tfrac{2l\rho-c_{1j}+c_{2j}}{k_{j}T_{j}} ,\;\;\;\;\;\;j=1,2,\\
&&t_{1j}=t_{2j}  + \tfrac{\tau_j + n_{j}}{k_{j}},\;\;\;\;\;\;j=1,2,
\end{eqnarray}
for some $n_{j} \in Z_{k_{j}},j=1,2$.
 We construct pre-gluing map $u_{(r_{j})}$ as in subsection \S\ref{relative_node_M_R_M}.  Then    we have
$$u_{(r_{1})}=u_{(r_{2})}=u_1,\;\mbox{ in }\;\bigcap_{j=1}^2\{|s_{1j}| \leq
\frac{lr_{j}}{4T_{j}k_{j}}\},\;\;\;\;\;\;\;\;u_{(r_{1})}=u_{(r_{2})}= u_{2},\;\mbox{ in }\;\bigcap_{j=1}^2\{|s_{2j}| \leq
\frac{lr_{j}}{4T_{j}k_{j}}\}.$$
The pre-gluing above can be generalized immediately to the case of gluing several nodal points.

\v
\v
\subsection{Norms on $C^{\infty}(\Sigma_{(r)};u_{(r)}^{\ast}TM_{(r)}^{+})$}

 We only consider the case of gluing one node, the other cases are the same.
For any $\eta \in
C^{\infty}(\Sigma_{(r)};u_{(r)}^{\ast}TM_{(r)}^{+} \otimes \wedge^{0,1}))$,
let $\eta_{i}$ be its restriction to the part
$\Sigma_{i0}\bigcup\{(s_i,t_i)|\;|s_i| < \frac{3lr}{2kT}\}$, extended by zero
to yield a section over $\Sigma_i$. Define
\begin{equation}\|\eta\|_{p,\alpha,r}=\|\eta_{1} \|_{\Sigma_1,p,\alpha} +
\|\eta_{2} \|_{\Sigma_2,p,\alpha}.
\end{equation}
Denote the resulting completed spaces by $L_{r}^{p,\alpha}$.

We define a norm $\|\cdot\|_{1,p,\alpha,r}$ on
$C^{\infty}(\Sigma_{(r)};u_{(r)}^{\ast}TM_{(r)}^{+}).$ For any section
$h\in C^{\infty}(\Sigma_{(r)};u_{(r)}^{\ast}TM_{(r)}^{+})$, denote
  \begin{eqnarray}
 && h_0
=\int_{\mathbb{S}^1}h\left(\frac{lr}{Tk},t\right)dt, \\
&&h_1 = (h-\hat h_{0})\beta\left(\frac{3}{2}-\frac{Tks_1}{lr}\right)\\
&&h_2 = (h-\hat h_{0})\left[1-\beta\left(\frac{3}{2}-\frac{Tks_1}{lr}\right)\right].
 \end{eqnarray}
 We define
 \begin{equation}\| h\|_{1,p,\alpha,r}=\|h_1\|_{\Sigma_1,1,p,\alpha} +
\|h_2\|_{\Sigma_2,1,p,\alpha}+|h_{0}|.
\end{equation}
Denote the resulting completed spaces by ${\mathcal W}_{r}^{1,p,\alpha}$.

\v
We introduce some notations. Choose the cylinder coordinates near puncture points and nodal points.
 Denote
$$\mathfrak D_i(R_0)=\Sigma_{i0}\cup \{(s_{i},t_{i})|\;|s_{i}| \leq R_0\},\;\;i=1,\;2,$$$$\mathfrak D(R_0)=\mathfrak D_1(R_0)\bigcup \mathfrak D_2(R_0).$$
 Choose $R_{0}$ such that
\begin{equation}
 \sum_{i=1}^{2}\widetilde{E}(u_i; |s_{i}|\geq \frac{R_{0}}{2})\leq \frac{\min\{\hbar,T\}}{8},\;\;\;\;\;
 4\mathcal{C}_{1}\fc_{1}^{-1}e^{-\frac{\fc R_{0}}{4}}\leq \frac{\hbar_{1}}{8},\;\;\;\;
 \end{equation}

We let $\frac{lr}{Tk}>>4R_0.$ 
\v\n
\begin{lemma}\label{approximation_estimates}
There exists $\epsilon>0$ such that for any $J$-holomorphic map $v:\Sigma_{(r)}\rightarrow M_{(r)}^{+}$ with  $v=\exp_{u_{(r)}}(h) ,$ if
\begin{equation}\label{neighbor_gluing}
|h\mid_{\mathfrak D_i(R_0)}|_{C^{1}}\leq   \epsilon \;\;for\;\;i=1,2,
\end{equation}
then for any $0<\alpha<\fc$
$$
\|h\|_{1,p,\alpha,r}\leq \frac{\hbar_{1} }{4},
$$
where $\fc$ is the constant in Theorem \ref{tube_exponential_decay}.
\end{lemma}
{\bf Proof.} Note that
 $$\widetilde{E}(v; {|s_{i}|\geq R_0/2})=\int_{ R_0/2 \leq s_{1} \leq \frac{2lr}{Tk}-R_0/2}v^{*} d \lambda =\int_{v(\frac{2lr}{Tk}-R_0/2,S^{1})}\lambda-\int_{v( R_0/2,S^{1})}\lambda. $$
Denote $s_{0}=R_{0}/2.$ A direct calculation gives us
\begin{align}
&\left|\int_{v(s_{0},S^{1})}\lambda-\int_{u_{1}(s_{0},S^{1})}\lambda\right|
=\left|\int_{S^1}\lambda(v_t)(s_{0},t)dt
-\int_{S^1}\lambda((u_1)_t)(s_{0},t)dt\right| \nonumber\\
&\leq\left|\int_{S^1}\lambda(P_{u_{1},v}(u_1)_t))(s_{0},t)dt-\int_{S^1}\lambda((u_1)_t)(s_{0},t)dt\right|+
\left|\int_{S^1}\lambda(d\exp_{u_1}h_{t})(s_{0},t)dt\right|\nonumber\\
&\leq C|h\mid_{\mathfrak D_{1}(R_0)}|_{C^{1}}
\end{align}
Then
$$\widetilde{E}(v; {|s_{i}|\geq R_0/2})\leq  \widetilde{E}(u_{i}; {|s_{i}|\geq R_0/2})+C\sum|h\mid_{\mathfrak D_{i}(R_0)}|_{C^{1}}.$$
It follows that
$$ \widetilde{E}(v; \Sigma_{(r)}-\mathfrak D(R_0/2) )\leq \frac{\min\{\hbar,T\}}{4}, $$
  when $0<\epsilon<\frac{\min\{\hbar,T\}}{8C}$. Then by Theorem \ref{tube_exponential_decay} we have
\begin{equation}|\nabla v\mid_{ \Sigma_{(r)}-\mathfrak D(R_0 )}|  \leq \mathcal{C}_{1} e^{-\fc  R_{0}/2 }
\end{equation}
for all $R_0 \leq s_1 \leq \frac{lr}{Tk}-R_{0}.$ Together with \eqref{neighbor_gluing} we have
$$  \|h\mid_{\mathfrak D(R_0)}\|_{1,p,\alpha}\leq \frac{\hbar_{1} }{8}$$
when $\epsilon$ small enough.
 Then lemma follows. $\;\;\;\;\Box$
\v\n

The lemma can be generalized to the case of gluing several Riemann surfaces and several nodal points.

\v

\section{Gluing theory--Regularization}\label{gluing_regularization}

\subsection{Local regularization}\label{Local regularization}
We fist discuss the top strata, then consider lower stratas.
\vskip 0.1in
\noindent
\subsubsection{Top strata.}
\vskip 0.1in
\noindent
 Let $b = (u,(\Sigma,j);{\bf y},{\bf p}) \in \overline{\mathcal{M}}_{A}(M^{+},C;g,m+\nu,(\mathbf{k},\mathfrak{e}))$.
Here $\Sigma$ is a smooth Riemann surface of genus $g$, $j$ is a complex structure (including marked points ), which
is standard near each puncture point. We introduce the holomorphic cylindrical coordinates on $\Sigma$ near each puncture points
$p_i$.
\vskip 0.1in
\noindent
We discuss two different cases:
\vskip 0.1in
\noindent
{\bf (1.1) $\Sigma$ is stable.}
\vskip 0.1in
\noindent
Denote by $O_{j}$ a neighborhood of complex structures on $(\Sigma,j)$.
Note that we change the complex structure in the compact set $\mathfrak D(R_0)$
of $\Sigma$ away from the puncture points.
A neighborhood ${\mathcal U}_b$ of $b$ can be described as
$$O_{j}\times \{\exp_u(h + \hat{h}_0); h\in {\mathcal C}(\Sigma;u^{\ast}TM^{+}),
h_0\in  \bh , \|h\|_{1,p,\alpha} + |h_0| <\epsilon \}/stab_b,$$
where $\bh= \bigoplus_{j=1}^\nu (T_{p_{j}}(\mathcal{F}_{\fe_j})\oplus (span\{\frac{\p}{\p a}\}) ).$
\vskip 0.1in
\noindent
{\bf (1.2) $\Sigma$ is unstable.}
\vskip 0.1in
\noindent
In this case the automorphism group $Aut_{\Sigma}$ is infinite. One must construct a slice of the action $Aut_{\Sigma}$ and
construct a neighborhood ${\mathcal U}_b$ of $b$. This is done in \cite{CLW}. We omit it here.
\vskip 0.1in
\noindent
There is a neighborhood $U$ of $b$  such that $\E$ is trivialized over $U$,
more precisely , $P(u,v)$ gives a
isomorphism $\E|_u \rightarrow \E|_v$ for any $v\in U$.
Choose $K_b \subset \E|_u$ to be a finite dimensional subspace such that
\begin{equation}\label{regularization_operator}
K_b + image D_u = \E|_u,
\end{equation}
where $D_u:W^{1,p,\alpha}(u^*TM^+)\rightarrow \E|_{u}.$
Without loss of generality we may assume that every element of $K_b$ is smooth along $u$ and supports in the common compact
subset $\mathfrak D(R_0)$.
Define a thickned Fredholm system $(K_b\times U, K_b\times \E|_U, S_b)$, where
$${\mathcal S}(\kappa,v) = \bar{\partial}_{J}v + P_{u;v}\kappa\in \E_v,\;\;\;\forall\;(\kappa,v)\in K_b\times U. $$
By \eqref{regularization_operator}, there exists a smaller neighborhood of $b$, still denoted by $U$, such that $D{\mathcal S}_{(\kappa,v)} :K_{b}\times W^{1,p,\alpha}
\rightarrow L^{p,\alpha}$  is surjective for any $(\kappa,v)\in K_b\times U$.
Then for each $b'=(\kappa,v)\in K_b\times U$ there exists a right inverse
$$Q_{b'} : L^{p,\alpha}\rightarrow K_{b}\times W^{1,p,\alpha}$$ such that
$\|Q_{b'}\|\leq C_{1}$.
Obviously,  $D{\mathcal S}_{b'} :K_{b}\times \mathcal{W}^{1,p,\alpha} \rightarrow L^{p,\alpha}$  is also surjective, $Q_{b'} : L^{p,\alpha}\rightarrow K_{b}\times \mathcal{W}^{1,p,\alpha}$ is also a right inverse of $D{\mathcal S}_{b'}$.

 \v
 A pair $(\kappa,v)\in K_b\times U$ is called a perturbed $J$-holomorphic map if
$$\bar{\partial}_{J}v + P_{u;v}\kappa =0.$$

\v

\vskip 0.1in
\noindent
\subsubsection{ Lower strata.}\label{Lower strata}

\vskip 0.1in
\noindent

We shall consider two cases for simplicity, the discussion for general cases are the same.
Let $D$ be a strata whose domain has two components $(\Sigma_1,j_1)$ and $(\Sigma_2,j_2)$
joining at $p$.

\vskip 0.1in
\noindent
{\bf Case 1}. Let $b=(u_1,u_2; \Sigma_1 \wedge \Sigma_2,j_1,j_2)$,
where $(\Sigma_1,j_1)$ and $(\Sigma_2,j_2)$ are smooth Riemann surfaces of genus $g_1$ and $g_2$ joining at $p$ and $u_i: \Sigma_i \rightarrow M^+$ are $J$-holomorphic maps with $u_1(p)=u_2(p).$ Suppose that both $(u_i,\Sigma_i,j_i)$ are stable. A neighborhood ${\mathcal U}^D_b$ of $b$ in the strata $D$ can be described as
$$O_{j_1}\times O_{j_2}\times
\left\{\left(\exp_{u_1}(h_1 + \hat{h}_{10}),\exp_{u_2}(h_2 + \hat{h}_{20})\right)|
h_i\in W^{1,p,\alpha}(\Sigma;{u_i}^{\ast}TM^{+})\right.,$$
$$\left.h_{10}=h_{20}\in T_{u_1(p)}M^+, \|h_i\|_{1,p,\alpha} + |h_0| <\epsilon \right\}/stab_u.$$
The fiber of ${\mathcal E}_{D}(M^+,g,m+\nu)\rightarrow {\mathcal B}_{D}(M^+,g,m+\nu)$
at $b$ is $L^{p,\alpha}(u_1^{\ast}TM^+\otimes \wedge^{0,1})\times
L^{p,\alpha}(u_2^{\ast}TM^+\otimes \wedge^{0,1}).$

\vskip 0.1in
\noindent
{\bf Case 2}. Let $b=(u_1,u_2; \Sigma_1 \wedge \Sigma_2,j_1,j_2)$,
where $(\Sigma_1,j_1)$ and $(\Sigma_2,j_2)$ are smooth Riemann surfaces of genus $g_1$ and $g_2$ joining at $q$, $u_1: \Sigma_1 \rightarrow M^+$ and
$u_2 :\Sigma_2\rightarrow {\mathbb{R}}\times \widetilde{M}$
are $J$-holomorphic maps such that
$u_1$ and $u_2$ converge to a same $kT$-periodic orbit $x(kTt)\subset \mathcal F_{j}$ for some $j$ when $z_1$ and $z_2$
converge to $q$. Suppose that both $(u_i,\Sigma_i,j_i)$ are stable.
\v
We first describe the neighborhood of $u_2 :\Sigma_2\rightarrow {\mathbb{R}}\times \widetilde{M}$ in the space of $\mathcal{W}^{1,p,\alpha}$-maps. The $R$-action on $\bR\times \widetilde{M}$ and $S^1$-action on every periodic orbit naturally induce the actions on $W^{1,p,\alpha}(\Sigma_{2},\bR\times \widetilde{M})$, which we need mod.
Write $u(s,t)=(a(s,t),\tilde{u}).$ For any $C $, the $\mathbb R$-action is defined as
$$ \mathbb R_{C}\circ u(s,t)=(C+a(s,t),\tilde{u}).$$

We construct a local R-action slice as follows. We fix a point $p\in \Sigma_{2}$ and choose coordinates $a$ such that the $u_2(p)=(0, \star)$. Let $h$ be a vector field along $u_2$. We write
$$h=b_1\frac{\p}{\p a} + \tilde{h},$$
where $\tilde{h}\in W^{1,p,\alpha}(\tilde{u}_2^*(T\tilde{M}))$ is a vector field along $\tilde{u}_2$. Then $h= b_1\frac{\p}{\p a} + \tilde{h}$ can be considered in a natural way as a $\mathbb R$ equivalent vector field along $\mathbb R\circ u$. We can construct a equivalent tubular neighborhood around the orbit $\mathbb R\circ u$. The point-wise exponential map give us a Banach manifold structure on the space of $\mathbb R$ equivalence class of $W^{1,p,\alpha}(\Sigma_{2},\bR\times \widetilde{M})$. To simplify notations we will use $W^{1,p,\alpha}(\Sigma_{2},\bR\times \widetilde{M})$ to denote the space of $\mathbb R$ equivalence class if no danger of confusion.
\v
On the other hand, we need mod $S^1$-action on every periodic orbit. Let $x$ be a periodic orbit. We choose a local pseudo-Darboux coordinate system $(a ,\vartheta ,{\bf w})$, near $x$. If we choose another origin in $S^1$ we get another local pseudo-Darboux coordinate system $(a' ,\vartheta' ,{\bf w})$, which differ by a canonical coordinates transformation \eqref{canonical}.  The vector $h=b_1\frac{\p}{\p a} + \tilde{h}$ is independent of the  local pseudo-Darboux coordinates, so independent of the choice of the origin.
\vskip 0.1in
\noindent
 For simplicity we assume that $\Sigma_{i},i=1,2$ has no puncture points. A neighborhood ${\mathcal U}^D_b$ of $b$ in $D$ can be described as
$$O_{j_1}\times O_{j_2}\times
\left\{(\exp_{u_1}(h_1 + \hat{h}_{10}),\exp_{u_2}(h_2 + \hat{h}_{20}))|
h_1\in W^{1,p,\alpha}(\Sigma;{u_1}^{\ast}TM^{+}),\right.$$
$$\left.h_2\in W^{1,p,\alpha}(\Sigma;{u_2}^{\ast}T({\mathbb{R}}\times \widetilde{M})),
h_{10}=h_{20}\in {\bh} ,\;\;\;\; \Sigma_{i=1}^2\|h_i\|_{1,p,\alpha} + |h_{10}|
<\epsilon \right\}/stab_u,$$
where $\bh=  T_{q}(\mathcal{F}_{j})\oplus (span\{\frac{\p}{\p a}\}).$
The fiber of ${\mathcal E}_{D}\rightarrow {\mathcal B}_{D}$
at $b$ is $L^{p,\alpha}(u_1^{\ast}TM^+\otimes \wedge^{0,1})\times
L^{p,\alpha}(u_2^{\ast}({\mathbb{R}}\times \widetilde{M})\otimes \wedge^{0,1}).$
\vskip 0.1in
\noindent

We use the gluing argument to describe the neighborhoods ${\mathcal U}_b$ of $b$ in ${\mathcal B}_{A}(M^{+},C;g,m+\nu,(\mathbf{k},\mathfrak{e}))$.
Denote $Glu_{T_0}=\{(r,\tau)|0\leq \tau<1, T_0\leq r\leq \infty\}$. Then
$${\mathcal U}_b =\bigcup_{T_0\leq r<\infty} B(u_{(r)},\epsilon) $$
where $B(u_{(r)},\epsilon)=\{ u\in{\mathcal B} |\; u=\exp_{u_{(r)}}(h+\hat h_{0}),\;\;\;\|h\|_{1,p,\alpha}+|h_{0}|< \epsilon \}.$
\vskip 0.1in
\noindent
Choose $K_b=(K_1,K_2) \subset \E|_b=(\E_{u_1},\E_{u_2})$ to be a finite dimensional subspace such that
$$K_1 + image D_{u_1} = \E_{u_1},\;\;\;K_2 + image D_{u_2} = \E_{u_2}.$$
where $D_{u_{i}}:W^{1,p,\alpha}\rightarrow \E_{u_{i}}.$ We may assume that every element of $K_i$ is smooth along $u_i$ and supports in the compact subset $\mathfrak D(R_0)$ of
$\Sigma_i$. For any $(r,\tau)\in Glu_{T_0}$ there is a neighborhood $U$ of $(\Sigma_{(r)}, u_{(r)})$ in ${\mathcal B}_{A}(M^{+},C;g,m+\nu,(\mathbf{k},\mathfrak{e}))$ such that $\E$ is trivialized over $U$.
We choose $T_0>R_0$ so large that for $r>T_0$
$$K_b\subset \E|_{U}.$$
Note that $u_{(r)}|_{|s_i|\leq R_0}=u_{i}|_{|s_i|\leq R_0}$ and $supp\;(K_b)\subset \{|s_i|\leq R_0\}$.
Then $P(u_{(r)},v)|_{|s_i|\leq R_0}$ is independent of $r$. Hence $\kappa_{u}$ is well defined for any $u\in B(u_{(r)},\epsilon).$ We can naturally identify $K_{b}$ and $K_{b_{(r)}}$.
 We define a thickned Fredholm system and regularization equation as in {\bf (1)}. Then
there exists a neighborhood $U$ of $b$, such that $D{\mathcal S}_{b'}:K_{b}\times W^{1,p,\alpha}
\rightarrow L^{p,\alpha}$ is surjective for any $b'=(\kappa,v)\in K_b\times U $.
Then   $D{\mathcal S}_{b'}:K_{b}\times \mathcal{W}^{1,p,\alpha} \rightarrow L^{p,\alpha}$  is also surjective.

\vskip 0.1in
\noindent

\subsection{Global regularization}\label{global_r}
Denote
$$
U_b^{\epsilon}=\{v\in {\mathcal{B}} | \; v=\exp_{u}(h+\hat h_{0}),\;\;\;\|h\|_{1,p,\alpha}+|h_{0}|<\epsilon\}.
$$
As $ \overline{\mathcal{M}}_{A}(M^{+},C;g_1,m_1+\nu,({\bf k},\fe))$ is compact, there exist finite points $b_i$, $1\leq i \leq n$ such that the collection $\{U_{b_i}^{\epsilon}\}$ is an open cover of $ \overline{\mathcal{M}}_{A}(M^{+},C;g_1,m_1+\nu,({\bf k},\fe))$.
For any $I\subset \{1,2,...,n\}$, setting
$$K_I=\prod_{i\in I}K_{b_i},\;\;\;S_I((\kappa_i)_{i\in I},v)=\overline{\p}_{J}(v) + \sum_{i\in I}\beta_{b_i}P_{u,v}(\kappa_i)\in \E_v.$$
One can construct a global regularization
 $(C_I, F_I,S_I)$ of the original Fredholm system
$(\mathcal{B}, \mathcal{E}, \overline{\p}_{J})$. From the global regularization we can obtain that (see \cite{CT,CLW}):
\begin{lemma} \label{prop_3.4}
 There exists  a finite dimensional virtual  orbifold system  for $(\mathcal{B}, \mathcal{E}, \overline{\p}_{J})$ which is   a collection of triples
 \[
  \{(U_I,   E_I, \sigma_I)|  I\subset \{1, 2, \cdots, n\} \}
  \]
  indexed by a   partially ordered
set $( I =2^{\{1, 2, \cdots, n\}}, \subset )$, where
  \begin{enumerate}
\item $\{ U_I |  I\subset \{1, 2, \cdots, n\}  \}$ is a finite dimensional  proper  \'etale virtual groupoid, where $U_{I}=S^{-1}_{I}(0),$
\item $\{ E_I\}$ is a  finite rank virtual vector bundle over $\{ U_I\}$,
\item $\{\sigma_I\}$ is a section  of the  virtual vector bundle $\{ E_I\}$ whose zeros $\{ \sigma_I^{-1} (0)\}$  form a cover of $ \overline{\mathcal{M}}_{A}(M^{+},C;g_1,m_1+\nu,({\bf k},\fe))$.
 \end{enumerate}
\end{lemma}

Let $ U_{I,\epsilon}$ be the set of points $((\kappa_i)_{i\in I},v) \in U_I$ such that
$ \sum_{i\in I}\|\kappa_i\|_{p,\alpha}\leq \epsilon,\;\mbox{where}\;I\subset \{1,2,...,n\}.$
We define $U_{I,\epsilon}({\mathbb{R}}\times
\widetilde{M})$, $U_{I,\epsilon}(M^{+}\cup ({\mathbb{R}}\times \widetilde{M}))$ in the same way.


By the same argument of Theorem \ref{compact_moduli_space} we conclude that:

\v
\begin{theorem}  $\{\Sinverse\}$ is a compact virtual manifold.
\end{theorem}

\section{Gluing theory-- analysis estimates}\label{gluing_analysis}

\subsection{Gluing theory for 1-nodal case}

We consider the case of gluing one node in Subsection \S\ref{relative_node_M_R_M}. The general cases are similar.
\v\n
Recall that we have a global regularization $\{(\mbf C, \mbf F, \mbf S)\}$. Consider the regular Fredholm system $(C_I, F_I,S_I)$. Let $b= ((u_{1},u_{2}),\Sigma_{1}\wedge\Sigma_{2};
j_{1} ,j_{2} ) \in C_I$, where $(\Sigma_1,j_1)$ and $(\Sigma_2,j_2)$ are smooth Riemann surfaces of genus $g_1$ and $g_2$ joining at $q$ and
$u_1: \Sigma_1 \rightarrow M^+$, $u_2: \Sigma_2 \rightarrow {\mathbb{R}}\times \widetilde{M}$ are $J$-holomorphic maps such that
$u_i(z)$ converge to the same $kT$-periodic orbit $x(kTt)\subset \mathcal{F}_{j}$ for some $j$ as $z\rightarrow p$.
\v\n

For any $(\kappa,h,h_{0})\in Ker D \mathcal {S}_{b},$
 where $h\in W^{1,p,\alpha}(\Sigma;u^{*}TN),$  we define $$ \|(\kappa,h)\|_{1,p,\alpha}=\|\kappa\|_{p,\alpha}+\|h\|_{1,p,\alpha},\;\;\;\;
\|(\kappa,h,h_{0})\|=\|(\kappa,h)\|_{1,p,\alpha}+|h_{0}|.$$
 For any $(\kappa,h_{(r)})\in Ker D \mathcal {S}_{b_{(r)}},$ we define
$$\|(\kappa,h_{(r)})\|=\|\kappa\|_{p,\alpha}+\|h_{(r)}\|_{1,p,\alpha,r}.
$$

By using the exponential decay of $u_i$ one can
easily prove that $u_{(r)}$ are a family of approximate
$J$-holomorphic map, precisely the following lemma holds. \vskip
0.1in \noindent
\begin{lemma} For any $r>R_{0},$  we have
\begin{equation}\label{estimate_p_u_r}
\|\bar{\partial}_{J}(u_{(r)})\|_{p,\alpha,r}\leq
Ce^{-(\mathfrak{c} -\alpha)r} .
\end{equation}  The constants C in the
above estimates are independent of $r$.
\end{lemma}

\subsection{Estimates of right inverse}

\begin{lemma}\label{right_inverse_after_gluing}
Let  $D{\mathcal
S}_{b} :K_{b}\times \mathcal{W}^{1,p,\alpha}\rightarrow L^{p,\alpha}$ be a Fredholm operator defined in  section \S\ref{gluing_regularization}.
 Suppose that $D{\mathcal S}_{b}|_{K_{b}\times  {W}^{1,p,\alpha}}:K_{b}\times  {W}^{1,p,\alpha}\rightarrow L^{p,\alpha}$ is surjective.
Denote by $Q_b:L^{p,\alpha}\rightarrow K_{b}\times W^{1,p,\alpha}$ a right inverse of $D{\mathcal
S}_{b} .$  Then $D{\mathcal
S}_{b_{(r)}}$ is surjective for $r$ large enough. Moreover, there are a right inverses $Q_{b_{(r)}}$
such that
\begin{equation}
\label{right_verse}
D{\mathcal S}_{b_{(r)}}Q_{b_{(r)}}=Id \end{equation}
\begin{equation}
\label{right_estimate}
\|Q_{b_{(r)}}\|\leq  {C}
\end{equation}
 for some constant $C>0$ independent of $ r $.
  \end{lemma}
  \vskip
0.1in
\vskip
0.1in

\noindent {\bf Proof:} We first construct an approximate right inverse $Q'_{b_{(r)}}$ such that the following estimates holds
\begin{eqnarray}
\label{approximate_right_inverse_estimate_1}
\|Q'_{b_{(r)}}\|\leq C_1 \\
\label{approximate_right_inverse_estimate_2}
\|D{\mathcal S}_{b_{(r)}}\circ Q'_{b_{(r)}}-Id\|\leq \frac{1}{2}.
\end{eqnarray}
 Then the operator $D{\mathcal S}_{b_{(r)}}\circ Q'_{b_{(r)}}$ is invertible and a right inverse $Q_{b_{(r)}}$  of $D{\mathcal S}_{b_{(r)}}$ is given by
 \begin{equation}
 \label{express_right_inverse}
Q_{b_{(r)}}=Q'_{b_{(r)}}( D{\mathcal S}_{b_{(r)}}\circ Q'_{b_{(r)}})\inv
 \end{equation}
Denote $\beta_1=\beta(3/2-\frac{Tks_1}{lr}).$ Let $\beta_{2}\geq 0$ be a smooth function such that $\beta_2^2=1-\beta_1^2.$ Given
$\eta\in L_{r}^{2,\alpha}$, we have a pair $(\eta_1,\eta_2)$, where
$$\eta_1=\beta_{1}\eta,\;\;\;\; \eta_2=\beta_{2}\eta.$$ Let
$Q_b(\eta_1,\eta_2)=(\kappa_b,h).$ We may write $h$ as
$(h_1,h_2)$, and define
\begin{equation}
\label{def_h_r}
h_{(r)}=
h_1\beta_1+ h_2\beta_{2}.
\end{equation}
Note that on $\{\frac{lr}{2Tk}\leq s_1\leq \frac{3lr}{2Tk}\},$ $\kappa=0$ and on $\{|s_i|\leq \frac{lr}{2Tk}\}$ we have $u_{(r)}=u_i$, $\kappa_{(r)}=\kappa_b$,
so  along $u_{(r)}$ we have $\kappa_{(r)}=\kappa_b$. Then we define
\begin{equation}
\label{def_approximate_right_inverse}
Q_{b_{(r)}}^{\prime}\eta = (\kappa_{(r)},
h_{(r)})=(\kappa_b,
h_{(r)}).\end{equation}
Since $|\beta_1|\leq 1$ and $|\frac{\partial \beta_1}{\partial s_1}|\leq \frac{CTk}{lr},$ \eqref{approximate_right_inverse_estimate_1} follows from $\|Q_{b}\|\leq C.$
 We prove \eqref{approximate_right_inverse_estimate_2}. Since $\kappa_b+D_{u}h=\eta$ we have
\begin{equation}\label{app_DS_right_inverse}
 D{\mathcal S}_{b_{(r)}}\circ Q'_{b_{(r)}}\eta=\eta,\;\;\;\;\;\;|s_{i}|\leq \frac{lr}{2Tk}.
 \end{equation}
 It suffices to estimate the left hand side in the left annulus $\frac{lr}{2Tk}\leq |s_i|\leq \frac{3lr}{2Tk}.$
 Note that in this annulus
 $$\beta_1^2+\beta_2^2=1,\;\;\kappa_{b}=0,\;\;\; D_{u_{i}}h_{i}=\eta_{i},$$
 $$\beta_1 D_{u_{1}}h_{1}+\beta_2 D_{u_{2}}h_{2}=(\beta_1^2+\beta_2^2)\eta. $$
 Since near $u_{1}(p)=u_{2}(p)$ ( or near the periodic orbit $x(kTt)$), $D_{u_{i}}=\bar{\partial}_{J_{0}}+S_{u_{i}}$, we have
\begin{eqnarray}
&&D \mathcal{S}_{b_{(r)}}\circ Q'_{b_{(r)}}\eta-(\beta_1^2+\beta_2^2)\eta = \kappa_{b_{(r)}} + D_{u_{(r)}}h_{(r)}-(\beta_1^2+\beta_2^2)\eta \nonumber\\\label{approximate_difference}
&=&(\bar{\partial}\beta_{1}) h_{1} +\beta_{1}(S_{u_{(r)}}-S_{u_{1}})h_{1}+(\bar{\partial}\beta_{2}) h_{2} +\beta_{2}(S_{u_{(r)}}-S_{u_{2}})h_{2}.
\end{eqnarray}
By the exponential decay of $S$ and $\beta_1^2+\beta_2^2=1$ we get
\begin{align}
& \left\|D \mathcal{S}_{b_{(r)}}\circ Q'_{b_{(r)}}\eta-\eta\right\|_{p,\alpha,r} =  \left\|D \mathcal{S}_{b_{(r)}}\circ Q'_{b_{(r)}}\eta-(\beta_{1}^2+\beta_{2}^2)\eta\right\|_{p,\alpha,r} \nonumber \\
&\leq \frac{C_1}{r}\left(\|h_1\|_{p,\alpha}+ \|h_2\|_{p,\alpha} \right)  \leq  \frac{C_2}{r} \|\eta\|_{p,\alpha,r}
\end{align}
In the last inequality we used that  $\|Q_{b}\|\leq C$ and $(h_{1},h_{2})=\pi_{2}\circ Q_{b}(\eta_{1},\eta_{2}).$
 Then \eqref{approximate_right_inverse_estimate_2} follows by choosing $r$ big enough. The estimate \eqref{approximate_right_inverse_estimate_2} implies that
\begin{equation}\label{app_DS_right_inverse_bound}
\frac{1}{2}\leq \|D{\mathcal S}_{b_{(r)}}\circ Q'_{b_{(r)}}\|\leq \frac{3}{2}.
\end{equation}
Then \eqref{right_estimate} follows.  $\Box$
\v

\subsection{Isomorphism between $Ker D \mathcal{S}_{b}$ and $Ker D \mathcal{S}_{b_{(r)}}$}
\v
Put
$$E_1:=Ker D \mathcal{S}_{u_1},\;\;E_2:=Ker D \mathcal{S}_{u_2}, \;\;\;\;\bh=  T_{q}(\mathcal{F}_{j})\oplus (span\{\frac{\p}{\p a}\}),$$
Denote
$$Ker D \mathcal{S}_{b}:=E_1\bigoplus_{\bh}E_2.$$

\vskip 0.1in
\noindent
For a fixed gluing parameter $(r)=(r,\tau)$
we define a map
$I_r: Ker D \mathcal{S}_{b}\longrightarrow Ker D \mathcal{S}_{b_{(r)}}.$ For any $(\kappa,h,h_{0})\in Ker D \mathcal {S}_{b},$
 where $h\in W^{1,p,\alpha}(\Sigma;u^{*}TN),$  we write $h=(h_1  ,h_2  ),$ and define
\begin{equation}
\label{definition_h_ker}
h_{(r)}= \hat{h}_{0} + h_{1}\beta_1+h_{2}\beta_2,
\end{equation}
\begin{equation}\label{definition_I}
I_{r}(\kappa,h,h_{0})=(\kappa,h_{(r)})-Q_{b_{(r)}}\circ D\mathcal{S}_{b_{(r)}}(\kappa,h_{(r)}).
\end{equation}

\begin{lemma}$I_r$ is an isomorphisms for $r$ big enough.
\end{lemma}
 \vskip 0.1in
\noindent {\bf Proof:} The proof is basically a similar gluing
argument as in \cite{D}. The proof is
devides into 2 steps. \vskip 0.1in \noindent {\bf Step 1}.  We define a map $I'_{r}:Ker D\mathcal{S}_{b_{(r)}}\longrightarrow Ker D\mathcal{S}_{b} $
 and show that $I'_{r}$ is injective for $r$ big enough. For any $(\kappa,h)\in Ker D\mathcal{S}_{b_{(r)}}$ we denote by $h_{i}$ the restriction of $h$
 to the part $|s_{i}|\leq \frac{lr}{kT} +\frac{1}{\alpha},$ we get a pair $(h_{1},h_{2}).$ Let
 \begin{equation}
 h_{0}=\int_{S^1}h\left(\frac{lr}{kT},t\right)dt.
 \end{equation}
  We denote
  $$\beta [h] =\left((h_{1}- \hat h_{0})\beta\left(\frac{\alpha lr}{kT} +1-\alpha s_1\right) + \hat h_{0},\;\;(h_{2}- \hat h_{0})\beta\left(\frac{\alpha lr}{kT}+1+\alpha s_2\right) + \hat h_{0}\right)$$
and define
$I'_{r}:Ker D\mathcal{S}_{b_{(r)}}\longrightarrow  Ker D\mathcal{S}_{b}$ by
\begin{equation}
\label{definition_I'}
I'_{r}(\kappa,h)= (\kappa,\beta [h])-Q_{b}\circ D\mathcal{S}_{b}(\kappa,\beta [h]),
\end{equation}
 where  $Q_{b}$ denotes the right inverse of $D\mathcal{S}_{b}|_{K_{b}\times W^{1,p,\alpha}}: {K}_{b } \times  W^{1,p,\alpha}\rightarrow L^{p,\alpha}.$  Since
$D\mathcal{S}_{b}\circ Q_{b}=D\mathcal{S}_{b}|_{K_{b}\times W^{1,p,\alpha}}\circ Q_{b}=I,$ we have
$I'_{r}(Ker D\mathcal{S}_{b_{(r)}})\subset  Ker D\mathcal{S}_{b}.$
 \v
Since $\kappa$ and $D_{u} (\beta(h-\hat{h}_{0}))$ have compact support and $S_{u}\in L^{p,\alpha}$, we have
$D\mathcal{S}_{b}(\kappa,\beta [h]) \in L^{p,\alpha}.$
Then $Q_{b}\circ D\mathcal{S}_{b}(\kappa,\beta [h])\in K_{b}\times W^{1,p,\alpha}.$

\v
Let $(\kappa,h)\in Ker D\mathcal{S}_{b_{(r)}}$ such that $I'_{r}(\kappa,h)=0.$ Since $\beta(h-\hat{h}_{0})\in W^{1,p,\alpha}$ and $Q_{b}\circ D\mathcal{S}_{b}(\kappa,\beta [h])\in K_{b}\times W^{1,p,\alpha},$ then  $I'_{r}(\kappa,h)=0$ implies that $h_{0}=0.$ From \eqref{definition_I'} we have
$$\|I'_{r}(\kappa,h) -  (\kappa, \beta h) \|_{1,p,\alpha}
 \leq C_1  \| \kappa + D_u(\beta h)\|_{p,\alpha}  $$
$$= C_1 \left \| \kappa  + \beta \left(D_u h
  + D_{u_{(r)}} h+\kappa- D_{u_{(r)}} h- \kappa\right) + (\bar{\partial}\beta) h \right\|_{p,\alpha} $$
Since $(\kappa,h)\in Ker D\mathcal {S}_{b_{(r)}}$, we have $\kappa+D_{u_{(r)}}h=0.$  We choose $\frac{lr}{2kT}>R_{0}$.
As $\kappa|_{|s_i|\geq R_{0}}=0$ and $\beta |_{|s_{i}|\leq \frac{lr}{kT}}=1$ we have $\kappa = \beta\kappa$.
Therefore
\begin{eqnarray*}
 \left\|I'_{r}(\kappa,h)-(\kappa,\beta h)\right\|_{1,p,\alpha}
 \leq    C_1  \|(\bar{\partial}\beta)  h +
\beta  (S_u - S_{u(r)}) h  \|_{p,\alpha} .
\end{eqnarray*}
Note that
$$S_u = S_{u(r)} \;\;\;\;if \; \; s_1\leq \; \frac{lr}{2kT}, \; or \;\; s_2 \geq -\frac{lr}{2kT}. $$
By exponential decay of $S$ we have
$$\|(S_u - S_{u(r)})\beta h\|_{p,\alpha} \leq Ce^{-\mathfrak{c}\frac{ lr}{2kT}}\|\beta h\|_{1,p,\alpha}$$
for some constant $C>0$.
Since $(\bar{\partial}\beta(\frac{\alpha lr}{kT} + 1 -\alpha s_1))h_{1}$ supports in
$\frac{lr}{kT}\leq s_1 \leq \frac{lr}{kT}+ \frac{1}{\alpha}$, and over this part
$$|\bar{\partial}\beta(\frac{\alpha lr}{kT}+ 1 -\alpha s_1)|\leq 2|\alpha|$$
$$\beta(\frac{\alpha lr}{kT}+ 1 + \alpha s_2)= 1,\;\;\;e^{2\alpha|s_1|}\leq e^4e^{2\alpha|s_2|},$$
we obtain
$$\|(\bar{\partial}\beta(\frac{\alpha lr}{kT} + 1 -\alpha s_1))h_{1}\|_{p,\alpha}\leq2
|\alpha|e^4 \|h_2\|_{p,\alpha} \leq 2|\alpha|e^4 \|\beta h\|_{p,\alpha}.$$
Similar inequality for $(\bar{\partial}\beta(\frac{\alpha lr}{kT} + 1 + \alpha s_2))h_{2}$
also holds. So we have
$$\|(\bar{\partial}\beta)h\|_{p,\alpha} \leq 4|\alpha|e^4 \|\beta h\|_{1,p,\alpha}.$$
Hence
\begin{equation}\label{delta_I'}
\|I'_{r}(\kappa,h) - (\kappa, \beta h)\|_{1,p,\alpha} \leq C_3(|\alpha| +
e^{- \frac{\mathfrak{c}   lr}{2kT}})\|\beta h\|_{1,p,\alpha} \leq   1/2  \|\beta h\|_{1,p,\alpha}  \end{equation}
for some constant $C_3>0$, here we choosed $0<\alpha<\frac{1}{4C_{4}}$ and  $r$   big enough  such
that $\frac{lr}{kT}>\frac{1}{\alpha}$ and  $C_3 e^{- \frac{\mathfrak{c}   lr}{2kT}} < 1/4$.

Then $I'_{r}(\kappa,h)=0$ and
\eqref{delta_I'} gives us $$\|\kappa\|_{p,\alpha}=0 , \;\;  \|\beta h\|_{1,p,\alpha}=0 .$$
It follows that $\kappa = 0, \;\; h=0$. So $I'_r$ is injective.
\vskip 0.1in
\noindent

{\bf Step 2}.
Since $\|Q_{b_{(r)}}\|$ is uniformly bounded, from \eqref{definition_I}  and \eqref{right_estimate}, we have
$$\|I_r((\kappa,h),h_{0}) - (\kappa, h_{(r)})\|_{1,p,\alpha,r}\leq C\|D\mathcal {S}_{b_{(r)}} (\kappa, h_{(r)})\|.$$
By a similar culculation as in the proof of Lemma \ref{right_inverse_after_gluing} we obtain
\begin{equation}\label{delta_I}
\|I_r((\kappa,h),h_{0}) - (\kappa, h_{(r)})\|_{1,p,\alpha,r}\leq \frac{C}{r}
(\|h\|_{p,\alpha} + |h_0|).\end{equation}
In particular, it holds for $p=2$. It remains to show that
$\|h_{(r)}\|_{2,\alpha,r}$ is close to $\|h\|_{2,\alpha}$. Denote $\pi$ the
projection into the second component, that is, $\pi((\kappa,h),h_{0})=h$. Then
$\pi(ker D{\mathcal S}_{b})$ is a finite dimentional space. Let $f_i,\;i=1,..,d$
be an orthonormal basis. Then $F=\sum f_i^2e^{2\alpha|s|}$ is an integrable
function on $\Sigma$. For any $\epsilon >0$, we may choose $R_0$ so big that
$$\int_{|s_i|\geq R_0}F \leq \epsilon.$$
Then the restriction of $h$ to $|s_i|\geq R_0$ satisfies
$$\|h|_{|s_i|\geq R_0}\|_{2,\alpha}\leq \epsilon\|h\|_{2,\alpha},$$
therefore
\begin{equation}\label{lower_bound_h_r}
\|h_{(r)}\|_{2,\alpha,r}\geq \|h|_{|s_i|\leq R_0}\|_{2,\alpha} + |h_0|
\geq (1 - \epsilon)\|h\|_{2,\alpha} + |h_0|,
\end{equation}
for $r>R_0$. Suppose that $I_r((\kappa,h),h_{0})=0$. Then \eqref{delta_I} and \eqref{lower_bound_h_r} give us $h=0$ and $h_{0}=0$,
and so $\kappa=0$. Hence $I_r$ is injective.
\vskip 0.1in
\noindent
The {\bf step 1} and {\bf step 2} together show that both $I_{r}$ and $I'_r$
are isomorphisms for $r$ big enough.  $\Box$

\subsection{Gluing maps}

There is a neighborhood $U$ of $b_{(r)}$ in $ {\mathcal B}_{A}(M^+,C;g,m+\nu,({\bk},\fe))$ such that $\E$ is trivialized over $U$,
more precisely, $P_{u_{(r)},v}$ gives a
isomorphism $\E|_{u_{(r)}} \rightarrow \E|_{v}$ for any $v \in U$.

Consider a map
\begin{align*}& F_{(r)}: {K}_{b_{(r)}} \times \mathcal W^{1,p,\alpha}_{r}(\Sigma_{(r)};u_{(r)}
^{\ast}TM^+)
\rightarrow L^{p,\alpha}_{r}(u_{(r)}^{\ast}TM^+\otimes \wedge^{0,1})\\
&F_{(r)}(\kappa,h)=P_{exp_{u_{(r)}}h, u_{(r)}}\left(\bar{\partial}_{J}exp_{u_{(r)}}h
+ \kappa\right).
\end{align*}
Let $(\kappa_{\tau},h_{\tau}):[0,\delta)\longrightarrow {K}_{b_{(r)}} \times \mathcal W^{1,p,\alpha }_{r}(\Sigma_{(r)};u_{(r)}
^{\ast}TM^+)$ be a  curve  satisfying
$$
\kappa_{0}=\kappa,\;\;\;h_{0}=h,\;\;\;\;\frac{d}{d\tau}\kappa_{\tau} \mid _{\tau=0}=\eta,\;\;\;\;\frac{d}{d\tau}h_{\tau} \mid _{\tau=0}=g
$$

  By the same method of \cite{MS} we can prove
\begin{lemma} \label{difference_DF_estimates_gw}
\begin{description}
\item[1.] $d F_{(r)}(0)=D\mathcal{S}_{ b_{(r)} }.$
\item[2.] There is a constant $ \hbar_2>0 $ such that for any $(\kappa,h)\in {K}_{b_{(r)}} \times \mathcal W^{1,p,\alpha}_{r}(\Sigma_{(r)};u_{(r)}
^{\ast}TM^+)$ with $\|(\kappa,h)\|<\hbar_2$ and $\|F_{u_{(r)}}(\kappa,h)\|_{p,\alpha,r}\leq 2\hbar_2,$ the following inequality holds
\begin{equation} \label{difference_of_DF_estimates}
\|dF_{(r)}(\kappa,h)-dF_{(r)}(0)\| \leq \frac{1}{4\mathcal C_2}.
\end{equation}
\end{description}
\end{lemma}
\v

 Let $\hbar_{3}=\min\{\hbar_{1},\hbar_{2} \}.$ Note that $M^{+}=M^{+}_{0}\times ([0,\infty)\times \widetilde M)$. By the definition of the metrics $\langle,\rangle$ (see  \eqref{omega_forms_on_M0} and \eqref{omega_forms_on_tubes}) the curvature  $|Rm |$ is uniform bounded above. Then there exists a constant $C_{2}>0$ such that $|dexp_{p}|\leq C_{2}$ for any $p\in M^+$. For any  $\|(\kappa,h)\|\leq \frac{\hbar_{3}}{8C_{2}\mathcal{C}_{2}}$,
\begin{align}
\|F_{(r)}(\kappa,h)\|_{p,\alpha,r}&=\|\bar{\partial}_{J}exp_{u_{(r)}}h
+ \kappa\|_{p,\alpha,r} \nonumber\\
&\leq \|\overline{\partial}_{J}u_{(r)}\|_{ p,\alpha,r} +\|dexp_{u_{(r)}}\partial_{s} h\|_{p,\alpha,r} +\|dexp_{u_{(r)}}\partial_{t} h\|_{p,\alpha,r}  + \| \kappa  \|_{p,\alpha}\nonumber\\
& \leq \|\overline{\partial}_{J}u_{(r)}\|_{ p,\alpha,r} +C_{2}\|(\kappa,h)\|\leq C_1e^{-(\mathfrak{c}-\alpha) r}+ \frac{\hbar_{3}}{8\mathcal{C}_{2}} \leq \frac{\hbar_{3}}{4\mathcal{C}_{2}}  \label{uniform_bound_of_F}
\end{align}   when  $r$ big enough.
It follows from (2) of Lemma \ref{difference_DF_estimates_gw} and \eqref{uniform_bound_of_F} that
\begin{equation}\label{uniform_bound_difference_dF}
 \|dF_{(r)}(\kappa,h)-D\mathcal{S}_{b_{(r)}}\| \leq \frac{1}{4\mathcal C_{2}}.
\end{equation}
 Then $F_{(r)}$ satisfies the conditions in Lemma \ref{details_implicit_function_theorem}. It follows that for any $(\kappa,h)\in ker D \mathcal{S}_{b_{(r)}}$ with $\|(\kappa,h)\|\leq \frac{\hbar_{3}}{8C_{2}\mathcal{C}_{2}}$ there exists a unique $(\kappa_{v},\zeta)$ such that $F_{(r)}(\kappa_{v},\zeta)=0.$ Since $ {K}_{b_{(r)}} \times \mathcal W^{1,p,\alpha}_{r}(\Sigma_{(r)};u_{(r)}
^{\ast}TM^+)=im Q_{b_{(r)}}+ker D \mathcal{S}_{b_{(r)}} $ and $Q_{b_{(r)}}$ is injective, then there exists a unique smooth map
$$f_{(r)}:   Ker\;D{\mathcal S}_{b_{(r)}}\rightarrow
L^{p,\alpha}_{r}(u_{(r)}^{\ast}TM^+\otimes \wedge^{0,1})$$ such that $(\kappa_{v},\zeta)=((\kappa,h) + Q_{b_{(r)}}\circ f_{(r)}(\kappa,h)).$
\v
On the other hand,
suppose that  $(\kappa _{v},\zeta)$  satisfies
\begin{equation}\label{normal_implicit}
 F_{(r)}(\kappa_{v},\zeta)=0,\;\;\;\;\;\;\;
\|(\kappa_{v},\zeta)\| \leq \frac{ \hbar_3}{8(1+\mathcal{C}_{2}  \|D\mathcal{S}_{b_{(r)}}\|)C_{2}\mathcal{C}_{2}}.
\end{equation}
We write $(\kappa _{v},\zeta)=(\kappa,h)+Q_{b_{(r)}}\circ\eta$ where $(\kappa,h)\in ker D \mathcal{S}_{b_{(r)}}$ and $\eta\in L^{p,\alpha}_{r}(u_{(r)}^{\ast}TM^+\otimes \wedge^{0,1}).$
Since $(\kappa,h)=(I-Q_{b_{(r)}}\circ D\mathcal{S}_{b_{(r)}})(\kappa_{v},\zeta),$
 we have $\|(\kappa,h)\|\leq \frac{\hbar_3}{8C_{2}\mathcal{C}_{2}}$, then $\eta= f_{(r)}( \kappa,h)$, i.e.,
\begin{equation}\label{expression_s}
(\kappa_{v},\zeta)=   (\kappa,h)+ Q_{b_{(r)}}\circ f_{(r)}( \kappa,h)).
\end{equation}
Hence the zero set of $F_{(r)}$ is locally the form
$((\kappa,h), Q_{b_{(r)}}\circ f_{(r)}(\kappa,h))$,\;i.e
\begin{equation}\label{j_holomorphic_curve}
F_{(r)}((\kappa,h) + Q_{b_{(r)}}\circ f_{(r)}(\kappa,h))=0\end{equation}
where $(\kappa,h) \in Ker\;D{\mathcal S}_{b_{(r)}}$.
\v
For fixed $(r)$ denote
$$U_{  b_{(r)} }=\left\{(\kappa,h)\in {K}_{b_{(r)}} \times \mathcal W^{1,p,\alpha,r}(\Sigma_{(r)};u_{(r)}
^{\ast}TM^+)\mid \bar{\partial}_{J}exp_{u_{(r)}}h
+ \kappa=0\right\},$$
$$(\kappa_{(r)},\zeta_{(r)}):=I_{r}(\kappa,\zeta) + Q_{b_{(r)}}\circ f_{(r)}\circ I_{r}(\kappa,\zeta)).$$
Define $\pi_{1}$, $\pi_{2}$ by
$$\pi_{1}(\kappa,h)=\kappa,\;\;\;\;\mbox{and}\;\;\;\;  \pi_{2}(\kappa,h)=h,\;\;\mbox{ for any }(\kappa,h).$$
Since $I_r: ker D{\mathcal S}_{b}\rightarrow KerD{\mathcal S}_{b_{(r)}}$ is a isomorphism, we have proved the following

\begin{lemma}\label{coordinates_on_stratum_nodal}
 There is a neighborhood $O$ of $0$ in $ker D{\mathcal S}_{b}$
and a neighborhood $O_j$ of $(j_1,j_2)$ and $R_0 > 0$ such that
$$glu_{(r)}: O_j\times O \times \mathbb Z_k \rightarrow
U_{b_{(r)}}$$
defined by
$$glu_{(r)}(j_1,j_2,\kappa,\zeta,n)=(\pi_{1}(\kappa_{(r)},\zeta_{(r)}),\exp_{u_{(r)}}(\pi_{2}(\kappa_{(r)},\zeta_{(r)})))$$
where
for $r>R_0$ is a family of orientation preserving local diffeomorphisms.
\end{lemma}
\vskip 0.1in
\noindent
We may choose $((j_1,j_2),r, \tau,\kappa, \zeta)$ as a local coordinate system around $b$ in $U_{ I,\epsilon}$.
We write $f_{(r)}\circ I_r$ as a function $f(j_1,j_2,r, \tau,\kappa, \zeta)$.
As $I_r$ is a smooth map we can see that $f(j_1,j_2,r, \tau,\kappa, \zeta)$ is smooth (see p.131 in \cite{SL}).

\v
\v

\subsection{Surjectivity  and  injectivity  }

\v\n
\begin{prop}\label{surjective_injective_gluing_maps}
The gluing map in Lemma \ref{coordinates_on_stratum_nodal} is surjective and injective in the sense of the Gromov-Uhlenbeck convergence.
\end{prop}
\v\n
{\bf Proof. } The injectivity follows immediately from the implicit function theorem.
We prove the surjectivity. Let $\Gamma=((\kappa_{v},v),\Sigma_{(r)})\in U_{  b_{(r)} }$ and $\zeta \in W^{1,p,\alpha}_{r}$ such that $\exp_{u_{(r)}}(\zeta)=v.$ Suppose that $\|\kappa_{v}\|_{p,\alpha}+|\zeta|_{C^{1}(\mathfrak{D}_{i} (R_{0}))}\leq  \epsilon .$     By the same argument of Lemma \ref{approximation_estimates} we have
 $(\kappa,\zeta)$ satisfies \eqref{normal_implicit}.
Then   there exists a unique  $(\kappa,h)\in Ker\;D\mathcal S_{b}$ such that
$$
(\kappa_{v},\zeta)=I_{r}(\kappa,h) + Q_{b_{(r)}}\circ f_{(r)}\circ I_{r}(\kappa,h)).
$$
The surjectivity follows.
\v

\v\n

\v

\section{Estimates of differentiations for gluing parameters}

In the section, let $\alpha<<\mathfrak{c}$ be a constant. We estimate the differentiations for gluing parameters. This is another key point of this paper.

\subsection{Linear analysis on weight sobolov spaces}

Let $\Sigma=\mathbb R\times S^{1},$ $E= \Sigma\times \mathbb R^{2n}.$ Let
$\pi:E\longrightarrow\Sigma$ be the trivial vector bundle.
Consider the Cauchy-Riemannian operator $D$ on $E$ defined by
 $$D =\frac{\partial}{\partial
s}+J_0\frac{\partial} {\partial t}+S  = \bar{\partial}_{J_{0}} +
S ,$$
where $S(s,t) =(S_{ij}(s,t))_{2n\times 2n}$ is a matrix such that for any $k>0,$
\begin{equation}\label{expontent_decay_of_s}
 \sum_{i+j=k} \left|\frac{\p^{k}S}{\partial^ i s \partial ^j t}\right|\leq
C_{k}e^{-\fc s} .
\end{equation}
for some constant $C_{k}>0$. Therefore, the operator $H_s=J_0\frac{d}{dt}+S$ converges
to $H_{\infty}=J_0\frac{d}{dt}$ as $s\to \infty$.

We choose a small weight $\alpha$ for each
end such that $H_{\infty}- \alpha = J_0\frac{d}{dt}- \alpha $ is invertible. Let $W$ be the function in section \S\ref{weight_norm}.
  We can define the space $W^{1,p,\alpha}(\Sigma;E)$ and
$L^{p,\alpha}(E\otimes \wedge^{0,1})$. Then
the operator $D:   W^{1,p,\alpha}(\Sigma;E)\rightarrow L^{p,\alpha}(E\otimes \wedge^{0,1})$
is a Fredholm operator so long as $\alpha$ does not lie in the spectrum
of the operator $H_{i \infty}$ for all $i=1,\cdots,\nu$.
\v
By the observation of Donaldson we know that  the multiplication by $W$ gives isometries from $ W^{1,2,\alpha} $ to $W^{1,2} $ and $L^{2,\alpha} $ to $L^2,$
and
\begin{equation}
\label{weight_norm_relation}
\|f\|_{ {2,\alpha}}=\|Wf\|_{L^{2}},\;\;\;C\inv\|Wf\|_{W^{1,2}}\leq \|f\|_{{1,2,\alpha}}\leq C\|Wf\|_{W^{1,2}}.
\end{equation}
Consider  the operator
 $$D:W^{1,2,\alpha}\rightarrow
 L^{2,\alpha}. $$
 It is equivalent to the map (see \cite{D},P59)
$$D_{W}:=D-\frac{\partial}{\partial s}\log W= WD W^{-1}: W^{1,2}\rightarrow L^2.$$
Let $h\in W^{1,2,\alpha}$ be the solution of the equation
$D h =\eta,$ where $\eta\in L^{2,\alpha}.$ Obviously, $$D_{W}(Wh)=WD (W^{-1}W h )=W\eta.$$
Let $\rho =W\eta,\;\;f=Wh .$ Then $f$ satisfies  the equation $D_{W}f=\rho.$ 

When $s\geq R_{0}$ we consider  the operator $D-\alpha: W^{1,2}\rightarrow L^2.$
 We write $D:= D-\alpha=\frac{\p}{\p s} + L + S$, where $L:=J_{0}\frac{\partial }{\partial t} -\alpha$ is a invertible elliptic operator on $ S^1$ when $0<\alpha<1.$

The space $L^2( S^1)$ can be decomposed two infinite spaces $H^+$ and $H^-$(see P33 in \cite{D}), where
 $$H^+=\left\{\sum a_{n} e^{i n t}|n\geq 0\right\},\;\;\; H^-=\left\{\sum a_{n} e^{i n t}|n<0\right\}.$$
 Since $L$ is an invertible operator, there exists a complete eigen-function space decomposition such that  $L\phi_{\lambda} =\lambda \phi_{\lambda}$ where $|\lambda|>\delta$ for some $\delta>0.$ Obviously,
 $\lambda\in \mathbb Z-\alpha.$ By choosing $0<\alpha<1/2$   we have
 $\alpha=\min |\lambda|.$

Using the same method of Donaldson in \cite{D} we can prove the following lemmas (for details please see our next paper \cite{LS}.)

For $A> 1$ we consider the finite tube $ (-A, A )\times{S}^1$ and denote
$$
B^{-}_{A}=(-A,-A+1)\times S^1,\;\;\; B^{+}_{A}= (A-1,A)\times {S}^1.
$$

\begin{lemma}
Let $h$ be a solution of $D h = 0$ over the finite tube $(-A,A)\times S^1$. Then  for any $0<\alpha<\min\{\fc,\frac{1}{2}\}$, there exists a constant $C_{1}>0$ such that for any  $0<|s_{o}|+1<A$
\begin{equation} \label{tube_estimate_holomorphic_curve}
\|h\mid_{(-|s_{o}|+1,|s_{o}|-1)\times S^1}\|_{1,p,\alpha} \leq {C_{1}e^{- \alpha(A-|s_{o}|)}} \left( \|h|_{B^{+}_A}\|_{2,\alpha}+\|h|_{B^{-}_A}\|_{2,\alpha}\right).
\end{equation}
\end{lemma}

Using this lemma we can prove the following corollary.
\begin{corollary}\label{tube_estimate_c}
 Suppose that $D$ is surjective. Denote by $Q$ a bounded right inverse of $D$.  Let $h=Q\eta$ be a solution of $Dh = \eta$ over the finite tube $(0,\tfrac{2lr}{kT})\times S^1$. Then  for any $0<\alpha<\min\{\fc,\frac{1}{2}\}$, there exists a constant $C_{2}>0$ such that for any $0<A<\frac{lr}{2kT},$
\begin{equation} \label{tube_estimate_curve}
\|h\mid_{(\frac{lr}{2kT},\frac{3lr}{2kT})\times\mathbb{S}^1}\|_{1,2,\alpha} \leq C_{2}\left(e^{- \alpha A} \|\eta \|_{2,\alpha}+ \|\eta\mid_{\frac{lr}{2kT}-A\leq s \leq \frac{3lr}{2kT}+A}  \|_{2,\alpha}\right).
\end{equation}
\end{corollary}
{\bf Proof. } Denote by $\hat{\eta}$ the restriction of $\eta$
 to the part $ {\frac{lr}{2kT}-A\leq s \leq \frac{3lr}{2kT}+A} .$ Let $\hat{h}=Q\hat{\eta}.$
 Then $D(h-\hat{h})=0$ in $ {\frac{lr}{2kT}-A\leq s \leq \frac{3lr}{2kT}+A} .$ By  \eqref{tube_estimate_holomorphic_curve} and $$\|\hat{\eta}\|_{1,p,\alpha}=\|\eta\mid_{\frac{lr}{2kT}-A\leq s \leq \frac{3lr}{2kT}+A}  \|_{2,\alpha}\;\;\;\;\;\;\|\hat{h}\|_{1,p,\alpha}=\|Q\hat{\eta}\|_{1,p,\alpha}\leq C\|\hat{\eta}\|_{1,p,\alpha}$$ we have
 \begin{align*}
\|h\mid_{(\frac{lr}{2kT},\frac{3lr}{2kT})\times\mathbb{S}^1}\|_{1,2,\alpha} &\leq \|h-\hat{h}\mid_{(\frac{lr}{2kT},\frac{3lr}{2kT})\times\mathbb{S}^1}\|_{1,2,\alpha}+ \|\hat{h}\mid_{(\frac{lr}{2kT},\frac{3lr}{2kT})\times\mathbb{S}^1}\|_{1,2,\alpha} \\
&\leq C_{1}\left( {e^{- \alpha A}} \|h-\hat{h}\mid_{\frac{lr}{2kT}-A\leq s \leq \frac{3lr}{2kT}+A}  \|_{1,2,\alpha}+ \|\eta\mid_{\frac{lr}{2kT}-A\leq s \leq \frac{3lr}{2kT}+A}  \|_{2,\alpha}\right).
 \end{align*}
Then the corollary follows from $\|h-\hat{h}\|_{1,2,\alpha}=\|Q(\eta-\hat \eta)\|_{1,2,\alpha}\leq C\|\eta\|_{2,\alpha}$.   $\Box$

\subsection{Estimates of $\frac{\partial Q_{b_{(r)}}}{\partial r}$}
\begin{lemma} \label{exponential_right_inverse_lemma}   Let  $D{\mathcal
S}_{b} :K_{b}\times \mathcal{W}^{1,p,\alpha}\rightarrow L^{p,\alpha}$ be a Fredholm operator.
 Suppose that $D{\mathcal S}_{b}|_{K_{b}\times  {W}^{1,p,\alpha}}:K_{b}\times  {W}^{1,p,\alpha}\rightarrow L^{p,\alpha}$ is surjective.
Denote by $Q_b:L^{p,\alpha}\rightarrow K_{b}\times W^{1,p,\alpha}$ a right inverse of $D{\mathcal
S}_{b} .$ Let $Q_{b_{(r)}}$ be an  right inverse of D${\mathcal
S}_{b_{(r)}}$ defined in \eqref{express_right_inverse}.  Then there exists a constant $\mathcal C_{3}>0,$ independent of $r,$ such that
$\left\|\frac{\partial}{\partial r}Q_{b_{(r)}}\right\| \leq \mathcal C_{3}$
and  for any $0<s_{0}\leq \frac{lr}{4kT}$
  \begin{equation}\label{partial_right_estimate}\left\|\frac{\partial}{\partial r}\left(Q_{b_{(r)}}\eta_{r}\right)
\right\|
\leq \mathcal C_{3}\left[\|\frac{\p}{\p r}  \eta_{r}\|_{p,\alpha,r}+\frac{1}{r}\|\eta_{r}\mid_{\frac{lr}{4kT}\leq |s_{i}| \leq\frac{7lr}{4kT} }\|_{p,\alpha,r} +e^{-\alpha\frac{lr}{8kT}}\|\eta_{r}\|_{2,\alpha,r}\right],
 \end{equation}
for any $\eta_{r} \in L_{r}^{2,\alpha}$.
  \end{lemma}
{\bf Proof. }  Given $\eta_{r}\in L_{r}^{2,\alpha}$,
denote $\tilde \eta_{r}= (D{\mathcal S}_{b_{(r)}}\circ Q_{b_{(r)}}^ {\prime})^{-1}\eta_{r},$
 we have a pair $(\tilde \eta_1,\tilde \eta_2)$, where
 $\tilde \eta_1=\beta_{1}\tilde \eta_{r},\;  \tilde \eta_2=\beta_{2}\tilde \eta_{r}.$ Set
$$Q_b(\tilde \eta_1,\tilde \eta_2)=(\tilde{\kappa}_b,\tilde h_1,\tilde h_2),\;\; \tilde h_{(r)}=\tilde h_1\beta_1+ \tilde h_2\beta_{2}.$$
Then $Q_{b_{(r)}}\eta_{r}=Q_{b_{(r)}}^ {\prime}\tilde \eta_{r}=(\tilde{\kappa}_b,\tilde h_{(r)}).$ So
$$  \frac{\p}{\p r}\left(Q_{b_{(r)}}\eta_{r}\right)\mid_{|s_i|\leq \frac{lr}{2kT}}=\left(\frac{\p \tilde \kappa_{b}}{\p r},  \sum\frac{\p \tilde h_{i}}{\p r}\beta_{i} + \sum \tilde{h}_{i}\frac{\p \beta_{i}}{\p r}\right)\mid_{|s_i|\leq \frac{lr}{2kT}}=Q_b\left(\frac{\p}{\p r}\tilde \eta_1,\frac{\p}{\p r}\tilde \eta_2\right)\mid_{|s_i|\leq \frac{lr}{2kT}},$$
where we used $\beta\mid_{|s_i|\leq \frac{lr}{2kT}}=1.$
As $\frac{\p }{\p r}\tilde\eta_{i}=\beta_{i}\frac{\p }{\p r}\tilde\eta_{r}+\frac{\p \beta_{i}}{\p r}\tilde\eta_{r},$ we have
 \begin{equation}\label{app_right_inverse}
\left\|(\frac{\p \kappa}{\p r},\frac{\p\tilde h_{1}}{\p r},\frac{\p \tilde h_{2}}{\p r})\right\| = \left\|  Q_b\left(\frac{\p}{\p r}\tilde \eta_1,\frac{\p}{\p r}\tilde \eta_2\right) \right\| \leq C_1\left\|\frac{\p \tilde \eta_{r}}{\p r}\right\|_{p,\alpha,r}+\frac{C_1}{r} \left\| \tilde{\eta}_{r}\mid_{\frac{lr}{2kT}\leq |s_i|\leq  \frac{3lr}{2kT}} \right\|_{p,\alpha,r}.
 \end{equation}
By the same calculation as \eqref{approximate_difference} we have
\begin{eqnarray}\label{app_diff}
&&\eta_{r} -\tilde \eta_{r}= D \mathcal{S}_{b_{(r)}}\circ Q'_{b_{(r)}}\tilde{\eta}_{r}-(\beta_1^2+\beta_2^2)\tilde{\eta}_{r} = \kappa_{b_{(r)}} + D_{u_{(r)}}\tilde{h}_{(r)}-(\beta_1^2+\beta_2^2)\tilde{\eta}_{r} \nonumber\\
&=&\sum(\bar{\partial}\beta_{i}) \tilde{h}_{i} +\sum \beta_{i}(S_{u_{(r)}}-S_{u_{i}})\tilde{h}_{i}.
\end{eqnarray}
By the exponential decay $S_{u_{(r)}},S_{u_{i}} $ and
\begin{equation}\label{app_diff_1}
\bar{\p}\beta|_{|s_{i}|\leq \frac{lr}{2kT}}=0,\;|\bar{\p}\beta|\leq \frac{C}{r},S_{u_{(r)}}|_{|s_{i}|\leq \frac{lr}{2kT}}=S_{u_{i}} |_{|s_{i}|\leq \frac{lr}{2kT}},
\end{equation}   we conclude from \eqref{app_diff} that
\begin{equation}\label{app_est_1}
\|\tilde \eta_{r}\mid_{\frac{lr}{4kT}\leq |s_i|\leq  \frac{7lr}{4kT}} \|_{2,\alpha,r} \leq \left\|  \eta_{r}\mid_{\frac{lr}{4kT}\leq |s_i|\leq  \frac{7lr}{4kT}} \right\|_{2,\alpha,r} +\frac{C_{2}}{r}\|\tilde{h}_{i}\mid_{\frac{lr}{2kT}\leq |s_i|\leq  \frac{3lr}{2kT}}\|_{1,2,\alpha,r}.
\end{equation}
Applying Corollary \ref{tube_estimate_c} to \eqref{app_est_1}  we have
\begin{align}
\left(1-\frac{C_{3}}{r}\right)\|\tilde \eta_{r}\mid_{\frac{lr}{4kT}\leq |s_i|\leq  \frac{7lr}{4kT}} \|_{2,\alpha,r}  \leq \left\|  \eta_{r}\mid_{\frac{lr}{4kT}\leq |s_i|\leq  \frac{7lr}{4kT}} \right\|_{2,\alpha,r} +e^{-\alpha\frac{lr}{4kT}}\|\eta_{r}\|_{2,\alpha,r},\label{app_right_inverse_1}
\end{align}
where we used   $\|\tilde \eta_{r}\|_{p,\alpha,r}\leq 2\|  \eta_{r}\|_{p,\alpha,r}$.
Taking the derivative of \eqref{app_diff},  by the exponential decay of $S_{u_{(r)}},S_{u_{i}},$ and
$ |\frac{\p}{\p r}\bar{\p}\beta|\leq \frac{C}{r^2}, $ using \eqref{app_diff_1}, \eqref{tube_estimate_curve} and \eqref{app_diff},
we have
\begin{align}\label{app_right_inverse_0}
&\|\frac{\p}{\p r}\tilde \eta_{r}\|_{p,\alpha,r}\leq \|\frac{\p}{\p r}  \eta_{r}\|_{p,\alpha,r}+ \frac{C_{4}}{r^2}\|\tilde h_{i}\mid_{\frac{lr}{2kT}\leq |s_{i}| \leq\frac{3lr}{2kT} }\|_{1,p,\alpha}+\frac{C_{4}}{r}\left\|\frac{\p \tilde h_{i}}{\p r}\mid_{\frac{lr}{2kT}\leq |s_{i}| \leq\frac{3lr}{2kT} }\right\|_{1,p,\alpha}\nonumber\\
&\leq \|\frac{\p}{\p r}  \eta_{r}\|_{p,\alpha,r}+\frac{C_{5}}{r^2 }\|\tilde\eta_{r}\mid_{\frac{lr}{4kT}\leq |s_{i}| \leq\frac{7lr}{4kT} }\|_{1,p,\alpha} +e^{-\alpha\frac{lr}{4kT}}\|\eta_{r}\|_{2,\alpha}+\frac{C_{5}}{r}\left\|\frac{\p \tilde \eta_{r}}{\p r}\right\|_{p,\alpha,r} .
\end{align}
Following from \eqref{app_right_inverse_0} and \eqref{app_right_inverse_1} we obtain
\begin{align}
 &C_{1}\|\frac{\p}{\p r}\tilde \eta_{r}\|_{p,\alpha,r}+\frac{C_{1}}{r}\|\tilde \eta_{r}\mid_{\frac{lr}{4kT}\leq |s_i|\leq  \frac{7lr}{4kT}} \|_{2,\alpha,r} \nonumber\\
 & \leq C_{6}\left[\|\frac{\p}{\p r}  \eta_{r}\|_{p,\alpha,r}+\frac{1}{r}\left\|  \eta_{r}\mid_{\frac{lr}{4kT}\leq |s_i|\leq  \frac{7lr}{4kT}} \right\|_{2,\alpha,r} +e^{-\alpha\frac{lr}{4kT}}\|\eta_{r}\|_{2,\alpha}\right].\label{app_est_right_0}
\end{align}
  \eqref{app_right_inverse} and \eqref{app_est_right_0} gives us
\begin{equation}\label{partial_right_estimate_1}\left\|\frac{\partial}{\partial r}\left(Q_{b_{(r)}}\eta_{r}\right)
\mid_{|s_i|\leq s_{0}}\right\|_{1,2,\alpha,r}
\leq \|\frac{\p}{\p r}  \eta_{r}\|_{p,\alpha,r}+\frac{C_{7}}{r}\|\eta_{r}\mid_{\frac{lr}{4kT}\leq |s_{i}| \leq\frac{7lr}{4kT} }\|_{1,p,\alpha} +e^{-\alpha\frac{lr}{4kT}}\|\eta_{r}\|_{2,\alpha}.
 \end{equation}
Using   \eqref{tube_estimate_curve} and \eqref{app_right_inverse_1}, we have
\begin{align}
\|\sum \tilde{h}_{i}\tfrac{\p \beta_{i}}{\p r}\|_{1,p,\alpha}  &\leq \frac{C_{8}}{r}\left({e^{- \alpha\frac{lr}{8kT}}} \|\tilde{\eta} \|_{1,2,\alpha}+ \|\tilde{\eta}\mid_{\frac{3lr}{8kT}\leq |s_{i}| \leq\frac{13lr}{8kT} } \|_{2,\alpha}\right)\nonumber\\ \label{app_right_inverse_3}
& \leq \frac{C_{9}}{r}\left({e^{- \alpha\frac{lr}{8kT}}} \|\eta \|_{1,2,\alpha}+ \|\eta\mid_{\frac{lr}{4kT}\leq |s_{i}| \leq\frac{7lr}{4kT} } \|_{2,\alpha}\right)
\end{align}
Note that $\frac{\p}{\p r}\left(Q_{b_{(r)}}\eta_{r}\right) =\left(\frac{\p \tilde \kappa_{b}}{\p r},  \sum\frac{\p \tilde h_{i}}{\p r}\beta_{i} + \sum \tilde{h}_{i}\frac{\p \beta_{i}}{\p r}\right) .$ Then
\eqref{partial_right_estimate} follows from \eqref{app_right_inverse}, \eqref{app_est_right_0} and \eqref{app_right_inverse_3}.
  $\Box$

\subsection{Estimates of $\frac{\partial I_{(r)}}{\partial r}$}

\begin{lemma}Let $I_r:Ker D\mathcal {S}_{b}\longrightarrow Ker D\mathcal {S}_{b_{(r)}}$ be an isomorphisms defined in \eqref{definition_I}.
Then there exists a constant $C>0,$ independent of $r,$ such that
 \begin{equation}
 \label{partial_I_estimate} \left\|\frac{\partial}{\partial r}\left(I_{r}(\kappa,h,h_{0})\right)
 \right\|_{1,2,\alpha,r} \leq
\frac{ C}{r}\|(\kappa,h)\mid_ {\frac{lr}{2kT}\leq |s_i|\leq \frac{3lr}{2kT}}\|_{1,2,\alpha}\leq \frac{ C}{r}e^{-(\fc-\alpha)\frac{lr}{2kT}}\|(\kappa,h)\|,
 \end{equation}
 for any $(\kappa,h,h_{0})\in Ker D\mathcal {S}_{b}$.
\end{lemma}
{\bf Proof. }
By definition $\kappa$ is independent of $r$. Then $\frac{\p}{\p r}h_{(r)}=\sum \frac{\p}{\p r}\beta_{i}h_{i}$ and
$$\frac{\partial }{\partial r}I_{r}( \kappa,h ,h_{0})=(0,\frac{\p h_{(r)}}{\p r})+\frac{\partial }{\partial r}\left(Q_{b_{(r)}}\right)D\mathcal{S}_{b_{(r)}}(\kappa,h_{(r)})+Q_{b_{(r)}}\frac{\partial }{\partial r}\left(D\mathcal{S}_{b_{(r)}}(\kappa,h_{(r)})\right).$$
A direct calculation gives us
\begin{equation}
D\mathcal{S}_{b_{(r)}}(\kappa,h_{(r)})=(\bar{\partial}\beta_{1}) h_{1} +\beta_{1}(S_{u_{(r)}}-S_{u_{1}})h_{1}+(\bar{\partial}\beta_{2}) h_{2} +\beta_{2}(S_{u_{(r)}}-S_{u_{2}})h_{2}.
\end{equation}
By \eqref{app_diff_1} we have  $supp \; D\mathcal{S}_{b_{(r)}}(\kappa,h_{(r)})\subset \{  {\frac{lr}{2kT}\leq |s_i|\leq \frac{3lr}{2kT}}\}.$  Since
$\left(|S_{u_{(r)}}|+|S_{u}|\right)|_{ {\frac{lr}{2kT}\leq |s_i|\leq \frac{3lr}{2kT}}}\leq Ce^{-\fc\frac{lr}{2kT}} , $
the lemma follows from
   $
 \|Q_{b_{(r)}}\|\leq C,\; \|\frac{\p}{\p r}Q_{b_{(r)}}\|\leq C.
 $
    $\Box$

\subsection{Estimates of $\frac{\partial }{\partial r}\left[I_r(\kappa,\zeta)+Q_{b(r)}f_{(r)}(I_r(\kappa,\zeta))\right]$}

First we prove the following lemma.
\begin{lemma}
Let $(\kappa,\zeta)\in ker D\mathcal{S}_{b}$ with $\|(\kappa,\zeta)\|\leq \hbar_{3}$ such that
\begin{align}
\label{exp_J_map}
(\kappa,h_{(r)})=I_{r}(\kappa,\zeta)+Q_{b_{(r)}}\circ f_{r}\circ I_{r}(\kappa,\zeta)\\
\bar{\p}_{J}\exp_{u_{(r)}}(h_{(r)})+\kappa=0.
\end{align}
Then there exists a positive constants $C$ such that
\begin{equation}\label{decay_exp_J_map}
\|f_{r}\circ I_{r}(\kappa,\zeta)|_{\frac{lr}{4kT}\leq |s_{i}|\leq \frac{7lr}{4kT}}\|_{p,\alpha,r}\leq Ce^{-(\fc'-\alpha)\frac{lr}{4kT}}\|(\kappa,\zeta)\|,\;\;\;\;\forall \;\frac{lr}{8kT}\geq R_{0}.
\end{equation}
\end{lemma}
\n{\bf Proof. } From  \eqref{exp_J_map} and $\kappa|_{|s_{i}|\geq R_0}=0$
we have
\begin{equation}\label{express_f_r}
f_{(r)}(I_r(\kappa,\zeta))=\left(\frac{\p}{\p {s}}+J_{0}\frac{\p}{\p {t}}\right) h_{(r)}+S_{u_{(r)}}h_{(r)},\;\;\;\;\;\;\;\bar{\p}_{J} \exp_{u_{(r)}}h_{(r)}=0,
\end{equation}
 for any  $|s_{i}|\geq R_{0}.$ By the assumption $u_{i}$ converges to $kT$-periodic orbit $x(kTt)$ as $z\rightarrow p.$ In the given  local pseudo-Dauboux coordinate $(a,\vartheta,{\bf w})$, we  denote $$  \exp_{u_{(r)}}(h_{(r)})=(\hat a_{(r)}(s,t),\hat \vartheta_{(r)}(s,t),\hat {\bf w}_{(r)}(s,t)).$$ 
By Theorem \ref{exponential_estimates_theorem}, \eqref{gluing_a_r} and \eqref{gluing_t_u_r} we have for any $ |s_{i}|  \geq R_{0}$
  \begin{align} \label{exp_J_map_3}
|\nabla( a_{(r)} -kTs_{i})|\leq Ce^{-\fc \frac{lr}{4kT}},\;\;\;\; &
|\nabla( \vartheta_{(r)}-kt_{i})|\leq Ce^{-\fc \frac{lr}{4kT}},\;\;\;\;\\
\label{exp_J_map_4}
&|\nabla({\bf w}_{(r)})|\leq Ce^{-\fc \frac{lr}{4kT}}.
\end{align}

Let $$  a_{(r)}^{\diamond}=  a_{(r)} - kTs_{i},\;\;  \vartheta_{(r)}^{\diamond}=  \vartheta_{(r)} - kt_{i},\;\;\;\;\hat a_{(r)}^{\diamond}=\hat a_{(r)} - kTs_{i},\;\;\hat \vartheta_{(r)}^{\diamond}=\hat \vartheta_{(r)} - kt_{i}.$$
By \eqref{exp_J_map_3} and \eqref{exp_J_map_4} we have
\begin{equation}\label{a_energy}
\int_{|s_{i}|\geq s_{o}}\|\nabla a^{\diamond}_{(r)}\|^2+\|\nabla \vartheta^{\diamond}_{(r)}\|^2+\|\nabla {\bf w}_{(r)}\|^2dsdt\leq Ce^{-\fc s_{o}},\;\;\;\;\;\;\; \forall  s_{o}\geq R_{0}.
\end{equation}
Since $\|(\kappa,h_{(r)})\|\leq C\hbar_{3}$ and $\kappa|_{|s_{i}|\geq R_{0}}=0$, we have
\begin{equation}
\|h_{(r)}|_{|s_{i}|\geq R_{0}}\|_{1,2,\alpha}\leq C\hbar_{3}.
\end{equation}
 By the elliptic regularity we can get
 $$|h_{(r)}|\leq C e^{-\alpha \min\{|s_1|,|s_{2}|\}},\;\;\;\forall |s_{i}|\geq R_{0} .$$
  Then $\hat {\bf w}_{(r)}(|s_{i}|\geq R_{0})$ is also in  the pseudo-Dauboux coordinates when $ \frac{lr}{8kT}>R_{0}$ and $\hbar_{3}$ small.
 By \eqref{a_energy}
 \begin{equation}\label{a_energy_1}
\int_{|s_{i}|\geq s_{o}}\|\nabla \hat a^{\diamond}_{(r)}\|^2+\|\nabla \hat \vartheta^{\diamond}_{(r)}\|^2+\|\nabla \hat {\bf w}_{(r)}\|^2dsdt\leq Ce^{-\alpha s_{o}},\;\;\;\; \;\;\;\forall s_{o}\geq R_{0}.
\end{equation}
Then by Theorem \ref{tube_exponential_decay} we have for any $\frac{lr}{4kT}\leq |s_{i}|\leq \frac{7lr}{4kT}$
\begin{align}\label{exp_J_map_1}
&|\nabla( \hat a_{(r)}-kTs_{i})|\leq Ce^{-\fc \frac{lr}{4kT}},\;\;\;\;
|\nabla( \hat \vartheta_{(r)}-kt_{i}) )|\leq Ce^{-\fc \frac{lr}{4kT}},\;\;\;\;\\
\label{exp_J_map_2}
&|\nabla(\hat{\bf w}_{(r)}) |\leq Ce^{-\fc \frac{lr}{4kT}}.
\end{align}
Since  $\nabla  \exp_{u_{(r)}}(h_{(r)})=P_{u_{(r)},\exp_{u_{(r)}}(h_{(r)})}(\nabla u_{(r)})+d\exp_{u_{(r)}}(\nabla h_{(r)}),$   \eqref{exp_J_map_3}, \eqref{exp_J_map_4}, \eqref{exp_J_map_1} and \eqref{exp_J_map_2} gives us
\begin{equation}
|d\exp_{u_{(r)}} \nabla h_{(r)}|_{\frac{lr}{4kT}\leq |s_{i}|\leq \frac{7lr}{4kT}}|\leq Ce^{-\fc \frac{lr}{4kT}}.
\end{equation}
Note that $d\exp_{u_{(r)}}^{-1}$ is uniform bounded as $\|h_{(r)}\|_{1,p,\alpha,r}$ small.
Then
\begin{equation}\label{a_nabla_h}
\|\nabla h_{(r)} |_{\frac{lr}{4kT}\leq |s_{i}|\leq \frac{7lr}{4kT}}\|_{p,\alpha}\leq Ce^{-(\fc-\alpha)\frac{lr}{4kT} }.\end{equation}
The lemma follows from \eqref{express_f_r}, \eqref{a_nabla_h} and the exponential decay of $S_{u_{(r)}}.$   $\Box$

\begin{lemma}\label{coordinate_decay} There exists a constant $C>0$ such that, \begin{equation}
\label{partial_r_kernel_global_estimate}
\left\|\frac{\partial }{\partial r}\left(I_r(\kappa,\zeta)+Q_{b_{(r)}}\circ f_{(r)} \circ I_{r}(\kappa,\zeta) \right) \right\|_{1,2,\alpha,r}\leq  {C}e^{-\alpha\frac{lr}{8kT} } \|(\kappa,\zeta\|_{1,2,\alpha}.
\end{equation}
  Restricting to the compact set $\{|s_i|\leq R_{0}\}$,  we have
\begin{equation}
\label{partial_r_kernel_local_estimate}
\left |\frac{\partial }{\partial r}\left(I_r(\kappa,\zeta)+Q_{b_{(r)}}\circ f_{(r)} \circ I_{r}(\kappa,\zeta) \right) \right | \leq  {C}e^{-\alpha\frac{lr}{8kT} } \|\zeta\|_{L^{\infty}} .
 \end{equation}
  \end{lemma}
{\bf Proof.} Without loss of generality, we will fix a complex structure $(j_1,j_2)$ to simplify notations. We write
\begin{equation} \label{kernel_express}
I_r(\kappa,\zeta)+Q_{b(r)}\circ f_{(r)}(I_r(\kappa,\zeta))=(\kappa_{r},h_{(r)})
\end{equation}
where $\kappa_{r}\in K_{b(r)}, h_{(r)}\in C^{\infty}(\Sigma; u^{\star}_{(r)}TM^{+}).$ We have
\begin{equation}
\label{local_kernel_express}
\bar{\partial }_{J} exp_{u_{(r)}}h_{(r)}+\kappa_{r}=0.
\end{equation}
To simplify the expressions of formulas we denote
$$\frac{\partial }{\partial r}\exp_{u_{(r)}}(h_{(r)})=G,\;\;\;\;\;\frac{\partial }{\partial r}\kappa_{r}=\Psi.$$
Taking differentiation of \eqref{local_kernel_express} with respect to $r$ we get
\begin{equation} \label{partial_kernel_express}
 D_{\exp_{u_{(r)}}(h_{(r)})} G+\Psi =0.
\end{equation}
 We may write \eqref{partial_kernel_express} as
\begin{equation}
\label{local_calculation_1}
P_{exp_{u_{(r)}}h_{(r)}}\left(D\mathcal{S}_{exp_{u_{(r)}}h_{(r)}}(\Psi,G) \right)-
 D\mathcal{S}_{b_{(r)}}\left(\frac{\partial}{\partial r}\left(\kappa_{r}, h_{(r)}\right)\right)
\end{equation}
$$ +D\mathcal{S}_{b_{(r)}}\left(\frac{\partial}{\partial r}\left(\kappa_{r}, h_{(r)}\right)\right)=0.$$
 We estimate the difference
 $$ P_{exp_{u_{(r)}}h_{(r)}}\left(D\mathcal{S}_{exp_{u_{(r)}}h_{(r)}}(\Psi, G)  \right)-
 D\mathcal{S}_{b_{(r)}}\left(\frac{\partial}{\partial r}\left(\kappa_{r}, h_{(r)}\right)\right).$$
 We choose $\|\eta\|_{1,p,\alpha,r}$ very small. From  the Implicity function Theorem we have
\begin{equation}
\label{local_calculation_2}
 \|f_{(r)}(I_r(\kappa,\zeta))\|_{p,\alpha,r}\leq C\|(\kappa,\zeta)\|_{1,p,\alpha}+\|\bar{\p}_{J}u_{(r)}\|
\end{equation}
 for some constant $C>0.$ For any small $\zeta\in ker D{\mathcal S}_{b},exp_{u_{i}}\zeta$ converges to a periodic orbit as
  $|s_i|\to \infty.$ It follows that $S_{exp_{u}\zeta}$ converges to zero exponentially. Since $(\kappa,h_{(r)})$
  satisfies \eqref{local_kernel_express} with $\kappa|_{|s_1|\geq \frac{lr}{2kT}}=0$, by elliptic estimate we conclude that in the part $ \frac{lr}{2kT}\leq |s_i|\leq \frac{3lr}{2kT} $
\begin{equation}
\label{local_calculation_3}
 |S_{exp_{u_{(r)}}h_{(r)}}|\leq C_0 e^{-\mathfrak{c} \frac{lr}{2kT}}, \;\;\;\;|S_{u_{(r)}}|\leq C_0 e^{-\mathfrak{c}  \frac{lr}{2kT}}
\end{equation}
 for $r$  big enough. Moreover, we may choose $r$ very large and $| (\kappa,\zeta)|$ very small such that
$$
\|Q_{b_{(r)}}\|\left\|P_{exp_{u_{(r)}}h_{(r)}}\left(D\mathcal{S}_{exp_{u_{(r)}}h_{(r)}}(\Psi, G) \right)-D\mathcal{S}_{b_{(r)}}
\left(\frac{\partial}{\partial r}\left(\kappa_{r}, h_{(r)}\right)\right)\right\|_{1,p,\alpha,r}
$$\begin{equation}
\label{local_calculation_4} \leq \epsilon\left\|\frac{\partial}{\partial r}\left(\kappa_{r}, h_{(r)}\right)\right\|_{1,p,\alpha,r} +Ce^{-(\fc-\alpha)\frac{lr}{2kT} } ,
 \end{equation}
 where $\epsilon>0$ is a constant such that $4\epsilon \mathcal{C}_{3}\leq 1$.
  Note that
\begin{equation}
\label{partial_r_kernel_estimate}
 \frac{\partial}{\partial r} (\kappa_{r}, h_{(r)})= \frac{\partial}{\partial r} I_{r}(\kappa,\zeta) + \frac{\partial}{\partial r}\left (Q_{b_{(r)}}\right)\circ f_{(r)}(I_{r}(\kappa,\zeta))+\left(Q_{b_{(r)}}\frac{\partial}{\partial r}  \circ f_{(r)}(I_{r}(\kappa,\zeta))\right).
\end{equation}
Then
 $$D\mathcal{S}_{b_{(r)}}\left(\frac{\partial }{\partial r}(\kappa_{r},h_{(r)})\right)=D\mathcal{S}_{b_{(r)}}\left(\frac{\partial }{\partial r}I_{r}(\kappa,\zeta)\right) +D\mathcal{S}_{b_{(r)}}\left(\frac{\partial}{\partial r} (Q_{b_{(r)}})f_{(r)}(I_{r}(\kappa,\zeta))\right)+\frac{\partial}{\partial r}(f_{(r)}(I_{r}(\kappa,\zeta))).$$
It follows together with \eqref{local_calculation_2},  \eqref{local_calculation_4} that
$$
\left\|\frac{\partial}{\partial r}(f_{(r)}(I_{r}(\kappa,\zeta)))\right\|_{p,\alpha,r}
$$
\begin{equation}\label{operator_partial_right_inverse_kernel}
\leq \epsilon\left\|\frac{\partial}{\partial r} (\kappa_{r},h_{(r)})\right\|_{1,p,\alpha,r} +\left\| D\mathcal{S}_{b_{(r)}}\left(\frac{\partial}{\partial r}I_{r}(\kappa,\zeta)\right)\right\|_{p,\alpha,r}+\| D\mathcal{S}_{b_{(r)}} \circ \frac{\partial}{\partial r} (Q_{b_{(r)}})\circ f_{(r)}(I_{r}(\kappa,\zeta))\|_{p,\alpha,r}.
\end{equation}
For any $(\kappa_{1},h_{1}),$  we have
$$\frac{\partial }{\partial r}\left( D\mathcal{S}_{b_{(r)}}\right)(\kappa_{1},h_{1})=
\frac{\partial }{\partial r}\left( D_{u_{(r)}}h_{1}+\kappa_{1}\right)=\frac{\partial S_{u_{(r)}}}{\partial r}h_{1}=\left(0,\frac{\partial S_{u_{(r)}}}{\partial r}\right)\cdot(\kappa_{1},h_{1}). $$
Then
\begin{equation}\label{exp_decay_DS}
\frac{\partial }{\partial r}\left( D\mathcal{S}_{b_{(r)}}\right)=\left(0,\frac{\partial S_{u_{(r)}}}{\partial r}\right)
\end{equation}
 As $I_{r}(\kappa,\zeta)\in ker D\mathcal{S}_{b_{(r)}}$ we have
$$D\mathcal{S}_{b_{(r)}}\circ I_{r}(\kappa,\zeta) =0. $$
Then
 $$D\mathcal{S}_{b_{(r)}}\left(\frac{\partial}{\partial r}I_{r}(\kappa,\zeta)\right)+\left(0,\frac{\partial S_{u_{(r)}}}{\partial r}\right)I_{r}(\kappa,\zeta)=0.$$
 Since $\frac{\partial S_{u_{(r)}}}{\partial r}$ supports in the part $ {\frac{lr}{2kT}\leq |s_i|\leq \frac{3lr}{2kT}}$, by the exponential decay of $S_{u_{i}}$ we get
\begin{equation}\label{operator_partial_I}
\left\| D\mathcal{S}_{b_{(r)}} \left( \frac{\partial}{\partial r}I_{r}(\kappa,\zeta)\right) \right\|_{p,\alpha,r}\leq C_1 e^{-(\mathfrak{c} -\alpha)\frac{lr}{2kT}}.\end{equation}
 Taking the differentiation of the equality
 $$D\mathcal{S}_{b_{(r)}}\circ Q_{b_{(r)}}=I$$
 we get
 $$D\mathcal{S}_{b_{(r)}}\left( \frac{\partial }{\partial r}(Q_{b_{(r)}})\right) =-\left(0,\frac{\partial S_{u_{(r)}}}{\partial r}\right) \cdot Q_{b_{(r)}}. $$
 Together with \eqref{local_calculation_2}
 and \eqref{local_calculation_3} we get \begin{equation}\label{operator_partial_right_inverse_I}
 \|Q_{b_{(r)}}\|\left\|D\mathcal{S}_{b_{(r)}}  \frac{\partial}{\partial r}(Q_{b_{(r)}})  f_{(r)}(I_{r}(\kappa,\zeta))  \right\|_{1,p,\alpha,r}\leq C_2 e^{-(\mathfrak{c} -\alpha)\frac{lr}{2kT}} \end{equation} and
\begin{equation}\label{estimate_operator_partial_I}
\left\| \frac{\partial}{\partial r}(f_{(r)}(I_{r}(\kappa,\zeta)))\right\|_{p,\alpha,r}\leq C_2 e^{-(\mathfrak{c} -\alpha)\frac{lr}{8kT} }  +\epsilon\left\| \frac{\partial}{\partial r}\left(\kappa_{r}, h_{(r)}\right)\right\|_{1,p,\alpha,r}.
\end{equation}
By \eqref{kernel_express} we have
\begin{align}\label{partial_kernel_estimate_1}
&\left\|\frac{\partial}{\partial r}\left(\kappa_{r}, h_{(r)}\right)\right\|_{1,p,\alpha,r} \leq   \left( \left\|\frac{\partial }{\partial r}I_{r}(\kappa,\zeta)\right\|_{1,p,\alpha,r}+\left\|\frac{\partial}{\partial r}\left(Q_{b_{(r)}} \circ f_{(r)}(I_{r}(\kappa,\zeta))  \right)\right\|_{1,p,\alpha,r}\right)\end{align}
Together with  \eqref{partial_right_estimate} we have
\begin{align}
&\left\|\frac{\partial}{\partial r}\left(\kappa_{r}, h_{(r)}\right)\right\|_{1,p,\alpha,r}   \nonumber\\
&\leq \mathcal C_3 \left(e^{ -\alpha\frac{lr}{8kT}}+\left\|\frac{\partial }{\partial r}I_{r}(\kappa,\zeta)\right\|+\left\|\frac{\partial}{\partial r}\left(  f_{(r)}(I_{r}(\kappa,\zeta))  \right)\right\|_{1,p,\alpha,r}
+\|f_{r}\circ I_{r}(\kappa,\zeta)\mid_{\frac{lr}{4kT}\leq |s_{i}| \leq\frac{7lr}{4kT} }\|_{p,\alpha,r}  \right)\nonumber\\\label{exp_ker_coord_est}
 &\leq C_{4}e^{-\alpha\frac{lr}{8kT}}+\mathcal C_3\epsilon\left\|\frac{\partial}{\partial r}\left(\kappa_{r}, h_{(r)}\right)\right\|_{1,p,\alpha,r} . \end{align}
 where we used   \eqref{partial_I_estimate}, \eqref{decay_exp_J_map} and \eqref{estimate_operator_partial_I} in last inequality.
Then the Lemma follows from $4\mathcal C_3\epsilon<1$ and \eqref{exp_ker_coord_est}.    $\Box$
\vskip
0.1in \noindent

\section{Contact invariants and Open string invariants}

\vskip 0.1in
\noindent
\begin{prop} $U_{I}$ has the property
\begin{itemize}
\item[{\bf(1)}.] Each strata of $U_{I}$ is a smooth manifold.
\item[{\bf(2)}.] If $U_{I,D'}\subset U_{I,D}$ is a lower stratum,
it is a submanifold of codimension
at least 2.
\end{itemize}
\end{prop}
{\bf Proof.} The proof of {\bf(1)} for the top strata is standard, we omit it here. The proof for the lower strata will be given in our next paper \cite{LS}.
It is well-known that an interior node corresponds to codimension 2 stratum. It suffices to consider those nodal corresponding to periodic orbits. Note that:
\begin{itemize}
\item[{\bf(1)}.] we mod the freedoms of choosing the origin in $\mathbb R$ and the origin of the periodic orbits;
\item[{\bf(2)}.] we choose the Li-Ruan's compactification in \cite{LR}, that is, we firstly let the Riemann surfaces degenerate in Delingne-Mumford space and then let $M^+$ degenerate compatibly. At any node, the Riemann surface degenerates independently with two parameters, which compatible with those freedoms of choosing the origins ( see section \S\ref{compact_theorem} for degeneration and section \S\ref{gluing_pregluing} for gluing);
\end{itemize}
Then both blowups at interior and at infinity
lead boundaries of codimension 2 or more. $\Box$

\vskip 0.1in
\noindent
From section \S\ref{global_r} we have a finite dimensional virtual orbifold system $\{U_{I}, E_{I},\sigma_{I}\}$ and a finite dimensional virtual orbifold $\{U_{I}\}$  indexed by a partially ordered
set $(\mathcal I =2^{\{1, 2, \cdots, N\}}, \subset )$.
 By the same argument of Theorem \ref{compact_moduli_space} we conclude that $\{\Sinverse\}$ is a compact virtual manifold. As we have mod the $S^1$-action on periodic orbit, and $S^1$-action on the puncture points on Riemann surfaces, by the same method of \cite{MS} we can show that  $\{\Sinverse\}$ is oriented (see \cite{LS}).

 Let $\Lambda=\{\Lambda_I\}$ be a partition of unity and $\{\Theta_I\}$ be a virtual  Euler form of $\{\mathbb E_I\}$ such that
$\Lambda_I\Theta_I$ is compactly supported in $U_{I,\epsilon}$.

Recall that we have two natural maps
$$e_i: U_{I,\epsilon} \longrightarrow
M^{+} $$
$$ (u;\Sigma,{\bf y},{\bf p}) \longrightarrow
u(y_{i}) $$ for $i\leq m$ defined by evaluating at marked points
and $$e_j: U_{I,\epsilon}\longrightarrow Z_{\mathfrak{e}_j}$$
$$ (u;\Sigma,{\bf y},{\bf p}) \longrightarrow
u(p_{j}) $$
for $j>m$ defined by projecting to its
periodic orbit. The contact invariant can be
defined as
\begin{equation}\label{c_invariant}
\Psi^{(C)}_{(A,g,m+\nu,\mathbf{k},\mathfrak{e})}(\alpha_1,...,
\alpha_{m} ; \beta_{m + 1},..., \beta_{m+\nu})=\sum_I\int_{U_{I,\epsilon}}\prod_i e^*_i\alpha_i\wedge
\prod_j e^*_j\beta_j\wedge \Lambda_I \Theta_I.
\end{equation} for
$\alpha_i\in H^*(M^{+}, {\mathbb{R}})$ and $\beta_j\in H^*(Z_{\mathfrak{e}_j},
{\mathbb{R}})$ represented by differential form. Clearly, $\Psi^C=0$ if
$\sum \deg(\alpha_i)+ \sum \deg (\beta_i)\neq Ind^C$.
\v
Similarly, the open string invariant can be
defined as
\begin{equation}\label{l_invariant}
\Psi^{(L)}_{(A,g,m+\nu,\overrightarrow{\mu})}(\alpha_1,...,
\alpha_{m})=\sum_I\int_{U_{I,\epsilon}}\prod_i e^*_i\alpha_i\wedge \Lambda_I \Theta_I.
\end{equation} for
$\alpha_i\in H^*(M^{+}, {\mathbb{R}})$. Clearly, $\Psi^{(L)}=0$ if
$\sum \deg(\alpha_i)\neq Ind^L$.
\v
It is proved that these integrals are independent of the choices of $\Theta_I$ and the choices of the regularization ( see \cite{CT}). We must prove the convergence of the integrals \eqref{c_invariant} and \eqref{l_invariant} near each lower strata.
\v
Let $\alpha\in
H^*(M^{+}, {\mathbb{R}})$ and $\beta\in H^*(Z,{\mathbb{R}})$ represented by
differential form. We may write $$\prod_i
e^*_i\alpha_i\wedge\prod_j e^*_j\beta_j\wedge \Lambda_I \Theta_I =
ydr\wedge d\tau \wedge dj\wedge d\zeta,$$ where $d\zeta$ and $dj$ denote
the volume forms of $ker D{\mathcal S}_{b}$ and the space of complex
structures respectively and $y=y(r,\tau, j,\zeta)$ is a function. Then \eqref{partial_r_kernel_local_estimate} implies
that $|y|\leq C_1e^{-C_2r}$ for some constants $C_1>0,\;C_2>0$. Then the convergence of the integral
\eqref{c_invariant} follows. Similarly, we can prove the  convergence of the integral
\eqref{l_invariant}. \vskip 0.1in \noindent
Obviously, both $\Psi^{(C)}$ and $\Psi^{(L)}$ are generalizations of the relative GW invariants.
\v
One can easily show that
\begin{theorem}
\vskip 0.1in
\noindent
{\it (i). $\Psi^{(C)}_{(A,g,m+\nu,\mathbf{k},\mathfrak{e})}(\alpha_1,...,
\alpha_{m} ; \beta_{m + 1},..., \beta_{m+\nu})$ is well-defined, multi-linear and skew
symmetric.
\vskip 0.1in
\noindent
(ii). $\Psi^{(C)}_{(A,g,m+\nu,\mathbf{k},\mathfrak{e})}(\alpha_1,...,
\alpha_{m} ; \beta_{m + 1},..., \beta_{m+\nu})$ is independent of the choice of forms $\alpha_i, \beta_j$
representing the cohomology classes $[\beta_j], [\alpha_i]$,  and the choice
of virtual neighborhoods.
\vskip 0.1in
\noindent
(iii). $\Psi^{(C)}_{(A,g,m+\nu,\mathbf{k},\mathfrak{e})}(\alpha_1,...,
\alpha_{m} ; \beta_{m + 1},..., \beta_{m+\nu})$
is independent of the choice of $\widetilde{J}$ and $J $ over $M_{0}^{+}$.}
\end{theorem}

\begin{theorem}
\vskip 0.1in
\noindent
{\it (i). $\Psi^{(L)}_{(A,g,m+\nu,\overrightarrow{\mu})}(\alpha_1,...,
\alpha_{m})$ is well-defined, multi-linear and skew
symmetric.
\vskip 0.1in
\noindent
(ii). $\Psi^{(L)}_{(A,g,m+\nu,\overrightarrow{\mu})}(\alpha_1,...,
\alpha_{m})$ is independent of the choice of forms $\alpha_i$
representing the cohomology classes $[\alpha_i]$,  and the choice
of virtual neighborhoods.
\vskip 0.1in
\noindent
(iii). $\Psi^{(L)}_{(A,g,m+\nu,\overrightarrow{\mu})}(\alpha_1,...,
\alpha_{m})$
is independent of the choice of $\widetilde{J}$ and $J $ over $M_{0}^{+}$.}
\end{theorem}

\end{document}